\documentclass[a4paper,12pt]{article}
\usepackage{amsfonts,amsmath}

\def\comack{\mathsf{M}^c}

\def\per{\mathsf{perm}}
\def\fun{\mathsf{Fun}}
\def\iiota{\imath}
\def\jiota{\jmath}
\newcommand{\bin}[3]{\displaystyle{\binom{S_#1^#3}{S_#2^#3}}}
\def\To{\longrightarrow}

\makeatletter
\newbox\bk@bxb
\newbox\bk@bxa
\newif\if@bkcont
\newcount\bk@lcnt

\def\breakboxskip{2pt}
\def\breakboxparindent{1.8em}

\def\breakbox{\vskip\breakboxskip\relax
\setbox\bk@bxb\vbox\bgroup
\advance\linewidth -2\fboxrule
\hsize\linewidth\@parboxrestore
\parindent\breakboxparindent\relax}

\def\bk@split{%
\@tempdimb\ht\bk@bxb 
\advance\@tempdimb\dp\bk@bxb
\setbox\bk@bxa\vsplit\bk@bxb to\z@ 
\setbox\bk@bxa\vbox{\unvbox\bk@bxa}
\setbox\@tempboxa\vbox{\copy\bk@bxa\copy\bk@bxb}
\advance\@tempdimb-\ht\@tempboxa
\advance\@tempdimb-\dp\@tempboxa}

\def\bk@addfsepht{%
\setbox\bk@bxa\vbox{\vskip\fboxsep\box\bk@bxa}}

\def\bk@addskipht{%
\setbox\bk@bxa\vbox{\vskip\@tempdimb\box\bk@bxa}}

\def\bk@addfsepdp{%
\@tempdima\dp\bk@bxa
\advance\@tempdima\fboxsep
\dp\bk@bxa\@tempdima}

\def\bk@addskipdp{%
\@tempdima\dp\bk@bxa
\advance\@tempdima\@tempdimb
\dp\bk@bxa\@tempdima}

\def\bk@line{%
\hbox to \linewidth{%
\hskip-2\fboxsep\vrule \@width\fboxrule\hskip.5\fboxsep\vrule \@width\fboxrule\hskip1.5\fboxsep
\box\bk@bxa\hfil
}}%

\def\endbreakbox{\egroup
\ifhmode\par\fi{\noindent\bk@lcnt\@ne
\@bkconttrue\baselineskip\z@\lineskiplimit\z@
\lineskip\z@\vfuzz\maxdimen
\bk@split\bk@addfsepht\bk@addskipdp
\ifvoid\bk@bxb 
\def\bk@fstln{\bk@addfsepdp
\hskip-\parindent\vbox{\llap{\raisebox{-2ex}{\rule{1.5\fboxsep}{\fboxrule}\hskip.5\fboxsep}}\bk@line\llap{\rule{1.5\fboxsep}{\fboxrule}\hskip.5\fboxsep}}}

\else 
\def\bk@fstln{\vbox{\llap{\raisebox{-2ex}{\rule{1.5\fboxsep}{\fboxrule}\hskip.5\fboxsep}}\bk@line}\hfil%
\advance\bk@lcnt\@ne
\loop
\bk@split\bk@addskipdp\leavevmode
\ifvoid\bk@bxb 
\@bkcontfalse\bk@addfsepdp
\vtop{\bk@line\noindent\hskip-2\fboxsep{\rule{1.5\fboxsep}{\fboxrule}}}%

\else 
\bk@line
\fi
\hfil\advance\bk@lcnt\@ne
\if@bkcont\repeat}%
\fi
\leavevmode\bk@fstln\par}\vskip\breakboxskip\relax}

\def\smp{\smallskip\par}
\def\un{{\bf 1}}
\def\zero{\{0\}}
\def\pf{\noindent{\bf Proof~:}\ }
\def\findemo{~\leaders\hbox to 1em{\hss\  \hss}\hfill~\raisebox{.5ex}{\framebox[1ex]{}}\smp}
\def\spn{\bigskip\par\noindent}
\def\mpn{\medskip\par\noindent}
\def\smpn{\smallskip\par\noindent}
\def\normal{\mathop{\trianglelefteq}}

\def\smp{\smallskip\par}
\def\smpn{\smallskip\par\noindent}

\def\mpoint{\;\;.}
\def\mvirg{\;\;,}

\def\Res{{\rm Res}}
\def\Ind{{\rm Ind}}
\def\Inf{{\rm Inf}}

\def\Iso{{\rm Iso}}

\def\Hom{{\rm Hom}}
\def\End{{\rm End}}
\def\Ext{{\rm Ext}}
\def\Inf{{\rm Inf}}
\def\Im{{\rm Im}}

\def\Ker{{\rm Ker}}

\def\Id{{\rm Id}}

\def\Irr{{\rm Irr}}

\def\op{^{op}}

\def\dsp{\displaystyle}
\def\Z{\mathbb{Z}}
\def\N{\mathbb{N}}
\def\F{\mathbb{F}}

\newcommand{\dirsum}[1]{\mathop{\oplus}_{#1}\limits}
\newcommand{\romain}[1]{\uppercase\expandafter{\romannumeral #1}}

\newcommand{\flh}[2]{\mathop{\hbox to 12mm{\rightarrowfill}}_{\displaystyle #2}^{\displaystyle #1}\limits}
\newcommand{\sflh}[2]{\mathop{\hbox to 12mm{\rightarrowfill}}_{\scriptstyle #2}^{\scriptstyle #1}\limits}

\newcommand{\gMod}[1]{#1{\hbox{-}\mathsf{Mod}}}
\newcommand{\gmod}[1]{#1\hbox{-$\mathsf{mod}$}}

\newcommand{\sur}[1]{\,\overline{\! #1}}
\newcommand{\sumb}[2]{\mathop{\sum}_{{\scriptstyle #1}\atop {\scriptstyle #2}}\limits}
\newcommand{\oplusb}[2]{\mathop{\bigoplus}_{{\scriptstyle #1}\atop{\scriptstyle #2}}}
\newcommand{\dirsumb}[2]{\mathop{\oplus}_{{\scriptstyle #1}\atop{\scriptstyle #2}}\limits}
\newcommand{\oplusc}[3]{\mathop{\bigoplus}_{{\scriptstyle #1}\atop {{\scriptstyle #2}\atop {\scriptstyle #3}}}}
\newcommand{\dirsumc}[3]{\mathop{\oplus}_{{\scriptstyle #1}\atop {{\scriptstyle \rule{0ex}{1.5ex}#2}\atop {\scriptstyle \rule{0ex}{1.5ex}#3}}}}

\newcommand{\sumc}[3]{\sum_{{\scriptstyle #1}\atop {{\scriptstyle #2}\atop {\scriptstyle #3}}}}

\def\op{^{op}}

\newcommand{\carre}[8]{\begin{array}{ccc}
#1&\mathop{\hbox to 12mm{\rightarrowfill}}^{\displaystyle{#2}}\limits&#3\\
\llap{$\displaystyle{#4}$}\left\downarrow\vbox to 6mm{}\right. & & \left\downarrow\vbox to 6mm{}\right.\rlap{$\displaystyle{#5}$}\\
#6&\mathop{\hbox to 12mm{\rightarrowfill}}_{\displaystyle #7}\limits&#8\\
\end{array}}

\newcommand{\carrem}[8]{\begin{array}{ccc}
#1&\mathop{\hbox to 12mm{\rightarrowfill}}^{\displaystyle #2}\limits&#3\\
\llap{$\displaystyle #4$}\left\uparrow\vbox to 6mm{}\right. & & \left\uparrow\vbox to 6mm{}\right.\rlap{$\displaystyle #5$}\\
#6&\mathop{\hbox to 12mm{\rightarrowfill}}_{\displaystyle #7}\limits&#8\\
\end{array}}

\newenvironment{enonce}[1]{\pagebreak[2]\refstepcounter{subsection}\refstepcounter{prop}\smpn{{\bf \thesection.\arabic{prop}.\ \ #1~:}}\begin{it} }{\end{it}\smp}
\newenvironment{enonce*}[1]{\pagebreak[2]\smpn{#1~:}\begin{it} }{\end{it}\smp}
\newcommand{\result}[1]{\begin{enonce}{#1}}
\def\fresult{\end{enonce}}
\newcommand{\npar}{\smallskip\par\noindent\pagebreak[2]\refstepcounter{subsection}\refstepcounter{prop}{\bf \thesection.\arabic{prop}.\ \ }}
\newcommand{\masubsect}[1]{\medskip\par\noindent\pagebreak[3]\refstepcounter{subsection}\refstepcounter{prop}{\bf \thesection.\arabic{prop}.\ \ #1.\ }}

\newenvironment{mth}[1]{\begin{breakbox}\begin{enonce}{#1}}{\end{enonce}\end{breakbox}}
\newenvironment{mth*}[1]{\begin{breakbox}\begin{enonce*}{#1}}{\end{enonce*}\end{breakbox}}
\newenvironment{rem}[1]{\refstepcounter{subsection}\refstepcounter{prop} \mpn{{\bf \thesection.\arabic{prop}.}\ \ \bf#1\ :}}{\smp}

\def\dom{\backslash}
\makeatletter
\renewenvironment{enumerate}{\ifnum \@enumdepth >3 \@toodeep\else
      \advance\@enumdepth \@ne
      \edef\@enumctr{enum\romannumeral\the\@enumdepth}\list
      {\csname label\@enumctr\endcsname}{\setlength{\topsep}{1ex}\setlength{\itemsep}{0pt}\usecounter
        {\@enumctr}\def\makelabel##1{\hss\llap{##1}}}\fi}{\endlist}
\renewenvironment{itemize}{\ifnum \@itemdepth >3 \@toodeep\else \advance\@itemdepth \@ne
\edef\@itemitem{labelitem\romannumeral\the\@itemdepth}%
\list{\csname\@itemitem\endcsname}{\setlength{\topsep}{1ex}\setlength{\itemsep}{0pt}\def\makelabel##1{\hss\llap{##1}}}\fi}
{\endlist}
\makeatother
\makeatletter
\def\@sect#1#2#3#4#5#6[#7]#8{\ifnum #2>\c@secnumdepth
    \let\@svsec\@empty\else
    \refstepcounter{#1}\edef\@svsec{\csname the#1\endcsname .\hskip .5em}\fi
    \@tempskipa #5\relax
     \ifdim \@tempskipa>\z@
       \begingroup #6\relax
         \@hangfrom{\hskip #3\relax\@svsec}{\interlinepenalty \@M #8\par}%
       \endgroup
      \csname #1mark\endcsname{#7}\addcontentsline
        {toc}{#1}{\ifnum #2>\c@secnumdepth \else
                     \protect\numberline{\csname the#1\endcsname}\fi
                   #7}\else
       \def\@svsechd{#6\hskip #3\relax  
                  \@svsec #8\csname #1mark\endcsname
                     {#7}\addcontentsline
                          {toc}{#1}{\ifnum #2>\c@secnumdepth \else
                            \protect\numberline{\csname the#1\endcsname}\fi
                      #7}}\fi
    \@xsect{#5}}
\def\section{\@startsection {section}{1}{\z@}{-3.5ex plus-1ex minus
    -.2ex}{2.3ex plus.2ex}{\reset@font\Large\bf}}  

\makeatother
\renewenvironment{equation}{\refstepcounter{subsection}\refstepcounter{prop}$$}{\leqno{\bf (\theprop)}$$}

\def\mar[#1]{\ar@{-}[#1]|-{\object@{<}}}
\def\marb[#1]{\ar@{-}[#1]|{\object+{  }}}

\usepackage[all]{xy}
\begin{document}
\centerline{\Large\bf Complexity and cohomology}
\vspace{.2cm}
\centerline{\Large\bf of cohomological Mackey functors}
\vspace{.5cm}
\centerline{\bf Serge Bouc}
\vspace{1cm}
{\footnotesize{\bf Abstract~:} Let $k$ be a field of characteristic $p>0$. Call a finite group $G$ a {\em poco group over $k$} if any finitely generated cohomological Mackey functor for $G$ over $k$ has polynomial growth. The main result of this paper is that $G$ is a poco group over $k$ if and only if the Sylow $p$-subgroups of $G$ are cyclic, when $p>2$, or have sectional rank at most 2, when $p=2$.\par
A major step in the proof is the case where $G$ is an elementary abelian $p$-group. In particular, when $p=2$, all the extension groups between simple functors can be determined completely, using a presentation of the graded algebra of self extensions of the simple functor $S_\un^G$, by explicit generators and relations.\medskip\par
{\bf AMS Subject Classification~:} 16P90, 18G10, 18G15, 20J05.\par
{\bf Keywords~:} Cohomological, Mackey functor, complexity, growth.
}
\section{Introduction}\label{intro}
This paper addresses the question of the cohomology and rate of growth of cohomological Mackey functors for a finite group $G$ over a field $k$. The three main results are the following~:
\begin{mth}{Theorem} \label{main}Let $G$ be a finite group and $k$ be a field of positive characteristic $p$. Then the following conditions are equivalent~:
\begin{enumerate}
\item Every finitely generated cohomological Mackey functor for $G$ over $k$ has a projective resolution with polynomial growth.
\item Let $S$ be a Sylow $p$-subgroup of $G$. Then $S$ is cyclic if $p>2$, or $S$ has sectional 2-rank at most 2, if $p=2$.
\end{enumerate}
\end{mth}
The key argument in the proof of this theorem is a reduction to the case where $G$ is an elementary abelian $p$-group. The case of elementary abelian 2-groups can be described quite completely~: 
\begin{mth}{Theorem} \label{Poincare series}Let $G=(C_2)^m$ be an elementary abelian 2-group of rank~$m$, and $k$ be a field of characteristic $2$. Let $\comack_k(G)$ denote the category of cohomological Mackey functors for $G$ over $k$. Let $S_\un^G$ denote the simple cohomological Mackey functor defined by
$$\forall Q\leq G,\;\;S_\un^G(Q)=\left\{\begin{array}{cl}k&\hbox{if}\;Q=\un\\\zero&\hbox{otherwise}\end{array}\right.\mpoint$$
Then the algebra $\Ext_{\comack_k(G)}^*(S_\un^G,S_\un^G)$ is finitely generated by elements of degree~2, and its Poincar\'e series
$$P(t)=\sum_{j\in\N}\dim_k\Ext_{\comack_k(G)}^j(S_\un^G,S_\un^G)\;t^j$$
is equal to
$$P(t)=\frac{1}{(1-t^2)(1-3t^2)(1-7t^2)\ldots\big(1-(2^{m-1}-1)t^2\big)}\mpoint$$
\end{mth}
\begin{mth}{Theorem}\label{presentation} Let $G=(C_2)^m$ be an elementary abelian 2-group of rank~$m$, and $k$ be a field of characteristic $2$. Then the algebra $\mathcal{E}=\Ext^*_{\comack_k(G)}(S_\un^G,S_\un^G)$ admits the following presentation~:
\begin{itemize}
\item The generators $\gamma_x$, of degree 2, are indexed by the elements $x$ of $G-\{0\}$.
\item The relations are the following~:
\begin{enumerate}
\item Whenever $H$ is a subgroup of index 2 of $G$,
$$\sum_{x\notin H}\gamma_x=0\mpoint$$
\item For any distinct elements $x$ and $y$ of $G-\{0\}$,
$$[\gamma_x+\gamma_y,\gamma_{x+y}]=0\mvirg$$
where $[a,b]=ab+ba$ denotes the commutator of two elements $a$ and~$b$ in~$\mathcal{E}$.
\end{enumerate}
\end{itemize}
\end{mth}
\npar This paper is divided in two parts~: the first one focuses on complexity, and the second one on cohomology. Section~\ref{growth} of Part~\romain{1} quickly recalls the definitions and basic results on the rate of growth of a module for a finite dimensional algebra over a field. In Section~\ref{comack}, the categories of cohomological Mackey functors are introduced, from different points of view, which are equivalent thanks to Yoshida's theorem. Next some functors associated to bisets between categories of cohomological Mackey functors are defined, with nice adjunction properties. In Section~\ref{reduction}, it is shown how to reduce the question of complexity of cohomological Mackey functors for a finite group~$G$ over a field of characteristic $p$ to the same question for a Sylow $p$-subgroup of~$G$. Section~\ref{sketch} exposes a sketch of the proof of Theorem~\ref{main}. In Sections~\ref{simple functors},~\ref{some functors}, and ~\ref{extensions}, some simple cohomological functors and extensions between them are discussed. In Section~\ref{cyclic p-groups}, the case of cyclic $p$-groups is recalled from Samy Modeliar's thesis~(\cite{samymodeliar}). In Section~\ref{elemab}, the case of elementary abelian $p$-groups is settled, and Section~\ref{sectional 2-rank 2} handles the case of 2-groups with sectional 2 rank (at most) equal to~2.\par
The first section of Part~\romain{2} states further results on extensions of simple cohomological functors for elementary abelian $p$-groups. This leads to the proof of Theorem~\ref{presentation}, in Section~\ref{algebra}. In Section~\ref{algebra p odd}, a similar partial result is stated, for $p=3$, which is conjectured for any odd prime $p$. Finally, Section~\ref{more extensions} exposes some results on extensions of simple functors for an arbitrary finite $p$-group $G$, which show in particular how to reduce the computation of these extensions to the computation of self extensions for simple functors indexed by the trivial subgroup of some subquotients of $G$.\spn
{\bf Acknowledgments~:} I wish to thank the MSRI, where this work was completed during my stay there for the program on Representation Theory of Finite Groups and Related Topics, in spring 2008. I also thank Dave Benson for stimulating conversations about all this.
\vspace{1cm}\par
\centerline{\Large\bf \romain{1} - Complexity}
\section{Polynomial growth}\label{growth}
\begin{mth}{Definition} Let $A$ be a finite dimensional (unital) algebra over a field $k$. A finitely generated $A$-module $M$ is said to have {\em polynomial growth} if there exists a resolution
$$P_*: \cdots \To P_n\To P_{n-1}\cdots \To P_0\To M\To 0$$
of $M$ by projective $A$-modules, and constants $c$, $d$ and $e$ such that
$$\forall n\in\N,\;\;\dim_kP_n\leq c\, n^d+e\mpoint$$
The module $M$ is said to have {\em exponential growth} if for any projective resolution $P_*$ of $M$, there are constants $c$, $d$, and $e$, with $c>0$ and $d>1$, such that 
$$\forall n\in\N,\;\;\dim_kP_n\geq c\, d^n+e\mvirg$$
The module $M$ is said to have {\em intermediate growth} if $M$ has neither polynomial nor exponential growth.
\end{mth}
\begin{rem}{Remark} If $S$ is a generating set of $M$ as an $A$-module, then there is a surjective map of $A$-modules $A^{|S|}\To M$. In particular, the projective cover of $M$ has dimension at most $d\,{\dim_kM}$, where $d=\dim_kA$. By induction, this shows that there exists a projective resolution as above such that 
$$\dim_kP_n\leq (d-1)^{n-1}d\dim_kM\mpoint$$
In particular, this dimension is always bounded by some exponential function of $n$.
\end{rem}
\begin{mth}{Lemma} \label{simple enough}Let $A$ be a finite dimensional algebra over a field $k$, and~$M$ be a finitely generated $A$-module.
\begin{enumerate}
\item If 
$$\cdots \To P_n\To P_{n-1}\cdots \To P_0\To M\To 0$$
is a minimal projective resolution of $M$, then
$$P_n\cong \dirsum{S\in \Irr(A)}P_S^{\dim_k\Ext_A^n(M,S)/\dim_k\End_A(S)}\mvirg$$
where $\Irr(A)$ is a set or representatives of isomorphism classes of simple $A$-modules, and $P_S$ denotes a projective cover of $S$.
\item In particular $M$ has polynomial growth if and only if for any simple $A$-module $S$, there exists constants $c$, $d$ and $e$ such that 
$$\forall n\in \N,\;\;\dim_k\Ext_A^n(M,S)\leq c\,n^d+e\mpoint$$
\end{enumerate}
\end{mth}
\pf Assertion~1 follows by {\em d\'ecalage} from the case $n=0$, and from the fact that the largest semisimple quotient of $M$ is isomorphic to
$$\dirsum{S\in \Irr(A)}S^{\dim_k\Ext_A^0(M,S)/\dim_k\End_A(S)}\mpoint$$
Assertion~2 is a straightforward consequence of Assertion~1, since there are finitely many simple $A$-modules, up to isomorphism.\findemo
\begin{rem}{Remark} \label{direct summand}In particular, if $M$ has polynomial growth, then any direct summand of $M$ has polynomial growth. 
\end{rem}
Conversely, the class of modules with polynomial growth is closed under extensions. More precisely~:
\begin{mth}{Lemma} \label{2 out of 3}Let $A$ be a finite dimensional algebra over a field $k$.
\begin{enumerate}
\item If $L\stackrel{f}{\To} M\stackrel{g}{\To} N$ is an exact sequence of finite dimensional $k$-vector spaces, then
$$\dim_kM\leq \dim_kL+\dim_kN\mvirg$$
with equality if and only if $f$ is injective and $g$ is surjective.
\item Let $$0\To L\To M\To N\To 0$$
be a short exact sequence of finitely generated $A$-modules. If two of the modules $L$, $M$, and $N$ have polynomial growth, so does the third.
\end{enumerate}
\end{mth}
\pf For Assertion~1, there is a short exact sequence
$$0\To \Im\,f\To M\To \Im\,g\To 0\mvirg$$
hence $\dim_kM=\dim_k\Im\,f+\dim_k\Im\,g\leq \dim_kL+\dim_kN$. Equality holds if and only if $\dim_k\Im\,f=\dim_kL$ and $\dim_k\Im\,g=\dim_kN$, i.e. if $f$ is injective and $g$ is surjective.\par
For Assertion~2, let $S$ be a finitely generated $A$-module. Consider the long exact sequence of $\Ext$-groups
$$\cdots \to\Ext_A^n(N,S)\to \Ext_A^n(M,S) \to \Ext_A^n(L,S)\to\Ext_A^{n+1}(N,S)\to\cdots$$
It follows from Assertion~1 that
$$\dim_k\Ext_A^n(M,S)\leq \dim_k\Ext_A^n(N,S)+\dim_k\Ext_A^n(L,S)\mvirg$$
hence if $L$ and $N$ have polynomial growth, so does $M$. Similarly, 
$$\dim_k\Ext_A^n(N,S)\leq \dim_k\Ext_A^n(M,S)+\dim_k\Ext_A^{n-1}(L,S)\mvirg$$
It follows that $N$ has polynomial growth if $L$ ands $M$ have. Finally, 
$$\dim_k\Ext_A^n(L,S)\leq \dim_k\Ext_A^n(M,S)+\dim_k\Ext_A^{n+1}(N,S)\mvirg$$
hence $L$ has polynomial growth if $M$ and $N$ do.\findemo
\begin{mth}{Corollary} \label{poco simple}Let $A$ be a finite dimensional algebra over a field $k$. The following conditions are equivalent~:
\begin{enumerate}
\item Every finitely generated $A$-module has polynomial growth.
\item Every simple $A$-module has polynomial growth.
\end{enumerate}
\end{mth}
\pf Obviously Condition~1 implies Condition~2. The converse follows by induction on the length of a finitely generated $A$-module.\findemo
\section{Cohomological Mackey functors}\label{comack}
\masubsect{Definition}
Let $R$ denote an arbitrary commutative unital ring. Recall that {\em A Mackey functor $M$ for $G$ over $R$} consists of the assignment $H\mapsto M(H)$ of an $R$-module to each subgroup $H$ of $G$, together with maps of $R$-modules 
$$t_H^K:M(H)\To M(K)\mvirg$$
called {\em transfer maps}, and 
$$r_H^K:M(K)\To M(H)\mvirg$$
called {\em restriction maps}, whenever $H\leq K\leq G$, and maps of $R$ modules
$$c_{x,H}:M(H)\To M({^xH})\mvirg$$
for each $x\in G$, subject to a list of compatibility conditions, in particular the Mackey formula (cf. \cite{thevwebb} for details). A Mackey functor $M$ is called {\em cohomological} if the additional conditions
$$t_H^Kr_H^K=|K:H|\Id_{M(H)}$$
are fulfilled, for any $H\leq K\leq G$.\par
There is an obvious notion of morphism of (cohomological) Mackey functors, and this yields the category $\comack_R(G)$ of cohomological Mackey functors for $G$ over $R$. 
\begin{rem}{Example (fixed points functors)} In particular, when $V$ is an $RG$-module, then the fixed point functor $FP_V$ is the Mackey functor defined by $FP_V(H)=V^H$, i.e. the set of elements of $V$ which are invariant by $H$. When $H\leq K\leq H$, the transfer map $V^H\To V^K$ is {\em the relative trace map}~${\rm Tr}_H^K$, defined by
$$ {\rm Tr}_H^K(v)=\sum_{x\in [K/H]}yv\mvirg$$
for any $v\in V^H$, where $[K/H]$ is a transversal of $H$ in $K$. The restriction map $V^K\To V^H$ is the inclusion map, and for $x\in G$, the conjugation map $V^H\To V^{{^xH}}$ is the map $v\mapsto xv$. The functor $FP_V$ is obviously cohomological.\par
The correspondence $V\mapsto FP_V$ is a functor from $\gMod{RG}$ to $\comack_R(G)$. This functor is fully faithful~: if $V$ and $W$ are $RG$-modules, then any morphism $\varphi: FP_V\to FP_W$ in $\comack_k(G)$ is equal to $FP_f$, where $f$ is the morphism of $RG$-modules from $V=FP_V(\un)$ to $W=FP_W(\un)$ obtained by evaluating $\varphi$ at the trivial subgroup of $G$. 
\end{rem}
\pagebreak[3]
\masubsect{The cohomological Mackey algebra} It was shown by Th\'evenaz and Webb (cf. \cite{thevwebb}) that $\comack_R(G)$ is equivalent to the category of modules over the {\em cohomological Mackey algebra} $co\mu_R(G)$ for $G$ over $R$. More precisely, if~$M$ is a cohomological Mackey functor for $G$ over $R$, then the corresponding $co\mu_R(G)$-module is equal to $\dirsum{H\leq G}M(H)$.\par
The algebra $co\mu_R(G)$ is a finitely generated free $R$-module, so in particular when $R$ is a field, it is a finite dimensional $R$-algebra. In this case, there is a natural notion of projective cover of a finitely generated $co\mu_R(G)$-module, hence a natural notion of {\em minimal projective resolution} of a cohomological Mackey functor.
\masubsect{Yoshida's theorem}
Let $\per_R(G)$ denote the full subcategory of the category $\gmod{RG}$ of finitely generated $RG$-modules, consisting of {\em permutation} $RG$-modules, i.e. modules admitting a (globally) $G$-invariant $R$-basis. It is an $R$-linear category, which is naturally equivalent to the opposite category~: indeed, the dual $V^*=\Hom_R(V,R)$ of a finitely generated permutation $RG$-module is again a finitely generated permutation $RG$-module, and the correspondence $\delta: V\mapsto V^*$ is an equivalence of categories from $\per_R(G)$ to $\per_R(G)\op$, which is its own inverse (up to a slight abuse of notation).\par
When $M$ is a cohomological Mackey functor for $G$ over $R$, then the correspondence
$$V\mapsto \Hom_{\comack_R(G)}(FP_V,M)$$
is an $R$-linear contravariant functor $\tilde{M}$ from $\per_R(G)$ to $\gmod{R}$. This yields in turn an $R$-linear functor $M\mapsto \tilde{M}$ from $\comack_R(G)$ to $\fun_R(G)$, where $\fun_R(G)$ denotes the category of $R$-linear contravariant functors from $\per_R(G)$ to $\gmod{R}$.\par
Conversely, if $F$ is such a functor, and $H$ is a subgroup of $G$, one can define $\hat{F}(H)=F\big(R(G/H)\big)$, where $R(G/H)$ is the free $R$-module with basis the $G$-set $G/H$ of $H$-cosets in $G$. If $H\leq K\leq G$, then the projection map $p_H^K: G/K\To G/H$ gives a map of $RG$-modules $Rp_H^K:R(G/K)\To R(G/H)$, and taking image by $F$ gives a transfer map $t_H^K=F(p_H^K): \hat{F}(H)\To \hat{F}(K)$. Similarly, applying first the equivalence $\delta$ gives a map $r_H^K=F\big(\delta(Rp_H^K)\big):\hat{F}(K)\To \hat{F}(H)$. Finally, if $x\in G$, then the map $g{^xH}\mapsto gxH$ induce a map $R(G/{^xH})\To R(G/H)$, whose image by $F$ yields a conjugation map $c_{x,H}:\hat{F}(H)\To\hat{F}({^xH})$. With these definitions $\hat{F}$ becomes a cohomological Mackey functor for $G$ over $R$, and the correspondence $F\mapsto\hat{F}$ is a functor from $\fun_R(G)$ to $\comack_R(G)$. The following theorem is essential~:
\begin{mth}{Theorem} {\rm [Yoshida \cite{yoshida}]} \label{yoshida}The functors $M\mapsto\tilde{M}$ and $F\mapsto\hat{F}$ are mutual inverse equivalences of categories between  $\comack_R(G)$ and $\fun_R(G)$.
\end{mth}
\begin{rem}{Remark} One checks easily that if $V$ is an $RG$-module, then the functor $FP_V$ is mapped to the functor $\Hom_{RG}({-},V)$ by this equivalence. For this reason, this functor will also be denoted by $FP_V$. More generally, Yoshida's equivalence allows for an identification of $\comack_R(G)$ with $\fun_R(G)$, that will be used freely throughout the rest of this paper.
\end{rem}
\begin{rem}{Remark} In particular, if $V$ is a permutation $RG$-module, then the Yoneda functor $FP_V=\Hom_{\per_R(G)}({-},V)$ is a projective object in $\fun_R(G)$.  More precisely, if $M$ is a cohomological Mackey functor for $G$ over $R$, and $H$ is a subgroup of $G$, there is an isomorphism of $R$-modules
\begin{equation}\label{mult}
\Hom_{\comack_R(G)}(FP_{R(G/H)},M)\cong M(H)\mpoint
\end{equation}
It follows more generally that if $V$ is a direct summand of a permutation $RG$-module, then the functor $FP_V$ is a projective object in $\comack_R(G)$. Th\'evenaz and Webb (cf.~\cite{thevwebb}) have shown conversely that any projective object in $\comack_R(G)$ is isomorphic to $FP_V$, where $V$ is a direct summand of a permutation $RG$-module.
\end{rem}
\begin{rem}{Remark}\label{associated module}
It also follows that the category $\fun_R(G)$ is equivalent to the category of modules over the Hecke algebra
$$\mathsf{Y}_R(G)=\End_{RG}\big(\dirsum{H\leq G}R(G/H)\big)\mpoint$$
One can show easily that this algebra is actually isomorphic to the cohomological Mackey algebra $co\mu_R(G)$ (cf. \cite{green}). The $\mathsf{Y}_R(G)$-module corresponding to the object $F$ of $\fun_R(G)$ (resp. to the cohomological Mackey functor $M$ for $G$ over $R$) under these equivalences, is equal to $\dirsum{H\leq G}F\big(R(G/H)\big)$ (resp. to $\dirsum{H\leq G}M(H)$). In particular $F$ is finitely generated (resp. $M$ is finitely generated) if and only if $F(W)$ is a finitely generated $R$-module, for any finitely generated permutation $RG$-module $W$ (resp. $M(H)$ is a finitely generated $R$-module, for any $H\leq G$).
\end{rem}
\pagebreak[3]
\masubsect{The dual of a Mackey functor} If $M$ is a Mackey functor for $G$ over $R$, then the {\em dual} Mackey functor $M^*$ is defined by
$$\forall H\leq G,\;\;M^*(G)=\Hom_R\big(M(G),R\big)\mpoint$$
The transfer, restriction, and conjugation maps for $M^*$ are defined by
$$r_H^K={^\tau}(t_{H}^{K})\mvirg\;\;t_H^K={^\tau}(r_{H}^{K})\mvirg\;\;c_{x,H}={^\tau}(c_{x^{-1},{^xH}})\mvirg$$
for any $H\leq K\leq G$ and any $x\in G$, where the exponent $\tau$ denotes transposed maps.\par
If $M$ is cohomological, then $M^*$ is also cohomological. Through the equivalence given by Yoshida's Theorem~\ref{yoshida}, this duality maps the functor $F$ of $\fun_R(G)$ to the functor $F^*$ defined as the composition
$$\per_R(G)\stackrel{*}{\To}\per_R(G)\op\stackrel{F}{\To}\gMod{R}\stackrel{*}{\To}\gMod{R}\op\mpoint$$
\begin{rem}{Remark} \label{ext dual}The correspondence $M\mapsto M^*$ is a functor from $\comack_R(G)$ to the opposite category. The canonical morphism from $M$ to its bidual $(M^*)^*$ is functorial in $M$, and it is an isomorphism when $R$ is a field $k$ and $M$ is finitely generated. In other words, the correspondence $M\mapsto M^*$ induces an equivalence from the category of finitely generated cohomological Mackey functors for $G$ over $k$ to the opposite category. Thus, for any finitely generated cohomological Mackey functors $M$ and $N$ for $G$ over $k$, there is a natural isomorphism
$$\Hom_{\comack_k(G)}(M,N)\cong \Hom_{\comack_k(G)}(N^*,M^*)\mpoint$$
The functor $M$ is a finitely generated projective object in $\comack_k(G)$ if and only if $M^*$ is a finitely generated injective object in $\comack_k(G)$, and the previous isomorphism extends to natural isomorphisms
$$\Ext^n_{\comack_k(G)}(M,N)\cong \Ext^n_{\comack_k(G)}(N^*,M^*)\mvirg$$
for any $n\in \N$.
\end{rem}
\masubsect{Construction of functors}
Let $G$ and $H$ be finite groups. If $U$ is a finite $(H,G)$-biset consider the $R$-linear functor
$$\mathsf{t}_U: V\mapsto RU\otimes_{RG}V$$
from $\gmod{RG}$ to $\gmod{RH}$ maps permutation $RG$-modules to permutation $RH$-modules. By composition, this induces a functor
$$\mathsf{L}_U: F\mapsto F\circ\mathsf{t}_U$$
from $\fun_R(H)$ to $\fun_R(G)$. Since the right adjoint to the functor $\mathsf{t}_U$ is the functor
$$\mathsf{h}_U:W\mapsto \Hom_{RH}(RU,W)\mvirg$$
it follows from standard results of category theory that the functor
$$\mathsf{R}_U: F\mapsto F\circ \mathsf{h}_U$$
from $\fun_R(G)$ to $\fun_R(H)$, is right adjoint to $\mathsf{L}_U$.
\begin{rem}{Remark} The functors $\mathsf{L}_U$ and $\mathsf{R}_U$ are a generalization of functors considered by Tambara (see \cite{tambara} Section~4).
\end{rem}
\begin{mth}{Proposition} \label{composition}\begin{enumerate}
\item Let $G$ and $H$ be finite groups, and $U$ be a finite $(H,G)$-biset. Then the functors $\mathsf{R}_U$ and $\mathsf{L}_U$ are exact.
\item Let $G$ and $H$ be finite groups. If $U$ and $U'$ are isomorphic finite $(H,G)$-bisets, then there are isomorphisms of functors $\mathsf{L}_U\cong \mathsf{L}_{U'}$ and $\mathsf{R}_U\cong \mathsf{R}_{U'}$.
\item Let $G$, $H$, and $K$ be finite groups. Let $U$ be a finite $(H,G)$-biset, and~$V$ be a finite $(K,H)$-biset. Then there are isomorphisms of functors
$$\mathsf{R}_V\circ\mathsf{R}_U\cong \mathsf{R}_{V\times_HU}\mvirg\;\mathsf{L}_U\circ\mathsf{L}_V\cong \mathsf{L}_{V\times_HU}\mpoint$$
\item Let $G$ be a finite group, and let $\Id_G$ denote the identity biset for $G$, i.e. the set $G$ for its $(G,G)$-biset structure given by multiplication. Then the functor $\mathsf{L}_{\Id_G}$ and $\mathsf{R}_{\Id_G}$ are isomorphic to the identity functor.
\item Let $G$ and $H$ be finite groups. If $U$ and $U'$ are finite $(H,G)$-bisets, there are isomorphisms of functors
$$\mathsf{L}_{U\sqcup U'}\cong\mathsf{L}_U\oplus\mathsf{L}_{U'}\mvirg\;\mathsf{R}_{U\sqcup U'}\cong\mathsf{R}_U\oplus\mathsf{R}_{U'}\mpoint$$
\item Let $G$ and $H$ be finite groups, and $U$ be a finite $(H,G)$-biset. Then for any object $F$ of $\fun_R(G)$, there is an isomorphism 
$$\mathsf{R}_{U\op}(F)^*\cong \mathsf{L}_U(F^*)$$
in $\fun_R(H)$, which is functorial in $F$.
\end{enumerate}
\end{mth}
\pf Assertion~1 is obvious, since  the functors $\mathsf{R}_U$ and $\mathsf{L}_U$ are obtained by pre-composition with some functor (in other words, they are both {\em restriction functors} along a suitable functor).\par
Assertion~2 is a consequence of the isomorphisms of functors $\mathsf{t}_U\cong\mathsf{t}_{U'}$, and  $\mathsf{h}_U\cong\mathsf{h}_{U'}$, which both follow from the isomorphism $RU\cong RU'$ of $(RH,RG)$-bimodules.\par
For Assertion~3, the associativity of tensor product gives an isomorphism of functor
$$\mathsf{t}_V\circ\mathsf{t}_U\cong\mathsf{t}_{V\times_HU}\mvirg$$
which by adjunction, gives the isomorphism of functors
$$\mathsf{h}_U\circ\mathsf{h}_V\cong\mathsf{h}_{V\times_HU}\mpoint$$
The isomorphisms of Assertion~3 follow by composition.\par
Assertion~4 follows from the fact that the functors $\mathsf{t}_{\Id_G}$ and $\mathsf{t}_{\Id_G}$ are both isomorphic to the identity functor. Similarly, Assertion~5 follows from similar additivity properties of the functors $t_U$ and $h_U$ with respect to $U$.\par
For Assertion~6, let $W$ be an object of $\per_R(G)$. Then
\begin{eqnarray*}
\mathsf{R}_{U\op}(F)^*(W)&=&\Hom_R\big(\mathsf{R}_{U\op}(F)(W^*),R\big)\\
&=&\Hom_R\left(F\Big(\Hom_{RG}\big(RU\op,\Hom_R(W,R)\big)\Big),R\right)\\
&\cong&\Hom_R\Big(F\big(\Hom_{R}(RU\otimes_{RG}W,R)\big),R\Big)\\
&=&\mathsf{L}_U(F^*)(W)\mvirg
\end{eqnarray*}
and these isomorphisms are functorial with respect to $W$.
\findemo
\begin{mth}{Proposition} \label{proj to proj}Let $G$ and $H$ be finite groups, and $U$ be a finite $(H,G)$-biset. 
\begin{enumerate}
\item If $V$ is an $RH$-module, then
$$\mathsf{L}_U(FP_V)=FP_{\Hom_{RH}(RU,V)}\mpoint$$
In particular, the functor $\mathsf{L}_U$ maps projective objects to projective objects.
\item If $G$ acts freely on $U$, then $\mathsf{L}_U\cong\mathsf{R}_{U\op}$, where $U\op$ denotes the opposite biset of $U$.
\end{enumerate}
\end{mth}
\pf For Assertion~1, there is an isomorphism of functors
$$\Hom_{RG}\big({-},\Hom_{RH}(RU,{V})\big)\cong \Hom_{RH}(RU\otimes_{RG}{-},{V})\mpoint$$
In other words $FP_V\circ\mathsf{t}_U\cong FP_{\Hom_{RH}(RU,V)}$. \par
The last part of Assertion~1 follows from the facts that if $V$ is a permutation $RH$-module, then $\Hom_{RH}(RU,V)=\mathsf{h}_U(V)$ is a permutation $RG$-module. An alternative proof consists in observing that since $\mathsf{L}_U$ is left adjoint to an exact functor, it maps projective objects to projective objects.\par
For Assertion~2, for any $RG$-module $W$, there is an isomorphism of $RH$-modules
$$RU\otimes_{RG}W\To\Hom_{RG}(RU\op,W)\mvirg$$
defined by sending $u\otimes w$, for $u\in U$ and $w\in W$, to the map sending $v\in U\op$ (recall that $U=U\op$ as a set) to $gw$, if there exists an element $g\in G$ such that $v=ug^{-1}$ (and in this case there is a unique such $g$, since $G$ acts freely on $U$), and to 0 otherwise. This isomorphism is obviously functorial in $W$, and this completes the proof.
\findemo
\begin{mth}{Proposition}\label{ext adjunction} Let $G$ and $H$ be finite groups, and $U$ be a finite $(H,G)$-biset. Then for any $n\in \N$, the adjunction of the pair $(\mathsf{L}_U,\mathsf{R}_U)$ induces an isomorphism of bifunctors 
$$\Ext_{\fun_R(G)}^n\big(\mathsf{L}_U({-}),{-}\big)\cong\Ext_{\fun_R(H)}^n\big({-},\mathsf{R}_U({-})\big)\mpoint$$
\end{mth}
\pf This follows from the fact that $\mathsf{L}_U$ and $\mathsf{R}_U$ are both exact functors.\findemo
\begin{rem}{Remark} \label{Yoneda adjunction}In terms of Yoneda extensions, this isomorphism can be viewed as follows~: the adjunction of the pair $(\mathsf{L}_U,\mathsf{R}_U)$ is equivalent to the existence of natural transformation of functors
$$\eta:\Id \To\mathsf{R}_U\circ\mathsf{L}_U\;\;\hbox{and}\;\;\varepsilon:\mathsf{L}_U\circ\mathsf{R}_U\To\Id\mvirg$$
called respectively the {\em unit} and {\em counit} of the adjunction, with the property that for each object $M$ in $\fun_R(G)$ and each object $N$ in $\fun_R(H)$,
$$\varepsilon_{\mathsf{L}_U(N)}\circ\mathsf{L}_U(\eta_N)=\Id_{\mathsf{L}_U(N)}\;\;\hbox{and}\;\;\mathsf{R}_U(\varepsilon_M)\circ\eta_{\mathsf{R}_U(M)} =\Id_{\mathsf{R}_U(M)}\mpoint$$
The bijection $\Hom_{\fun_R(G)}\big(\mathsf{L}_U({N}),{M}\big)\To\Hom_{\fun_R(H)}\big({N},\mathsf{R}_U({M})\big)$ is given by $\alpha\mapsto \mathsf{R_U}(\alpha)\circ\eta_N$. The inverse bijection is given by $\beta\mapsto \varepsilon_M\circ\mathsf{L}(\beta)$.\par
In other words, these bijections consist in taking images by one of the functors, and then compose on the suitable side with unit or counit.\par
Now interpreting the extension group $\Ext_{\fun_R(G)}^n\big(\mathsf{L}_U({N}),{M}\big)$ a the set of equivalence classes of exact sequences in $\fun_R(G)$ of the form
$$
\xymatrix{
0\ar[r]& M\ar[r]& X_{n-1}\ar[r]&\cdots\ar[r]& X_1\ar[r]& X_0\ar[r]&\mathsf{L}_U({N})\ar[r]& 0\mvirg\\
}
$$
the procedure is the same~: first apply the functor $\mathsf{R}_U({-})$, to get an exact sequence
$$
\xymatrix@C=2.5ex{
0\ar[r]& \mathsf{R}_U(M)\ar[r]& \mathsf{R}_U(X_{n-1})\ar[r]&\cdots\ar[r]& \mathsf{R}_U(X_1)\ar[r]& \mathsf{R}_U(X_0)\ar[r]&\mathsf{R}_U\circ\mathsf{L}_U({N})\ar[r]& 0\mvirg\\
}
$$
and then compose with the map $\eta_N$, i.e. complete the cartesian square at the right of the following diagram
$$
\xymatrix@C=2.5ex{
0\ar[r]& \mathsf{R}_U(M)\ar[r]& \mathsf{R}_U(X_{n-1})\ar[r]&\cdots\ar[r]& \mathsf{R}_U(X_1)\ar[r]& \mathsf{R}_U(X_0)\ar[r]&\mathsf{R}_U\circ\mathsf{L}_U({N})\ar[r]& 0\\
0\ar[r]& \mathsf{R}_U(M)\ar[r]\ar@{=}[u]& \mathsf{R}_U(X_{n-1})\ar[r]\ar@{=}[u]&\cdots\ar[r]& \mathsf{R}_U(X_1)\ar[r]\ar@{=}[u]& Y_0\ar[r]\ar[u]&N\ar[r]\ar[u]_-{\eta_N}& 0\mpoint\\
}
$$
The inverse bijection is obtained similarly by first applying the functor $\mathsf{L}_U$, and composing with the map $\varepsilon_M$, i.e. completing a cocartesian square at the left of the resulting diagram. 
\end{rem}
\begin{rem}{Remark} \label{Yoneda compatible}It follows easily that the isomorphisms of functors 
$$\alpha_n: \Ext_{\fun_R(G)}^n\big(\mathsf{L}_U({-}),{-}\big)\stackrel{\cong}{\To}\Ext_{\fun_R(H)}^n\big({-},\mathsf{R}_U({-})\big)\mpoint$$
of Proposition~\ref{ext adjunction}, are compatible with the Yoneda product, in the following sense~: if $P$ is an object of $\fun_R(H)$, and if $M$ and $N$ are objects of $\fun_R(G)$, if $e\in \Ext_{\fun_R(G)}^n\big(\mathsf{L}_U({P}),{N}\big)$ and $f\in \Ext_{\fun_R(G)}^m({N},M)$, then
$$\alpha_{m+n}(f\circ e)=\mathsf{R}_U(f)\circ \alpha_n(e)\in\Ext_{\fun_R(H)}^{m+n}\big({P},\mathsf{R}_U({M})\big)\mpoint$$
\end{rem}
\begin{rem}{Example (induction and restriction)} Let $G$ be a subgroup of $H$. Set $U=H$, viewed as an $(H,G)$-biset by left and right multiplication. Then the functor $V\mapsto RU\otimes_{RG}V$ from $\gmod{RG}$ to $\gmod{RH}$ is isomorphic to the induction functor $V\mapsto\Ind_G^HV$. It follows that the functor $\mathsf{L}_U$ is isomorphic to the {\em restriction} functor $\Res_G^H:\comack_R(H)\To \comack_R(G)$ in this case. \par
In the same situation, the functor $W\mapsto \Hom_{RH}(RU,W)$ from $\gmod{RH}$ to $\gmod{RG}$ is isomorphic to the restriction functor $\Res_G^H$. It follows that the functor $\mathsf{R}_U$ is isomorphic to the {\em induction} functor $\Ind_G^H:\comack_R(G)\To \comack_R(H)$.
\end{rem}
\begin{rem}{Example (restriction and induction)}  Suppose now that $H$ is a subgroup of $G$, and consider $U=G$ as an $(H,G)$-biset by left and right multiplication. Then the functor $V\mapsto RU\otimes_{RG}V$ from $\gmod{RG}$ to $\gmod{RH}$ is isomorphic to the restriction functor $V\mapsto\Res_H^GV$. It follows that the functor $\mathsf{L}_U$ is isomorphic to the {\em induction} functor $\Ind_H^G:\comack_R(H)\To \comack_R(G)$ in this case. \par
In the same situation, the functor $W\mapsto \Hom_{RH}(RU,W)$ from $\gmod{RH}$ to $\gmod{RG}$ is isomorphic to the induction functor $\Res_G^H$. It follows that the functor $\mathsf{R}_U$ is isomorphic to the {\em induction} functor $\Ind_G^H:\comack_R(G)\To \comack_R(H)$.\par
This yields another proof of the well known fact that the induction and restriction functors between categories of cohomological Mackey functors are left and right adjoint to each other (cf. \cite{thevwebb} for details).
\end{rem}
\begin{rem}{Example (the functor $\rho_{G/N}^G$)} Let $H$ be a finite group, and $N$ be a normal subgroup of~$H$, and set $G=H/N$. Also set $U=G$, viewed as an $(H,G)$-biset in the obvious way. In this case, the functor $\mathsf{t}_U$ is the inflation functor from $\gmod{RG}$ to $\gmod{RH}$. Using the equivalences of categories of Theorem~\ref{yoshida}, the functor~$\mathsf{L}_U$ gives a functor denoted by $\rho_{H/N}^H$, from $\comack_R(H)$ to $\comack_R(H/N)$. One checks easily that if $M$ is a cohomological Mackey functor for $H$ over $R$, then, denoting by $x\mapsto\sur{x}$ the projection map $H\To H/N$
$$\big(\rho_{H/N}^H(M)\big)(\sur{K})=M(K)\mpoint$$
Similarly, the transfer, restriction, and conjugation maps for the functor $\rho_{H/N}^H(M)$ are obtained by just ``removing the bars'', i.e.
$$t_{\sur{K}}^{\sur{L}}=t_K^L\mvirg r_{\sur{K}}^{\sur{L}}=r_K^L\mvirg \;\;c_{\sur{x},\sur{K}}=c_{x,K}\mpoint$$
The right adjoint to $\rho_{H/N}^H$ is the functor $\mathsf{R}_U$, hereafter denoted by $\jiota_{G/N}^G$.\par
Moreover Assertion~3 of Proposition~\ref{proj to proj} shows that the functor $\rho_{G/N}^G$ is also equal to the functor $\mathsf{R}_{U\op}$. In particular, the left adjoint to $\rho_{G/N}^G$ is the functor $\mathsf{L}_{U\op}$, hereafter denoted by $\iiota_{G/N}^G$.
\end{rem}
\begin{rem}{Example (the functor $\iiota_{G/N}^G$)} The left adjoint to the functor $\rho_{G/N}^G$ will be denoted by $\iiota_{G/N}^G$. It is obtained as follows~: if $M$ is a cohomological Mackey functor for $G/N$, and $K$ is a subgroup of $G$, then
$$\big(\iiota_{G/N}^G(M)\big)(K)=M(KN/N)\mpoint$$
If $K\leq L\leq G$, then the transfer, restriction, and conjugation maps for the functor $\iiota_{G/N}^G(M)$ are given by
$$t_K^L=|L\cap N:K\cap N|t_{KN/N}^{LN/N}\mvirg \;\;r_K^L=r_{KN/N}^{LN/N}\mvirg \;\;c_{x,K}=c_{xN,KN/N}\mpoint$$
\end{rem}
\begin{rem}{Example (the functor $\jiota_{G/N}^G$)} The right adjoint to the functor $\rho_{G/N}^G$ will be denoted by $\jiota_{G/N}^G$. It is obtained as follows~: if $M$ is a cohomological Mackey functor for $G/N$, and $K$ is a subgroup of $G$, then
$$\big(\jiota_{G/N}^G(M)\big)(K)=M(KN/N)\mpoint$$
If $K\leq L\leq G$, then the transfer, restriction, and conjugation maps for the functor $\jiota_{G/N}^G(M)$ are given by
$$t_K^L=t_{KN/N}^{LN/N}\mvirg \;\;r_K^L=|L\cap N:K\cap N|r_{KN/N}^{LN/N}\mvirg \;\;c_{x,K}=c_{xN,KN/N}\mpoint$$
\end{rem}
\begin{rem}{Remark} The functors $\iiota_{G/N}^G$ and $\jiota_{G/N}^G$ should not be confused with the inflation functor $\Inf_{G/N}^G$~: recall (cf. \cite{thevwebb}) that if~$N$ is a normal subgroup of a finite group $G$, and $M$ is a Mackey functor for $G/N$ over $R$, then the functor $\Inf_{G/N}^GM$ is the Mackey functor for $G$ over $R$ defined by
$$(\Inf_{G/N}^GM)(H)=\left\{\begin{array}{cl}M(H/N)&\hbox{if}\;H\geq N\\\zero&\hbox{otherwise}\end{array}\right.\mpoint$$
The transfer, restriction, and conjugation maps are the obvious ones. In general, this inflation procedure does not preserve cohomological Mackey functors~: indeed, if $M$ is cohomological, and if $K<H$ are subgroups of $G$ such that $H\geq N$ but $K\not\geq N$, then the composition 
$$t_K^Hr_K^H:(\Inf_{G/N}^GM)(H)\To(\Inf_{G/N}^GM)(H)$$
is equal to 0, so it is not equal in general to the multiplication by $|H:K|$. \par
However, if $N$ is a $p$-group, for some prime number $p$,  then $p$ divides $|KN:K|=|N:K\cap N|$, hence $p$ divides $|H:K|$. If moreover $R$ has characteristic~$p$, then $|H:K|=0$ in $R$ in this situation, and one checks then that $\Inf_{G/N}^GM$ is cohomological. In this case, the inflation functor $\comack_R(G/N)\To \comack_R(G)$ corresponds via Yoshida's equivalence of categories to the functor $\fun_R(G/N)\To \fun_R(G)$ obtained by composition with the {\em Brauer quotient} functor $V\mapsto V[N]$ from $\per_R(G)$ to $\per_R(G/N)$.
\end{rem}
\begin{rem}{Example : group isomorphism} Let $f:G\to H$ be a group isomorphism. Then the set $U=H$ has a natural structure of $(H,G)$-biset, for which $h\in H$ acts by left multiplication by $h$, and $g\in G$ acts by right multiplication by $f(g)$. It is clear in this case that the functor $\mathsf{R}_U$ is the {\em transport by isomorphism} via $f$. This functor will be denoted by $\Iso(f)$. The functor $\mathsf{L}_U$ is isomorphic to $\Iso(f^{-1})$, by Assertion~3 of Proposition~\ref{proj to proj}.
\end{rem}
Recall that a {\em section} $(B,A)$ of $G$ is a pair of subgroups of $G$ such that $A\normal B$.
\begin{mth}{Lemma} \label{ind iota}Let $G$ be a finite group, and $(B,A)$ be a section of $G$. Then there is an isomorphism
$$\Ind_B^G\,\iiota_{B/A}^B\cong\mathsf{L}_{A\dom G}$$
of functors from $\comack_R(B/A)$ to $\comack_R(G)$, where $A\dom G$ is endowed with its natural $(B/A,G)$-biset structure. Similarly, there is an isomorphism 
$$\rho_{B/A}^B\Res_B^G\cong\mathsf{L}_{G/A}$$
of functors from $\comack_R(G)$ to $\comack_R(B/A)$, where $G/A$ is endowed with its natural $(G,B/A)$-biset structure.
\end{mth}
\pf The functor $\Ind_B^G\,\iiota_{B/A}^B$ is equal to the composition of $\Ind_B^G=\mathsf{L}_V$, where $V$ is the set $G$ for its natural $(B,G)$-biset structure, and $\iiota_{B/A}^B=\mathsf{L}_U$, where $U$ is the set $B/A$, for its natural $(B/A,B)$-biset structure. It follows from Proposition~\ref{composition} that $\Ind_B^G\iiota_{B/A}^B\cong\mathsf{L}_W$, where $W=(B/A)\times_BG$ is clearly isomorphic to the $(B/A,G)$-biset $A\dom G$.\par
By adjunction, it follows that $\rho_{B/A}^B\Res_B^G\cong\mathsf{R}_{A\dom G}$. But $\mathsf{R}_{A\dom G}\cong\mathsf{L}_{G/A}$ by Assertion~3 of Proposition~\ref{proj to proj}, since $B/A$ acts freely on $G/A$.\findemo
\begin{mth}{Proposition} \label{general Mackey}Let $G$ be a finite group, and $(B,A)$, $(D,C)$ be two sections of $G$. Then there is an isomorphism
$$\rho_{D/C}^D\;\Res_D^G\;\Ind_B^G\;\iiota_{B/A}^B\cong\dirsum{g\in[D\dom G/B]}\Ind_{\sur{D}_g}^{\sur{D}}\;\iiota_{\sur{D}_g/\sur{C}_g}^{\sur{D}_g}\;\Iso(f_g)\;\rho_{\sur{B}_g/\sur{A}_g}^{\sur{B}_g}\;\Res_{\sur{B}_g}^{\sur{B}}$$
of functors from $\comack_k(B/A)$ to $\comack_k(D/C)$, where
$$\sur{D}=D/C, \;\;\sur{D}_g=(D\cap{^gB})C/C, \;\;\sur{C}_g=(D\cap{^gA})C/C\mvirg$$
$$\sur{B}=B/A, \;\;\sur{B}_g=(D^g\cap B)A/A, \;\;\sur{A}_g=(C^g\cap B)A/A\mvirg$$
and where $f_g:\sur{B}_g/\sur{A}_g\To\sur{D}_g/\sur{C}_g$ is the group isomorphism sending $x\sur{A}_g$ to ${^gx}\sur{C}_g$, for $x\in D^g\cap B$.
\end{mth}
\pf (see Proposition7.1 of \cite{both2} for details) By Lemma~\ref{ind iota} and Proposition~\ref{composition}, there are isomorphisms of functors
\begin{eqnarray*}
\rho_{D/C}^D\Res_D^G\Ind_B^G\iiota_{B/A}^B&\cong&\mathsf{L}_{G/C}\circ \mathsf{L}_{A\dom G}\\
&\cong&\mathsf{L}_{A\dom G\times_GG/C}\cong\mathsf{L}_{A\dom G/C}\mpoint
\end{eqnarray*}
Now the $(B/A,D/C)$-biset $A\dom G/C$ splits as a disjoint union of transitive ones
$$A\dom G/C\cong\bigsqcup_{g\in [D\dom G/B]}A\dom Bg^{-1}D/C\mpoint$$
Moreover for each $g\in [D\dom G/B]$, with the above notation, there is an isomorphism of $(\sur{B},\sur{D})$-bisets
$$A\dom Bg^{-1}D/C\cong (\sur{B}/\sur{A}_g)\times_{\sur{B}_g/\sur{A}_g}\Iso(f_g)\times_{\sur{D}_g/\sur{C}_g}(\sur{C}_g\dom\sur{D})\mpoint$$
By Lemma~\ref{ind iota} and Proposition~\ref{composition} again, there is an isomorphism of functors
$$\mathsf{L}_{A\dom Bg^{-1}D/B}\cong \Ind_{\sur{D}_g}^{\sur{D}}\circ\iiota_{\sur{D}_g/\sur{C}_g}^{\sur{D}_g}\circ\Iso(f_g)\circ\rho_{\sur{B}_g/\sur{A}_g}^{\sur{B}_g}\circ\Res_{\sur{B}_g}^{\sur{B}}\mvirg$$
as was to be shown.\findemo
\section{Reduction to $p$-groups}\label{reduction}
\begin{mth}{Lemma} \label{pol to pol}Let $k$ be a field, let $G$ and $H$ be finite groups, and $U$ be a finite $(H,G)$-biset. 
\begin{enumerate}
\item If $F$ is an object of $\fun_k(H)$, then
$$\dim_k\mathsf{L}_U(F)\leq C_U\dim_kF\mvirg$$
where $C_U=\mathop{\sum}_{K\leq G}\limits|H\dom U/K|$.
\item If $F$ is an object of $\fun_k(H)$, and if $F$ has polynomial growth, then so does $\mathsf{L}_U(F)$.
\end{enumerate}
\end{mth}
\pf
Let $F$ be an object of $\fun_k(H)$. The $co\mu_k(G)$-module corresponding to $\mathsf{L}_U(F)$ is equal to $\dirsum{K\leq G}\mathsf{L}_U(F)\big(k(G/K)\big)$, by Remark~\ref{associated module}. Thus
\begin{eqnarray*}
\dim_k\mathsf{L}_U(F)&=&\sum_{K\leq G}\dim_kF\big(kU\otimes_{kG}k(G/K)\big)\\
&=&\sum_{K\leq G}\dim_kF\big(k(U/K)\big)\\
&=&\sum_{K\leq G}\sum_{u\in [H\dom U/K]}\dim_kF\big(k(H/H_{uK})\big)\mvirg
\end{eqnarray*}
where $H_{uK}=\{h\in H\mid\exists g\in K,\;hu=ug\}$. Now for each $K\leq G$ and each $u\in [H\dom U/K]$
$$\dim_kF\big(k(H/H_{uK})\big)\leq \dim_kF=\sum_{L\leq H}\dim_kF\big(k(H/L)\big)\mvirg$$
thus
$$\dim_k\mathsf{L}_U(F)\leq \Big(\sum_{K\leq G}\sum_{u\in [H\dom U/K]}1\Big)\dim_kF\mvirg$$
showing Assertion~1.\par
Now if
$$\cdots P_n\To P_{n-1}\To\cdots\To P_0\To F\To 0$$
is a projective resolution with polynomial growth, there are constants $c$, $d$, and $e$, such that $\dim_kP_n\leq cn^d+e$, for all $n\in\N$. Now the complex
$$\cdots \mathsf{L}_U(P_n)\To \mathsf{L}_U(P_{n-1})\To\cdots\To \mathsf{L}_U(P_0)\To \mathsf{L}_U(F)\To 0$$
 is a projective resolution  $\mathsf{L}_U(F)$, since $\mathsf{L}_U$ is exact and maps projective objects to projective objects. Moreover
$$\dim_k \mathsf{L}_U(P_n)\leq C_U\dim_kP_n\leq C_Ucn^d+C_Ue\mvirg$$
so $\mathsf{L}_U(F)$ has polynomial growth.\findemo
\begin{mth}{Proposition} \label{section}Let $G$ be a finite group, and $(B,A)$ be a section of $G$. Let $k$ be a field, and $N$ be a cohomological Mackey functor for $B/A$ over $k$. Then the following conditions are equivalent~:
\begin{enumerate}
\item The functor $N$ has polynomial growth.
\item The functor $\Ind_B^G\,\iiota_{B/A}^B(N)$ has polynomial growth.
\end{enumerate}
\end{mth}
\pf The functor $\Ind_B^G\,\iiota_{B/A}^B$ is isomorphic to $\mathsf{L}_{W}$, where $W$ is the $(B/A,G)$ biset $A\dom G$, by Lemma~\ref{ind iota}, so Condition~1 implies Condition~2, by Lemma~\ref{pol to pol}.\par
Conversely, if $L_{W}(N)$ has polynomial growth, then so does $\mathsf{L}_{W\op}\circ\mathsf{L}_W(N)\cong\mathsf{L}_{W\times_GW\op}(N)$, by Lemma~\ref{pol to pol}. But $W\times_GW\op\cong A\dom G/A$, as a $(B/A,B/A)$-biset. In particular, it is the disjoint union of the identity biset $\Id_{B/A}=B/A$ and some other $(B/A,B/A)$-biset. By Assertions~4 and~5 of Proposition~\ref{composition}, it follows that the identity functor is a direct summand of $\mathsf{L}_{W\times_GW\op}$. In particular $N$ is isomorphic to a direct summand of $\mathsf{L}_{W\times_GW\op}(N)$, and $N$ has polynomial growth, by Remark~\ref{direct summand}. Thus Condition~2 implies Condition~1.\findemo
\begin{mth}{Lemma} \label{IndRes}Let $G$ be a finite group and $M$ be a cohomological Mackey functor for $G$ over $R$. Then for any subgroup $H$ of $G$ the composition
$$M\To \Ind_H^G\Res_H^GM\To M$$
of the unit and counit morphisms of the adjoint pairs of functors $(\Res_H^G,\Ind_H^G)$ and $(\Ind_H^G,\Res_H^G)$ is equal to the multiplication by $|G:H|$.
\end{mth}
\pf This is because the same is true for the categories of $RG$-modules and $RH$-modules, and the induction and restriction functors between them.\findemo
\begin{mth}{Proposition} \label{Sylow}Let $k$ be a field of positive characteristic $p$. Let $G$ be a finite group, and $S$ be a subgroup of $G$, containing a Sylow $p$-subgroup of $G$. Then for any cohomological Mackey functor $M$ for $G$ over $k$, the following conditions are equivalent~:
\begin{enumerate}
\item The functor $M$ has polynomial growth.
\item The functor $\Res_S^GM$ has polynomial growth.
\end{enumerate}
\end{mth}
\pf Condition~1 implies Condition~2, by Lemma~\ref{pol to pol}. Conversely, if $\Res_S^GM$ has polynomial growth, so does $\Ind_S^G\Res_S^GM$, by Lemma~\ref{pol to pol} again. Now Lemma~\ref{IndRes} shows that $M$ is a direct summand of $\Ind_S^G\Res_S^GM$, since $|G:S|$ is non zero in $k$. Hence $M$ has polynomial growth, by Remark~\ref{direct summand}.\findemo
\begin{mth}{Definition}
Let $k$ be a field. A finite group $G$ is called a {\em poco} group over $k$ if every finitely generated cohomological Mackey functor for~$G$ over $k$ has polynomial growth.
\end{mth}
\pagebreak[3]
\begin{mth}{Corollary}\label{reduction to Sylow} Let $k$ be a field of characteristic $p$. Let $G$ be a finite group, and $S$ be a Sylow $p$-subgroup of $G$. Then $G$ is a poco group over $k$ if and only if $S$ is a poco group over $k$.
\end{mth}
\pf Indeed, if $G$ is a poco group, and $N$ is a cohomological Mackey functor for $S$ over $k$, then $\Ind_S^GN$ has polynomial growth, and $N$ has polynomial growth, by Proposition~\ref{section}, applied to the section $(S,\un)$ of $G$. So $S$ is a poco group.\par
Conversely, if $S$ is a poco group, and $M$ is a cohomological Mackey functor for $G$ over $k$, the $\Res_S^GM$ has polynomial growth, so $M$ has polynomial growth, by Proposition~\ref{Sylow}.\findemo
\begin{mth}{Proposition} \label{FP_V is enough}Let $k$ be a field, and $G$ be a finite group. The following conditions are equivalent~:
\begin{enumerate}
\item The group $G$ is a poco group over $k$.
\item For any finitely generated $kG$-module $V$, the functor $FP_V$ has a polynomial growth.
\end{enumerate}
\end{mth}
\pf Obviously Condition~1 implies Condition~2. Conversely, observe that since any morphism from $FP_L$ to $FP_N$, where $L$ and $N$ are $kG$-modules, is determined by a morphism of $kG$-modules from $L$ to $N$, for any cohomological Mackey functor $M$ for $G$ over $k$, there is a short exact sequence
$$0\To FP_V\To FP_L\To FP_N\To M\To 0\mvirg$$
where $L$ and $N$ are permutation $kG$-modules, and $V$ is the kernel of a morphism of $kG$-modules from $L$ to $N$. Now if $FP_V$ has polynomial growth, so does $M$.\findemo
\section{Proof of Theorem~\ref{main}}\label{sketch}
Corollary~\ref{reduction to Sylow} shows that $G$ is a poco group if and only if $S$ is a poco group, where $S$ is a Sylow $p$-subgroup of $G$. Now by Proposition~\ref{section}, if $(B,A)$ is a section of $S$, the factor group $B/A$ is also a poco group. \par
It will be shown in Sections~\ref{cyclic p-groups} and~\ref{elemab} that an elementary abelian $p$-group is a poco group if and only if it has rank at most 1, when $p>2$, or at most~2, when $p=2$. It will follow that if $G$ is a poco group, then $S$ has sectional rank at most 1 if $p>2$, which implies that $S$ is cyclic,
or at most 2 if $p=2$. In other words, Condition~1 of Theorem~\ref{main} implies Condition~2.\par
Conversely, it will be shown in Section~\ref{cyclic p-groups} that cyclic $p$-groups are poco groups. So Condition~2 of Theorem~\ref{main} implies Condition~1, when $p>2$. The corresponding assertion for $p=2$ will be stated in Section~\ref{sectional 2-rank 2}, completing the proof of Theorem~\ref{main}.
\pagebreak[3]
\section{Simple cohomological functors for $p$-groups} \label{simple functors}Recall (see~\cite{thevwebb} for details) that if $G$ is a $p$-group, the simple cohomological Mackey functors for $G$ over $k$, up to isomorphism, are in one to one correspondence with the subgroups of $G$, up to conjugation. The simple functor~$S_Q^G$ (also denoted by $S_Q$ if $G$ is clear from context) corresponding to the subgroup~$Q$ is defined by
$$\forall T\leq G,\;\;S_Q^G(T)=\left\{\begin{array}{rl}k&\hbox{if}\;\;T=_GQ\\\zero&\hbox{otherwise.}\end{array}\right.$$
The projective cover of the functor $S_Q^G$ is the fixed point functor $FP_{kG/Q}$. Moreover $\End_{\comack_k(G)}(S_Q^G)\cong k$, and $S_Q^G$ is self dual (i.e. $(S_Q^G)^*\cong S_Q^G$).
\begin{mth}{Lemma} \label{reso S_un}Let $k$ be a field of characteristic $p$, and $G$ be a $p$-group. There is a short exact sequence in $\comack_k(G)$
$$0\To FP_{\Omega_G}\To FP_{kG}\To S_\un^G\To 0\mvirg$$
where $\Omega_G$ is the kernel of the augmentation map $\varepsilon: kG\to k$.
\end{mth}
\pf Consider the exact sequence of $kG$-modules
$$0\To \Omega_G\To kG\stackrel{\varepsilon}{\To} k\To 0\mpoint$$
Since fixed point functors are left exact, the inclusion $\Omega_G\subseteq kG$ yields an inclusion $FP_{\Omega_G}\subseteq FP_{kG}$. Let $S$ denote the quotient functor. In particular $S(\un)\cong k$. And if $Q$ is a non-trivial subgroup of $G$, then $(kG)^Q={\rm Tr}_\un^Q(kG)\subseteq \Omega_{G}$. Thus $(\Omega_G)^Q=(kG)^Q$, and $S(Q)=\zero$. Hence $S\cong S_\un^G$.\findemo
\begin{mth}{Theorem} {\rm [Samy Modeliar \cite{samymodeliar}]} \label{ext1s1s1}Let $k$ be a field of characteristic $p$, and $G$ be a finite $p$-group. Then
$$\Ext^1_{\comack_k(G)}(S_\un^G,S_\un^G)\cong \Hom_\Z\Big(G/\big(\Phi(G)I(G)\big),k^+\Big)\mvirg$$
where $k^+$ is the additive group of $k$, where $\Phi(G)$ is the Frattini subgroup of~$G$, and $I(G)$ is the subgroup of $G$ generated by elements of order 2 (so $I(G)=\un$ if $p>2$).
\end{mth}
\pf By Lemma~\ref{reso S_un}, there is an exact sequence
$$0\To k\to k\To \Hom_{\comack_k(G)}(FP_{\Omega_G},S_\un^G)\To \Ext^1_{\comack_k(G)}(S_\un^G,S_\un^G)\To 0\mvirg$$
hence $\Ext^1_{\comack_k(G)}(S_\un^G,S_\un^G)\cong \Hom_{\comack_k(G)}(FP_{\Omega_G},S_\un^G)$. \par
Now a morphism $\varphi:FP_{\Omega_G}\to S_\un^G$ is entirely determined by its evaluation at the trivial group, which is a morphism of $kG$-modules from $\Omega_G$ to $k$. Conversely, a morphism of $kG$-modules $f:\Omega_G\to k$ is the evaluation at $\un$ of a morphism of Mackey functors from $FP_{\Omega_G}$ to $S_\un^G$ if and only if it maps $(\Omega_G)^Q=r_\un^QFP_{\Omega_G}(Q)$ to $r_\un^QS_\un^Q(Q)=\zero$, for any non-trivial subgroup $Q$ of~$G$, or equivalently, for any subgroup $Q$ of order $p$ of $G$.\par
In other words $f$ is a $G$-invariant linear form on the space 
$$\Omega_G/\sumb{Q\leq G}{|Q|=p}{\rm Tr}_\un^QkG\mvirg$$
i.e. a linear form on $\Omega_G$ such that
\begin{equation}\label{eq1}
f\big(h(g-1)-(g-1)\big)=0,\;\;\forall g,h\in G\mvirg
\end{equation}
and
\begin{equation}\label{eq2}
f(1+x+\cdots+x^{p-1})=0,\;\;\forall x\in G,\;|x|=p\mpoint
\end{equation}
Now $\Omega_G$ is generated as a $k$-vector space by the elements $d_g=g-1$, for $g\in G$, and the only relations between these generators are $d_1=0$. Hence a linear form $f$ on $\Omega_G$ is determined by the values $u(g)=f(g-1)$, subject to $u(1)=0$. Since $d_{hg}=(hd_g-d_g)+d_h+d_g$, for any $g,h\in G$, Equation~\ref{eq1} is equivalent to
$$u(gh)=u(g)+u(h),\;\;\forall g,h\in G\mpoint$$
In other words $u$ is a group homomorphism from $G$ to $k^+$ (note that this implies $u(1)=0$, as required). Equivalently $u$ factors through a group homomorphism $G/\Phi(G)\to k^+$.\par
Now if $x$ is an element of order $p$ of $G$
$$1+x+\cdots+x^{p-1}=d_1+d_x+d_{x^2}+\cdots+d_{x^{p-1}}\mvirg$$
so Equation~2 is equivalent to
$$u(1)+u(x)+u(x^2)+\cdots +u(x^{p-1})=0,\;\;\forall x\in G,\;|x|=p\mpoint$$
If $u$ is a group homomorphism from $G$ to $k^+$, this is equivalent to
$$\big(0+1+2+\cdots+(p-1)\big)u(x)=\binom{p}{2}u(x)=0,\;\;\forall x\in G,\;|x|=p\mpoint$$
Now if $p$ is odd, the integer $\binom{p}{2}$ is a multiple of $p$, so this condition is satisfied whenever $u$ is a group homomorphism to $k^+$. And if $p=2$, this condition is equivalent to 
$$u(x)=0,\;\;\forall x\in G,\;|x|=2\mpoint$$
This shows finally that $\Hom_{\comack_k(G)}(FP_{\Omega_G},S_\un^G)$ is isomorphic to the group of homomorphisms from $G/\big(\Phi(G)I(G)\big)$ to $k^+$, and this completes the proof.\findemo
\begin{mth}{Lemma} \label{simple from elemab}Let $G$ be a finite $p$-group, and $Q$ be a subgroup of $G$. Set $B=N_G(Q)$, and denote by $A=\Phi(Q)$ the Frattini subgroup of $Q$. Then
$$S_Q^G\cong\Ind_{B}^G\,\iiota_{B/A}^{B}(S_{Q/A}^{B/A})\cong\mathsf{L}_{A\dom G}(S_{Q/A}^{B/A})\mvirg$$
where $A\dom G$ is viewed as a $(B/A,G)$-biset.
\end{mth}
\pf The second isomorphism follows from Lemma~\ref{ind iota}. As for the first one, recall (see~\cite{thevwebb}) that $S_Q^G=\Ind_{B}^GS_Q^{B}$, so it suffices to show that $S_Q^G=\iiota_{G/A}^G(S_{Q/A}^{G/A})$, when $Q\normal G$. In this case, if $H\leq G$ 
$$
\iiota_{G/A}^G(S_{Q/A}^{G/A})(H)=S_{Q/A}^{G/A}(HA/A)\mpoint$$
This is equal to $k$ if $HA=Q$, and to zero otherwise. But $HA=Q$ if and only if $H=Q$, and this completes the proof.\findemo
\begin{mth}{Corollary} \label{simple from elemab 2}With the same notation,
$$S_Q^G\cong\Ind_{B}^G\,\jiota_{B/A}^{B}(S_{Q/A}^{B/A})\cong\mathsf{R}_{G/A}(S_{Q/A}^{B/A})\mpoint$$
\end{mth}
\pf This follows from Lemma~\ref{simple from elemab}, using Proposition~\ref{composition}, and the fact that $S_Q^G$ is self dual.\findemo
\begin{mth}{Lemma} \label{rho simple}Let $G$ be a $p$-group, let $Q$ be a subgroup of $G$, and $(B,A)$ be a section of $G$. Then
$$\rho_{B/A}^B\Res_B^G(S_Q^G)\cong\dirsum{Q'}S_{Q'/A}^{B/A}\mvirg$$
where $Q'$ runs through the set of $G$-conjugates of $Q$ which contain $A$ and are contained in $B$, modulo $B$-conjugation.
\end{mth}
\pf Let $H/A$ be a subgroup of $B/A$. Then
$$\big(\rho_{B/A}^B\Res_B^G(S_Q^G)\big)(H/A)=\big(\Res_B^G(S_Q^G)\big)(H)=S_Q^G(H)\mvirg$$
and this is equal to $k$ if $H$ is conjugate to $Q$ in $G$, and to $\zero$ otherwise.\findemo 
\section{Some cohomological functors for $p$-groups}\label{some functors}
\begin{mth}{Notation} Let $G$ be a finite group. A subset $\mathcal{S}$ of the set of subgroups of $G$ will be called {\em convex} if
$$\forall H\leq K\leq L\leq G,\;\;H,L\in\mathcal{S}\;\Rightarrow\;K\in\mathcal{S}\mpoint$$
The set $\mathcal{S}$ will be called {\em $G$-stable} if it is invariant by $G$-conjugation.
\end{mth}
The following is an extension of Assertion~(ii) of Corollary 15.3 of~\cite{thevwebb}~:
\begin{mth}{Proposition} \label{subquotients of FP_k}Let $k$ be a field of characteristic $p$, and $G$ be a finite group. The correspondence
$$M\mapsto {\rm Supp}_p(M)=\{H\leq G\mid H\;\hbox{is a $p$-group and}\;M(H)\neq \zero\}$$
induces a one to one correspondence between the set of isomorphism classes of subquotients of the functor $FP_k$ and the set of $G$-stable convex subsets of the set of $p$-subgroups of $G$. The inverse bijection maps the $G$-stable convex subset $\mathcal{S}$ of $p$-subgroups of $G$ to the class of the functor $k_\mathcal{S}$ defined by
$$\forall H\leq G,\;\;k_\mathcal{S}(H)=\left\{\begin{array}{cl}k&\hbox{if}\;H_p\in \mathcal{S}\\\zero&\hbox{otherwise}\mvirg\end{array}\right.$$
where $H_p$ denotes a Sylow $p$-subgroup of $H$. \par
When $H\leq K\leq G$, the restriction map $r_H^K$ is equal to 0, except if $H_p\in \mathcal{S}$ and $K_p\in \mathcal{S}$, in which case $r_H^K=1$. The transfer map $t_H^K$ is multiplication by $|K:H|$ for all $H\leq K\leq G$. The conjugation map $c_{x,H}$ is always the identity map $k_\mathcal{S}(H)\to k_\mathcal{S}({^xH})$.\par
\end{mth}
\pf By Corollary~15.3 of~\cite{thevwebb}, the correspondences
$$F\mapsto {\rm Supp}_p(F)=\{Q\mid Q\;\hbox{is a $p$-subgroup of $G$ and }\;F(Q)\neq\zero\}$$
and
$$\mathcal{T}\mapsto k_{\mathcal{T}}={<}FP_k(Q){>}_{Q\in\mathcal{T}}$$
are mutual inverse bijections between the set of subfunctors of $FP_k$ and the set of $G$-stable subsets of the set $s_p(G)$ of $p$-subgroups of $G$, which are closed under taking subgroups. \par
By Proposition~2.4 of~\cite{thevwebb}, the value of $k_{\mathcal{T}}$ at some subgroup $H$ of $G$ is equal to
$$k_{\mathcal{T}}=\sumb{Q\in\mathcal{T}}{Q\leq H}|H:Q|k\mpoint$$
This is equal to $k$ if $H_p\in\mathcal{T}$, and to zero otherwise. The restriction maps between non zero values of $k_{\mathcal{T}}$ are equal to the identity map of $k$. The transfer map $t_H^K$ is multiplication by $|K:H|$, and the conjugation maps are identity maps (possibly zero).\par
It follows that any subquotient $M$ of $FP_k$ is equal to $k_\mathcal{T}/k_{\mathcal{T}'}$, where $\mathcal{T}$ and~$\mathcal{T}'$ are $G$-stable subsets of $s_p(G)$, which are closed under taking subgroups, and such that $\mathcal{T}'\subseteq\mathcal{T}$. With the notation of Proposition~\ref{subquotients of FP_k}, this means that $M\cong k_{\mathcal{S}}$, where $\mathcal{S}=\mathcal{T}-\mathcal{T}'$ is a $G$-stable convex subset of $s_p(G)$, equal to ${\rm Supp}_p(M)$.\par
Conversely, let $\mathcal{S}$ be a $G$-stable convex subset of $s_p(G)$. Set
$$\mathcal{T}=\{Q\in s_p(G)\mid \exists S\in\mathcal{S},\;Q\leq S\}$$
$$\mathcal{T}'=\mathcal{T}-\{Q\in s_p(G)\mid \exists S\in\mathcal{S},\;S\leq Q\}\mpoint$$
Then  $\mathcal{T}$ and $\mathcal{T}'$ are $G$-stable, and closed under taking subgroups. Moreover
$$\mathcal{T}-\mathcal{T}'=\{Q\in s_p(G)\mid \exists S, S'\in\mathcal{S},\;\;S'\leq Q\leq S\}\mpoint$$
Thus $\mathcal{T}-\mathcal{T}'=\mathcal{S}$, and $k_\mathcal{S}=k_\mathcal{T}/k_{\mathcal{T}'}$ is a subquotient of $FP_k$.\par
These correspondences $M\mapsto {\rm Supp}_p(M)$ and $\mathcal{S}\mapsto k_\mathcal{S}$ are clearly mutual inverse bijections between the set of isomorphism classes of subquotients of~$FP_k$ and the set of $G$-stable convex subsets of $s_p(G)$.\findemo
\begin{mth}{Corollary} \label{SigmaQR}Let $G$ be a finite $p$-group, and $Q\leq R$ be subgroups of~$G$. Then there exists a unique object $\Sigma_{Q,R}^G$ of $\comack_k(G)$ such that
$$\forall H\leq G,\;\;\Sigma_{Q,R}^G(H)=\left\{\begin{array}{cl}k&\hbox{if}\;\;Q\leq_G H\leq_G R\\\zero&\hbox{otherwise}\mpoint\end{array}\right.$$
and such that 
$$\forall H<K\leq G,\;\;t_H^K=0,\;\;r_H^K=\left\{\begin{array}{cl}1&\hbox{if}\;Q\leq_G H\:\hbox{and}\;K\leq_GR\\0&\hbox{otherwise}\mpoint\end{array}\right.$$
Moreover the socle of $\Sigma_{Q,R}^G$ is isomorphic to $S_Q^G$, and the head of $\Sigma_{Q,R}^G$ is isomorphic to $S_R^G$.
\end{mth}
\pf Let $\mathcal{S}=\{S\leq G\mid Q\leq_GS\leq_GR\}$. Then $\mathcal{S}$ is a $G$-stable convex set of subgroups of $G$, and the corresponding functor $k_\mathcal{S}$ fulfills the required conditions, so the functor $\Sigma_{Q,R}^G$ exists. The uniqueness follows from the fact that the values of $\Sigma_{Q,R}^G$ are given, as well as the transfer and restriction maps. The non zero conjugation maps $c_{x,Q}$ are determined by elements of $H^1(G,\Ind_Q^Gk^\times)\cong\Hom_\Z(Q,k^\times)$, which is equal to zero since $G$ is a $p$-group. So the conjugation maps are all identity maps (possibly zero).\par
Now a morphism $\varphi$ from $\Sigma_{Q,R}^G$ to some functor $N$ is entirely determined by its evaluation $\varphi_R:\Sigma_{Q,R}^G(R)=k\to N(R)$~: indeed, if $S\leq G$ and $\Sigma_{Q,R}^G(S)\neq\zero$, then $S^g\leq R$ for some $g\in G$, and then the map $\psi=c_{g,S^g}r_{S^g}^Q$ is an isomorphism from $\Sigma_{Q,R}^G(R) $ to $\Sigma_{Q,R}^G(S)$, such that $\varphi_Q=c_{g,S^g}r_{S^g}^Q\varphi_R\psi^{-1}$. Hence the only simple quotient of $\Sigma_{Q,R}^G$ is $S_R^G$, with multiplicity one.\par
By a similar argument, a morphism from $N$ to $\Sigma_{Q,R}^G$ is determined by its value at $Q$, so the only simple subfunctor of $\Sigma_{Q,R}^G$ is $S_Q^G$, with multiplicity one.\findemo
\begin{mth}{Notation} \label{gamma_X} Let $G$ be a finite $p$-group. If $Q<R$ are subgroups of $G$ such that $|R:Q|=p$, the functor $\Sigma_{Q,R}^G$ will be denoted by $\bin{R}{Q}{G}$. The dual functor $\big(\Sigma_{Q,R}^G\big)^*$ will be denoted by $\bin{Q}{R}{G}$.
With this notation, there are non split exact sequences in $\comack_k(G)$
\begin{equation}\label{EQR}
D_{Q,R}:0\To S_Q^G\To\bin{R}{Q}{G}\To S_R^G\To 0\mpoint
\end{equation}
\vspace{-3ex}
\begin{equation}\label{EQR*}
D_{Q,R}^*:0\To S_R^G\To\bin{Q}{R}{G}\To S_Q^G\To 0\mpoint
\end{equation}
In particular, if $X$ is a subgroup of order $p$ of $G$, set $D_X=D_{\un,X}$ and $D_X^*=D_{\un,X}^*$, and denote by $\gamma_X^G$ (or $\gamma_X$ if $G$ is clear from the context) the element of $\Ext^2_{\comack_k(G)}(S_\un^G,S_\un^G)$ represented by the exact sequence
$$\Gamma_X:\;\;0\To S_\un^G\To\bin{X}{\un}{G}\To \bin{\un}{X}{G}\To S_\un^G\To 0\mvirg$$
obtained by splicing the sequences $D_X$ and $D_X^*$.
\end{mth}
\section{Extensions of simple functors for $p$-groups}\label{extensions}
\begin{mth}{Lemma} \label{D_Z}Let $G$ be a $p$-group, and $Z$ be a central subgroup of order~$p$ of $G$. Then there are isomorphisms
$$\iiota_{G/Z}^G(S_\un^{G/Z})\cong \bin{Z}{\un}{G}\mvirg \;\;\jiota_{G/Z}^G(S_\un^{G/Z})\cong \bin{\un}{Z}{G}$$
in $\comack_k(G)$. In particular, there are non-split exact sequences in $\comack_k(G)$
$$D_Z: \;\;0\To S_\un^G\To\iiota_{G/Z}^G(S_{\un}^{G/Z})\To S_Z^G\To 0\mpoint$$
$$D_Z^*: \;\;0\To S_Z^G\To\jiota_{G/Z}^G(S_{\un}^{G/Z})\To S_\un^G\To 0\mpoint$$
\end{mth}
\pf Let $H$ be a subgroup of $G$. Then $\iiota_{G/Z}^G(S_{\un}^{G/Z})(H)=S_{\un}^{G/Z}(HZ/Z)$ is equal to $k$ if $HZ=Z$, i.e. if $H\leq Z$, and to zero otherwise. The conjugation maps are all identity maps (possibly zero), the restriction map $r_\un^Z$ is the identity map of $k$, and the transfer map $t_\un^Z$ is zero. Thus $\iiota_{G/Z}^G(S_\un^{G/Z})\cong \bin{Z}{\un}{G}$. The other isomorphism follows by duality, and the two exact sequences are special cases of Sequences~\ref{EQR} and~\ref{EQR*}.\findemo
\begin{mth}{Corollary} \label{central p}For any $n\in\N$, composition by $D_Z$ induces a group isomorphism 
$$\Ext^{n-1}_{\comack_k(G)}(S_\un^G,S_\un^G)\cong \Ext^{n}_{\comack_k(G)}(S_Z^G,S_\un^G)$$
\end{mth}
\pf Indeed, for $n\in\N$, 
$$\Ext^n_{\comack_k(G)}\big(\bin{Z}{\un}{G},S_\un^G\big)\cong \Ext^n_{\comack_k(G/Z)}\big(S_\un^{G/Z},\rho_{G/Z}^G(S_\un^G)\big)=\zero$$
by Proposition~\ref{ext adjunction} and Lemma~\ref{rho simple}.\findemo
\begin{mth}{Corollary} \label{ext1sxs1}Let $k$ be a field of characteristic $p$, and $G$ be a $p$-group. If $X$ is a subgroup of order $p$ of $G$, then
$$\Ext_{\comack_k(G)}^1(S_X^G,S_\un^G)\cong k\cong \Ext_{\comack_k(G)}^1(S_\un^G,S_X^G)\mpoint$$
\vspace{-1ex}
\nopagebreak
\end{mth}
\nopagebreak
\pf Indeed $N_G(X)$ is equal to the centralizer $C$ of $X$, since $|X|=p$, and $S_X^G=\Ind_{C}^GS_{X}^{C}$. Thus for any $n\in\N$
\begin{eqnarray*}
\Ext_{\comack_k(G)}^n(S_X^G,S_\un^G)&\cong& \Ext_{\comack_k(C)}^n(S_X^C,\Res_C^GS_\un^G)\\
&\cong&\Ext_{\comack_k(C)}^n(S_X^C,S_\un^C)\\
&\cong&\Ext_{\comack_k(C)}^{n-1}(S_\un^C,S_\un^C)
\end{eqnarray*}
by Corollary~\ref{central p}, since $X$ is a central subgroup of order $p$ of $C$. Corollary~\ref{ext1sxs1} follows, taking $n=1$.\par
The isomorphism $k\cong \Ext_{\comack_k(G)}^1(S_\un^G,S_X^G)$ follows from Remark~\ref{ext dual}, since the simple functors are self dual. \findemo
\begin{mth}{Notation} If $N$ is a normal subgroup of a group $G$, contained in the subgroup $H$ of $G$, denote by $K_H(N)$ the set of complements of $N$ in $H$, i.e. the set of subgroups $X$ of $G$ such that $NX=H$ and $N\cap X=\un$. \par
Denote by $[N_G(H)\dom K_H(N)]$ a set of representatives of $N_G(H)$-conjugacy classes of subgroups in $K_H(N)$. 
\end{mth}
The following proposition is a generalization of Lemma~\ref{D_Z}~:
\begin{mth}{Proposition} \label{iota S_Q}Let $G$ be a $p$-group, let $Q$ be a subgroup of $G$, and~$Z$ be a central subgroup of order~$p$ of $G$ contained in $Q$. Then there are non-split exact sequences in $\comack_k(G)$
$$0\To \dirsum{X\in [N_G(Q)\dom K_Q(Z)]}S_{X}^G\To\iiota_{G/Z}^G(S_{Q/Z}^{G/Z})\To S_Q^G\To 0\mvirg$$
$$0\To S_Q^G\To\jiota_{G/Z}^G(S_{Q/Z}^{G/Z})\To \dirsum{X\in [N_G(Q)\dom K_Q(Z)]}S_{X}^G\To 0\mpoint$$
\end{mth}
\pf Since the simple functors are self dual, the second sequence is obtained from the first by applying duality, by definition of the functor $\jiota_{G/Z}^G$ and Assertion~6 of Proposition~\ref{composition}. So it suffices to prove the existence of the first one.\par
Denote by $I$ the functor $\iiota_{G/Z}^G(S_{Q/Z}^{G/Z})$, and by $R$ its radical. If $H$ is a subgroup of $G$, then $I(H)=S_{Q/Z}^{G/Z}(HZ/Z)$ is equal to $k$ if $HZ=_GQ$, and to zero otherwise. Moreover $HZ=Q$ if and only if $H=Q$, or $H\in K_Q(Z)$, since $Z$ has order $p$. Equation~\ref{mult} shows that the composition factors of $I$ are the functors $S_Q^G$ and the functors $S_X^G$, for $X\in [N_G(Q)\dom K_Q(Z)]$, each with multiplicity~1.\par
Moreover, if $Y$ is any subgroup of $G$, consider
$$\mathcal{H}_Y=\Hom_{\comack_k(G)}\big(I,S_Y^G\big)\cong \Hom_{\comack_k(G/Z)}\big(S_{Q/Z}^{G/Z},\rho_{G/Z}^G(S_Y^G)\big)\mpoint$$
Now $\rho_{G/Z}^G(S_Y^G)=\zero$ if $Z\not\leq Y$, and $\rho_{G/Z}^G(S_Y^G)=S_{Y/Z}^{G/Z}$, by Lemma~\ref{rho simple}. Thus $\mathcal{H}_Y=\zero$ unless $Y=_GQ$, and $\mathcal{H}_Y\cong k$ in this case.\par
This means that $I/R$ is simple, isomorphic to $S_Q^G$. Moreover, if $H$ is a subgroup of $G$, then $R(H)=\zero$ except if $HZ=_GQ$ and $H\cap Z=\un$, i.e. if~$H$ is conjugate to some element of $K_Q(Z)$ in $G$, and $R(H)\cong k$ in this case. Hence $R$ is isomorphic to the direct sum of the simple functors $S_H^G$, where $H$ runs in a set of representatives of $G$-conjugacy classes of such subgroups, i.e. equivalently $H\in [N_G(Q)\dom K_Q(Z)]$.\findemo
\begin{mth}{Theorem} \label{quotient by central}Let $k$ be a field of characteristic $p>0$, let $G$ be a finite $p$-group. Let $Q$ and $R$ be subgroups of $G$, and let $Z$ be a subgroup of order $p$ of $Q\cap R\cap Z(G)$. Set $\sur{G}=G/Z$, $\sur{Q}=Q/Z$, $\sur{R}=R/Z$, and $\mathcal{K}=[N_G(Q)\dom K_Q(Z)]$. \par
Then, for any $n\in\N$, there are short exact sequences of vector spaces
$$0\to \dirsum{X\in\mathcal{K}}\Ext^{n-1}_{\comack_k(G)}(S_X^G,S_R^G)\to \Ext^{n}_{\comack_k(G)}(S_Q^G,S_R^G)\stackrel{\pi_n}{\to} \Ext^{n}_{\comack_k(\sur{G})}(S_{\sur{Q}}^{\sur{G}},S_{\sur{R}}^{\sur{G}})\to 0\mvirg$$
where the map $\pi_n$ is induced by $\rho_{G/Z}^G$.
\end{mth}
\pf Applying $\Hom_{\comack_k(G)}({-},S_R^G)$ to the first exact sequence of Proposition~\ref{iota S_Q} gives a long exact sequence
$$\xymatrix@C=3.5ex{
\cdots\ar[r]\ar[r]& *!U(0.3){\dirsum{X\in\mathcal{K}}\Ext^{n-1}_{\comack_k(G)}(S_X^G,S_R^G)}\ar[r]& *!U(0.1){\Ext^{n}_{\comack_k(G)}(S_Q^G,S_R^G)}\ar[r]^-{\pi_n}& *!U(0.1){\Ext^{n}_{\comack_k(\sur{G})}(S_{\sur{Q}}^{\sur{G}},S_{\sur{R}}^{\sur{G}})}\ar`r[d]`[l]`[dlll]`[dll][dll]\\
& *!U(0.3){\dirsum{X\in\mathcal{K}}\Ext^{n}_{\comack_k(G)}(S_X^G,S_R^G)}\ar[r]& *!U(0.1){\Ext^{n+1}_{\comack_k(G)}(S_Q^G,S_R^G)}\ar[r]^-{\pi_{n+1}}&\cdots\rule{1.5cm}{0cm}
}
$$
where the image of an extension $u\in \Ext^n_{\comack_k(G)}(S_Q^G,S_R^G)$ under the map $\pi_n$ is obtained by first taking the image under the map 
$$\varphi:\iiota_{G/Z}^G(S_{Q/Z}^{G/Z})\to S_Q^G\in\Ext^0_{\comack_k(G)}\big(\iiota_{G/Z}^G(S_{Q/Z}^{G/Z}),S_Q^G\big)\mvirg$$
and then using the adjunction isomorphisms
$$\Ext^n_{\comack_k(G)}\big(\iiota_{G/Z}^G(S_{Q/Z}^{G/Z}),S_R^G\big)\cong\Ext^{n}_{\comack_k(\sur{G})}\big(S_{\sur{Q}}^{\sur{G}},\rho_{G/Z}^G(S_R^G)\big)\cong\Ext^{n}_{\comack_k(\sur{G})}(S_{\sur{Q}}^{\sur{G}},S_{\sur{R}}^{\sur{G}})\mpoint$$
In other words, with the notation of Remark~\ref{Yoneda compatible}
$$\pi_n(u)=\alpha_n(u\circ\varphi)=\rho_{G/Z}^G(u)\circ \alpha_0(\varphi)=\rho_{G/Z}^G(u)\mvirg$$
since the map $\alpha_0(\varphi):S_{Q/Z}^{G/Z}\to \rho_{G/Z}^G(S_Q^G)\cong S_{Q/Z}^{G/Z}$ obtained from $\varphi$ by adjunction is the identity map. This shows that $\pi_n$ is induced by $\rho_{G/Z}^G$.\par
Now the inflation functor $\Inf_{G/Z}^G:\fun_k(G/Z)\to\fun_k(G)$ is an exact functor such that $\rho_{G/Z}^G\circ\Inf_{G/Z}^G$ is isomorphic to the identity functor. Since $S_Q^Z\cong\Inf_{G/Z}^G(S_{Q/Z}^{G/Z})$ and $S_R^Z\cong\Inf_{G/Z}^G(S_{R/Z}^{G/Z})$, this inflation functor induces a map
$$\sigma_n:\Ext^{n}_{\comack_k(\sur{G})}(S_{\sur{Q}}^{\sur{G}},S_{\sur{R}}^{\sur{G}})\to \Ext^{n}_{\comack_k(G)}(S_{Q}^{G},S_{R}^{G})$$
such that $\pi_n\circ\sigma_n=\Id$. In particular $\pi_n$ is surjective, so the long exact sequence above splits as a series of short exact sequences.\findemo

\begin{mth}{Proposition} \label{index p}Let $k$ be a field of characteristic $p>0$, let $G$ be a finite $p$-group, and $H$ be a subgroup of index~$p$ in $G$. Set $I=\Ind_H^GS_\un^H$, and let $R$ and $S$ denote respectively the radical and the socle of~$I$ as an object of $\comack_k(G)$. Then $I\supset R\supseteq S\supset \zero$. Moreover $I/R\cong S\cong S_\un^G$, and
$$R/S\cong L\oplus\dirsum{X\in[G\dom K_G(H)]}S_X^G\mvirg$$
where $L$ is a functor all of whose composition factors are isomorphic to $S_\un^G$, with multiplicity $p-2$.
\end{mth}
\pf Let $Q$ be a subgroup of $G$. Then by Proposition~\ref{ext adjunction}
$$\Hom_{\comack_k(G)}(I,S_Q^G)\cong\Hom_{\comack_k(H)}(S_\un^H,\Res_H^GS_Q^G)\mpoint$$
Moreover, by Lemma~\ref{rho simple}, the functor $\Res_H^GS_Q^G$ is isomorphic to the direct sum of functors $S_{Q'}^H$, where $Q'$ runs through a set of $H$-conjugacy classes of $G$-conjugates of $Q$ which are contained in $H$. \par
Thus $\Hom_{\comack_k(G)}(I,S_Q^G)=\zero$ if $Q\neq\un$, and $\Hom_{\comack_k(G)}(I,S_\un^G)\cong k$. This shows that $I/R\cong S_\un^G$. \par
A similar argument shows that $\Hom_{\comack_k(G)}(S_Q^G,I)$ is equal to $\zero$ if $Q\neq\nolinebreak\un$, and isomorphic to $k$ if $Q=\un$. Thus $S\cong S_\un^G$.\par
Let $X$ be a subgroup of $G$. Then
$$I(X)\cong\dirsum{x\in G/HX}S_\un^H(H\cap{^xX})\mpoint$$
This is isomorphic to $k^p$ if $X=\un$, to $k$ if $X$ is non trivial and $H\cap X=\un$ (which implies that $HX=G$), and to $\zero$ otherwise. In particular $I\neq S$ (since $I(\un)=k^p\neq k=S(\un)$), so $S\subseteq R$.\par
Moreover, the short exact sequence 
$$0\To R\To I\to S_\un^G\To 0$$
shows that $\Ext_{\comack_k(G)}^n(S_X^G,R)\cong \Ext_{\comack_k(G)}^{n-1}(S_X^G,S_\un^G))$ for any $n\in \N$, since $\Ext_{\comack_k(G)}^n(S_X^G,I)\cong \Ext_{\comack_k(H)}^n(\Res_H^GS_X^G,S_\un^H)$, and since $\Res_H^GS_X^G=\zero$. In particular $\Ext_{\comack_k(G)}^1(S_X^G,R)\cong \Hom_{\comack_k(G)}(S_X^G,S_\un^G)=\zero$.\par
Now the short exact sequence
$$0\To S_\un^G\To R\to R/S\To 0$$
yields a long exact sequence of $\Ext_{\comack_k(G)}(S_X^G,-)$ groups starting with
$$0\To 0\To 0\To \Hom(S_X,R/S)\to \Ext^1(S_X,S_\un)\to \Ext^1(S_X,R)=\zero\mpoint$$
It follows that $\Hom(S_X,R/S)\cong k$, by Corollary~\ref{ext1sxs1}. Moreover $I$ is self dual, since $S_\un^H$ is self dual and $\Ind_H^G$ commutes with duality, by Proposition~\ref{composition}. Thus 
$$\Hom(R/S,S_X)\cong \Hom(S_X,R/S)\cong k\mvirg$$
by Remark~\ref{ext dual}, since $S_X$ is also self dual.\par
In particular, there exists an injection $i:S_X\To R/S$ and a surjection $s:R/S\To S_X$. If $s\circ i=0$, then $(R/S)(X)$ has dimension at least equal to~2. But $(R/S)(X)\cong I(X)\cong k$, so $s\circ i\neq 0$. Hence $s\circ i$ is invertible, and $S_X$ is a direct summand of $R/S$.\par
It follows that $\dirsum{X\in[G\dom K_G(H)]}S_X^G$ is also a direct summand of $R/S$, for it is a direct sum of non isomorphic simple summands of $R/S$. Hence there exists a subfunctor $L$ of $R/S$ such that
$$R/S\cong L\oplus \dirsum{X\in[G\dom K_G(H)]}S_X^G\mpoint$$
Now if $Q$ is a subgroup of $G$, then $(R/S)(Q)$ is equal to $\zero$, unless $Q=\un$ or~$Q$ is a complement of $H$ in $G$, in which case $(R/S)(Q)\cong k$. It follows that $L(Q)$ is equal to $\zero$, except if $Q=\un$ (and $L(\un)$ has dimension $p-2$ over $k$). This completes the proof.\findemo
\section{The case of cyclic $p$-groups}\label{cyclic p-groups}
\npar It has been shown by M. Samy Modeliar (\cite{samymodeliar}) that cyclic $p$-groups are poco groups over a field $k$ of characteristic~$p$. This result relies on the construction of periodic projective resolutions for the fixed points functors, which can be seen has follows~:\pagebreak[3] if $G$ is cyclic of order $p^m$, then the group algebra $kG$ is isomorphic to a truncated polynomial algebra $A=k[X]/(X^{p^m})$, via the map sending $X\in A$ to $g-1\in kG$, where $g$ is some chosen generator of~$G$. \par
The indecomposable $A$-modules are the modules $E_d=k[X]/(X^d)$, for $1\leq d\leq p^m$, so $E_d$, where the chosen generator $g$ of $G$ acts by multiplication by $1+X$, is the unique indecomposable $kG$-module of dimension~$d$, up to isomorphism. The indecomposable permutation $kG$-modules have dimension equal to a power of $p$, so they are the modules $E_{p^j}$, for $0\leq j\leq m$.\par
This means that the functors $FP_{E_{p^j}}$ are the projective indecomposable objects, in $\comack_k(G)$, and they are their own projective resolution.\par
Now if $1\leq d\leq p^m$ and $d$ is not a power of $p$, there exists a unique integer~$h$ in $\{0,\cdots,m-1\}$ such that $p^h<d<p^{h+1}$. Let 
$$e_d:E_{p^h}\oplus E_{p^{h+1}}\To E_d$$
denote the morphism of $A$-modules induced by the map 
$$(Q,R)\in k[X]\times k[X]\mapsto X^{d-p^h}Q+R\in E_d\mpoint$$
This map is well defined, and surjective, and it is a morphism of $kG$-modules.\par
The unique subgroup $G_l$ of $G$ of index $p^l$, where $0\leq l\leq m$, is the subgroup generated by $g^{p^l}$. If $W$ is a $kG$-module, then the subspace $W^{G_l}$ of~$W$ on which $G_l$ acts trivially is equal to
$$W^{G_l}=\{w\in W\mid g^{p^l}w=w\}=\{w\in W\mid (g-1)^{p^l}w=0\}\mpoint$$
So viewing $W$ as an $A$-module, the space $W^{G_l}$ is equal to the kernel of $X^{p^l}$ on $W$. In particular, for any $j\in\{1,\cdots,m\}$, the space $(E_j)^{G_l}$ is equal to $k[X]X^r/(X^j)$, where $r=\max(j-p^l,0)$. This shows that
$$e_d\big((E_{p^h}\oplus E_{p^{h+1}})^{G_l}\big)=\big(k[X]X^{d-p^h+\max(0,p^h-p^l)}+k[X]X^{\max(0,p^{h+1}-p^l)}\big)/(X^d)\mpoint$$
If $l\leq h$, this is equal to $k[X]X^{d-p^l}/(X^d)$, and if $l>h$, this is equal to $k[X]/(X^d)$. In other words
$$e_d\big((E_{p^h}\oplus E_{p^{h+1}})^{G_l}\big)=k[X]X^{\max(0,d-p^l)}/(X^d)=(E_d)^{G_l}\mpoint$$
It follows that the restriction of the map $e_d$ to the spaces of fixed points by any subgroup of $G$ is surjective~: in the terminology of Samy Modeliar~(\cite{samymodeliar}), the map $e_d$ is {\em supersurjective}. Equivalently, it induces a surjection
$$FP_{e_d}:FP_{E_{p^h}}\oplus FP_{E_{p^{h+1}}}\To FP_{E_d}$$
in $\comack_k(G)$.\par
Now the kernel of $e_d$ consists of the images in $E_{p^h}\oplus E_{p^{h+1}}$ of the pairs $(Q,R)\in k[X]\times k[X]$ such that $X^{d-p^h}Q+R$ is a multiple of $X^d$, i.e. of the pairs $(-T,X^{d-p^h}T)$, for $T\in k[X]$. Hence the map
$$T\in E_{p^h+p^{h+1}-d}\mapsto (-T,X^{d-p^h}T)\in\Ker\,e_d$$
is well defined (note that $p^h<p^h+p^{h+1}-d<p^{h+1}$), and it is an isomorphism of $A$-modules. This yields a short exact sequence
$$0\To E_{p^h+p^{h+1}-d}\To E_{p^h}\oplus E_{p^{h+1}}\stackrel{e_d}{\To}E_d\To 0$$
of $kG$-modules, leading to a short exact sequence
$$0\To FP_{E_{p^h+p^{h+1}-d}}\To FP_{E_{p^h}}\oplus FP_{E_{p^{h+1}}}\stackrel{FP_{e_d}}{\To}FP_{E_d}\To 0$$
in $\comack_k(G)$. This shows in particular that $FP_{E_d}$ has a projective resolution in $\comack_k(G)$, which is periodic of period 2, of the form
$$\cdots \stackrel{FP_{\alpha_d}}{\To}FP_{E_{p^h}\oplus E_{p^{h+1}}}\stackrel{FP_{\beta_d}}{\To}FP_{E_{p^h}\oplus E_{p^{h+1}}}\stackrel{FP_{\alpha_d}}{\To}FP_{E_{p^h}\oplus E_{p^{h+1}}}\stackrel{FP_{e_d}}{\To}FP_{E_d}\To 0$$
where $\alpha_d$ and $\beta_d$ are the endomorphisms of $E_{p^h}\oplus E_{p^{h+1}}$ given by
\begin{eqnarray*}
\alpha_d(Q,R)&=&(-X^{p^{h+1}-d}Q-R,X^{p^{h+1}-p^h}Q+X^{d-p^h}R)\\
\beta_d(Q,R)&=&(-X^{d-p^{h}}Q-R,X^{p^{h+1}-p^h}Q+X^{p^{h+1}-d}R)\mpoint\\
\end{eqnarray*}
Hence all the functors $FP_V$, where $V$ is an indecomposable $kG$-module, have a periodic resolution. In particular $G$ is a poco group, by Proposition~\ref{FP_V is enough}.\par
\npar If $p^m\geq 3$, this applies to the case $d=p^m-1$ (which is not a power of~$p$), and $h=m-1$, consequently. In this case $E_d\cong \Omega_G$, and $E_{p^{h+1}}\cong kG$. Splicing the above resolution with the short exact sequence
$$0\To FP_{E_{p^m-1}}\To FP_{E_{p^m}}\to S_\un^G\To 0$$
of Lemma~\ref{reso S_un} gives a projective resolution of $S_\un^G$ of the form
$$\cdots \stackrel{FP_{\beta}}{\To}FP_{E_{p^{m-1}}\oplus E_{p^{m}}}\stackrel{FP_{\alpha}}{\To}FP_{E_{p^{m-1}}\oplus E_{p^{m}}}\stackrel{FP_{\gamma}}{\To}FP_{E_{p^m}}\To S_\un^G\To 0\mvirg$$
where $\alpha$ and $\beta$ are the above maps $\alpha_d$ and $\beta_d$, for $d=p^m-1$ (and $h=m-1$), and $\gamma:E_{p^{m-1}}\oplus E_{p^m}\to E_{p^m}$ is the map defined by
$$\gamma(Q,R)=X^{p^m-p^{m-1}}Q+XR\mpoint$$
Now $FP_{E_{p^{m-1}}\oplus E_{p^m}}\cong FP_{E_{p^{m-1}}}\oplus FP_{E_{p^m}}$. Moreover, by Equation~\ref{mult},
$$\Hom_{\comack_k(G)}(FP_{E_{p^{m-1}}},S_\un^G)\cong S_\un^G(G_{p^{m-1}})=\zero\mvirg$$
and 
$$\Hom_{\comack_k(G)}(FP_{E_{p^m}},S_\un^G)\cong S_\un^G(G_{p^m})=S_\un^G(\un)=k\mpoint$$
More precisely, a morphism $FP_{E_{p^m}}\To S_\un^G$ is determined by its evaluation at the trivial subgroup, which is a scalar multiple of the augmentation morphism~$\varepsilon$.\par
It follows easily that the groups $\Ext^n_{\comack_k(G)}(S_\un^G,S_\un^G)$ are the cohomology groups of the complex whose terms are all isomorphic to $k$, with zero differentials. Thus $\Ext^n_{\comack_k(G)}(S_\un^G,S_\un^G)\cong k$, for $n\in\N$.
\begin{mth}{Proposition} \label{ext cyclic}Let $k$ be a field of characteristic $p$, and $G$ be a cyclic $p$-group. Then~:
$$\forall n\in \N-\{0\},\;\;\Ext_{\comack_k(G)}^n(S_\un^G,S_\un^G)\cong\left\{\begin{array}{cl}k&\hbox{if}\;|G|\geq 3\\\zero&\hbox{if}\;|G|\leq 2\end{array}\right.\mpoint$$
\end{mth}
\pf The case $|G|\geq 3$ follows from the above discussion. The case $|G|\leq 2$ follows essentially from the fact that in this case, any $kG$-module is a permutation module~: if $G$ is trivial, there is nothing to prove. And if $G$ has order 2, then $S_\un^G$ has a projective resolution
$$0\To FP_k\To FP_{kG}\To S_\un^G\To 0\mvirg$$
and moreover $\Hom_{\comack_k(G)}(FP_k,S_\un^G)=\zero$, by Equation~\ref{mult}.\findemo
\section{The case of elementary abelian $p$-groups}\label{elemab}
\begin{mth}{Proposition} \label{index p elemab}Let $k$ be a field of characteristic $p$, and $G$ be an elementary abelian $p$-group. Let $H$ be a subgroup of index $p$ in $G$, and $T\in K_G(H)$. Then the functor $I=\Ind_H^GS_\un^H$ of $\comack_k(G)$ has a subfunctor~$J$ isomorphic to $\iota_{G/T}^G(S_\un^{G/T})$. Moreover $J$ is contained in the radical $R$ of $I$, and there is an isomorphism
$$R/J\cong L\oplus \dirsum{X\in K_G(H)-\{T\}}S_X^G\mvirg$$
where $L$ is a functor all of whose composition factors are isomorphic to $S_\un^G$, with multiplicity $p-2$.
\end{mth}
\pf Let $S$ be the socle of $I$. By Proposition~\ref{index p}, there is a filtration $I\supset R\supseteq S\supset \zero$, and $I/R\cong S\cong S_\un^G$. Moreover
\begin{equation}\label{RoverS}
R/S\cong L\oplus \dirsum{X\in K_G(H)}S_X^G\mvirg
\end{equation}
where $L$ is a functor all of whose composition factors are isomorphic to $S_\un^G$, with multiplicity $p-2$.\par
By Lemma~\ref{central p}, the functor $J=\iiota_{G/T}^G(S_\un^{G/T})$ has simple socle, equal to its radical, and isomorphic to $S_\un^G$, and simple top, isomorphic to $S_T^G$. Moreover
$$\Hom_{\comack_k(G)}(J,I)\cong \Hom_{\comack_k(H)}(\Res_H^GJ,S_\un^H)\cong k\mvirg$$
since $\Res_H^GJ\cong S_\un^H$, for the only non zero evaluation of this functor is at the trivial group, where it is equal to $k$.\par
Hence there is a non zero morphism $f:J\to I$. If $f$ is not injective, then its image is isomorphic to $S_T^G$, since this is the only proper non zero quotient of $J$. But $S_T^G$ is not isomorphic to a subfunctor of $I$, since the socle of $I$ is simple and isomorphic to $S_\un^G$.\par
It follows that $f$ is injective, and one can identify $J$ with a subfunctor of~$I$, containing $S$, and $J/S\cong S_T^G$. Moreover $J$ is a proper subfunctor of $I$, since $S_\un^G$ is not a quotient of $J$. Hence $J\subseteq R$, and
$$R/J\cong L\oplus \dirsum{X\in K_G(H)-\{T\}}S_X^G\mvirg$$
by (\ref{RoverS}).\findemo
\begin{mth}{Corollary} \label{long exact sequence}With the same notation, there is a long exact sequence of extension groups
$$\xymatrix{
\cdots\ar[r]&*!U(0.6){L(n-1)\oplus\dirsum{X\in\mathcal{X}}E_G(n-2)}\ar[r]&E_G(n)\ar[r]&E_H(n)\ar`r[d]`[l]`[dlll]`[dll][dll]\\
&*!U(0.6){L(n)\oplus\dirsum{X\in \mathcal{X}}E_G(n-1)}\ar[r]&E_G(n+1)\ar[r]&E_H(n+1) \cdots\mvirg\\
}
$$
where $E_G(n)=\Ext^n_{\comack_k(G)}(S_\un^G,S_\un^G)$, $E_H(n)=\Ext^n_{\comack_k(H)}(S_\un^H,S_\un^H)$, $L(n)=\Ext^n_{\comack_k(G)}(L,S_\un^G)$, and $\mathcal{X}=K_G(H)-\{T\}$.
\end{mth}
\pf There is a short exact sequence
$$0\To J\To I\To I/J\To 0$$
in $\comack_k(G)$.
Since $\Ext_{\comack_k(G)}^n(J,S_\un^G)\cong\Ext_{\comack_k(G/T)}^n\big(S_\un^{G/T},\rho_{G/T}^G(S_\un^G)\big)=\zero$ for any $n\in\nolinebreak\N$, by Proposition~\ref{ext adjunction} and Lemma~\ref{rho simple}, and since $\Res_H^GS_\un^G=S_\un^H$, it follows that
$$\Ext^n_{\comack_k(G)}(I/J,S_\un^G)\cong \Ext^n_{\comack_k(G)}(I,S_\un^G)\cong \Ext^n_{\comack_k(H)}(S_\un^H,S_\un^H)=E_H(n)\mpoint$$
Now by Proposition~\ref{index p elemab}, there is a short exact sequence
$$0\To R/J\To I/J\To S_\un^G\To 0\mvirg$$
and $R/J\cong L\oplus \dirsum{X\in K_G(H)-\{T\}} S_X^G$. Applying the functor $\Hom_{\comack_k(G)}({-},S_\un^G)$ gives the long exact sequence of extension groups in $\comack_k(G)$
$$\xymatrix{
0\ar[r]&\Hom(S_\un^G,S_\un^G)\ar[r]&\Hom(I/J,S_\un^G)\ar[r]&\Hom(R/J,S_\un^G)\ar`r[d]`[l]`[dlll]`[dll][dll]\\
&\Ext^1(S_\un^G,S_\un^G)\ar[r]&\Ext^1(I/J,S_\un^G)\ar[r]&\Ext^1(R/J,S_\un^G)\ar`r[d]`[l]`[dlll]`[dll][dll]\\
&\ \ \cdots&\cdots&\makebox[13ex]{$\cdots$}\ \ \ar`r[d]`[l]`[dlll]`[dll][dll]\\
&\Ext^n(S_\un^G,S_\un^G)\ar[r]&\Ext^n(I/J,S_\un^G)\ar[r]&\Ext^n(R/J,S_\un^G)\ar`r[d]`[l]`[dlll]`[dll][dll]\\
&\Ext^{n+1}(S_\un^G,S_\un^G)\ar[r]&\Ext^{n+1}(I/J,S_\un^G)\ar[r]&\Ext^{n+1}(R/J,S_\un^G)\cdots\\
}
$$
Now 
\begin{eqnarray*}
\Ext^n(R/J,S_\un^G)&\cong&\Ext^n(L,S_\un^G)\oplus\dirsum{X\in K_G(H)-\{T\}}\Ext^n(S_X^G,S_\un^G)\\
&\cong& \Ext^n(L,S_\un^G)\oplus\dirsum{X\in K_G(H)-\{T\}} \Ext^{n-1}(S_\un^G,S_\un^G)\\
&=&L(n)\oplus\dirsum{X\in\mathcal{X}}E_G(n-1)\mvirg
\end{eqnarray*}
by Corollary~\ref{central p}. This completes the proof.\findemo
\pagebreak[3]
\masubsect{The case $p>2$}
\begin{mth}{Theorem} Let $k$ be a field of odd characteristic $p$. Let $G$ be an elementary abelian $p$-group of rank 2. Then the simple functor $S_\un^G$ of $\comack_k(G)$ has exponential growth.
\end{mth}
\vspace{1ex}
\pf Let $H$ be a subgroup of index $p$ in $G$, and consider the long exact sequence of Corollary~\ref{long exact sequence}. The set $K_G(H)$ consists of $p$ subgroups of order~$p$ in $H$, so $\mathcal{X}$ has cardinality $p-1$. Moreover
$$E_H(n)= \Ext^n_{\comack_k(H)}(S_\un^H,S_\un^H)\cong k\mvirg$$
by Proposition~\ref{ext cyclic}, since $H$ is cyclic of order $p\geq 3$.
\medskip
\par 
This yields a long exact sequence of the form
$$\cdots\!\to\! k\!\to\! L(n)\oplus E_G(n-1)^{p-1}\!\to\! E_G(n+1)\!\to\! k\!\to\! L(n+1)\oplus E_G(n)^{p-1}\!\to\!\cdots\mpoint$$
In particular, by Lemma~\ref{2 out of 3}, this shows that
$$\dim_k\big(L(n)\oplus E_G(n-1)^{p-1}\big)\leq 1+\dim_kE_G(n+1)\mvirg$$
thus $\dim_kE_G(n+1)\geq (p-1)\dim_kE_G(n-1)-1$, for any $n\geq 1$. Hence, by induction on $n$, 
$$\dim_kE_G(2n)-\frac{1}{p-2}\geq (p-1)^{n}\big(\dim_kE_G(0)-\frac{1}{p-2}\big)= \frac{(p-1)^{n}(p-3)}{p-2}\mvirg$$
$$\dim_kE_G(2n+1)-\frac{1}{p-2}\geq (p-1)^{n}\big(\dim_kE_G(1)-\frac{1}{p-2}\big)= \frac{(p-1)^{n}(2p-5)}{p-2}\mvirg$$
since $\dim_kE_G(1)=2$ by Lemma~\ref{ext1s1s1}. Thus $S_\un^G$ has exponential growth if $p>3$. 
\medskip
\par
Now if $p=3$, the functor $L$ is simple, and isomorphic to $S_\un^G$, by Proposition~\ref{index p elemab}. It follows that in this case
$$\dim_kE_G(n)+(p-1)\dim_kE_G(n-1)-1\leq \dim_kE_G(n+1)\mpoint$$
In particular $\dim_kE_G(2)\geq 2+2-1=3$, since $\dim_kE_G(1)=2$ by Lemma~\ref{ext1s1s1}. Since $\dim_kE_G(n+1)\geq 2\dim_kE_G(n-1)-1$, it follows that
$$\dim_kE_G(2n)-1\geq 2^{n-1}\big(\dim_kE_G(2)-1\big)\geq 2^{n}\mvirg$$
$$\dim_kE_G(2n+1)-1\geq 2^{n}\big(\dim_kE_G(1)-1\big)\geq 2^{n}\mvirg$$
hence $S_\un^G$ has exponential growth in this case also.\findemo
\vspace{4ex}
\pagebreak[4]
\masubsect{The case $p=2$}
\begin{mth}{Theorem} \label{odd Ext}Let $k$ be a field of characteristic 2, and $G$ be an elementary abelian $2$-group. Then
$$\forall j\in\N,\;\;\Ext^{2j+1}_{\comack_k(G)}(S_\un^G,S_\un^G)=\zero\mpoint$$
\end{mth}
\pf By induction on the rank $m$ of $G$, starting with $m=0$, i.e. $G=\un$, where the result is trivial. Suppose $m>0$, and that the result holds for all elementary abelian 2-groups of rank smaller than $m$. Choose a subgroup $H$ of index $p$ in $G$, and a complement $T\in K_G(H)$, and consider the long exact sequence of Corollary~\ref{long exact sequence}.\par
In this case by Proposition~\ref{index p elemab}, since $p-2=0$, the functor $L$ is equal to zero. Thus $L(n)=\zero$ for any $n\in\N$.
Moreover $\Ext^1_{\comack_k(G)}(S_\un^G,S_\un^G)=\zero$, by Theorem~\ref{ext1s1s1}, since $G$ is generated by involutions. This starts an induction argument on $j$~: suppose that $\Ext^{2r+1}_{\comack_k(G)}(S_\un^G,S_\un^G)=\zero$, for $0\leq r<j$. The exact sequence of Corollary~\ref{long exact sequence} becomes
\begin{equation}\label{the sequence}\xymatrix{
\cdots\ar[r]&*!U(0.6){\dirsum{X\in\mathcal{X}}E_G(n-2)}\ar[r]&E_G(n)\ar[r]&E_H(n)\ar`r[d]`[l]`[dlll]`[dll][dll]\\
&*!U(0.6){\dirsum{X\in \mathcal{X}}E_G(n-1)}\ar[r]&E_G(n+1)\ar[r]&E_H(n+1) \cdots\mvirg\\
}
\end{equation}
where $E_G(n)=\Ext^{n}_{\comack_k(G)}(S_\un^G,S_\un^G)$ and $E_H(n)=\Ext^{n}_{\comack_k(H)}(S_\un^H,S_\un^H)$, and $\mathcal{X}=K_G(H)-\{T\}$.\par
Set $n=2j+1$ in this sequence. By induction hypothesis on $m$
$$E_H(n)=E_H(2j+1)=\Ext^{2j+1}_{\comack_k(H)}(S_\un^H,S_\un^H)\cong\zero\mvirg$$
since $H$ is elementary abelian of rank $m-1$. Also 
$$E_G(n-2)=E_G(2j-1)=\Ext^{2j-1}_{\comack_k(G)}(S_\un^G,S_\un^G)=\zero\mvirg$$
by induction hypothesis on $j$.\par
It follows that $E_G(n)=\Ext^{2j+1}(S_\un^G,S_\un^G)=\zero$, and this completes the inductive step on $j$, hence the inductive step on $m$.\findemo
\begin{mth}{Theorem} \label{ses p=2}Let $k$ be a field of characteristic 2. Let $G$ be an elementary abelian 2-group, and $H$ be a subgroup of index 2 of $G$. Choose $T\in K_G(H)$, and set $\mathcal{S}=K_G(H)-\{T\}$. \par
Then for any $j\in \N$, there is a short exact sequence of extension groups 
$$0\to\dirsum{X\in \mathcal{S}}\Ext^{2j-2}(S_\un^G,S_\un^G)\stackrel{\gamma}{\to}\Ext^{2j}(S_\un^G,S_\un^G)\stackrel{r_H^G}{\to} \Ext^{2j}(S_\un^H,S_\un^H)\To 0$$
where $r_H^G$ is induced by restriction to $H$, and $\gamma$ is the direct sum of the maps induced by the Yoneda product with the element $\gamma_X$ of $\Ext^{2}(S_\un^G,S_\un^G)$, for $X\in\mathcal{S}$.
\end{mth}
\pf This follows from the exact sequence~\ref{the sequence}, in which all the terms $E_G(2j-1)$, $E_G(2j+1)$, and $E_H(2j+1)$ are equal to zero. The exact sequence~\ref{the sequence} splits as a series of short exact sequences
$$0\to\dirsum{X\in \mathcal{S}}\Ext^{2j-2}(S_\un^G,S_\un^G)\stackrel{\gamma}{\to}\Ext^{2j}(S_\un^G,S_\un^G)\stackrel{r}{\to} \Ext^{2j}(S_\un^H,S_\un^H)\To 0\mvirg$$
for each $j\in\N$. It remains to show that the morphisms $\gamma$ and $r$ are as stated in Theorem~\ref{ses p=2}.\par
First, with the notation of Proposition~\ref{index p elemab}, the morphism 
$$t: \Ext^{2j-1}(R/J,S_\un^G)\To \Ext^{2j}(S_\un^G,S_\un^G)$$
is the transition morphism associated to the short exact sequence
$$T:\;\; 0\To R/J\To I/J\To S_\un^G\To 0\mvirg$$
so it it given by Yoneda product with $T$.\par
Via the isomorphism $R/J\cong\dirsum{X\in\mathcal{S}}S_X^G$, the morphism $t$ can be viewed as the direct sum
of morphisms $t_X:\Ext^{2j-1}(S_X^G,S_\un^G)\to \Ext^{2j}(S_\un^G,S_\un^G)$, where $t_X$ is given by composition with the sequence $T_X$ obtained from $T$ by the projection $\pi_X$ to the summand~$S_X^G$, i.e. by completing the diagram
$$\xymatrix{
T_{\phantom{X}}:\;\;0\ar[r]&*!U(0.4){\dirsum{X\in\mathcal{S}}S_X^G}\ar[r]\ar[d]^-{\pi_X}&I/J\ar[r]\ar[d]^-{f}&S_\un^G\ar[r]\ar@{=}[d]&0\\
T_X:\;\;0\ar[r]&S_X^G\ar[r]&Y\ar[r]&S_\un^G\ar[r]&0\mpoint\\
}
$$
Now $f$ is surjective, since $\pi_X$ and $\Id_{S_\un^G}$ are. Since $S_X^G$ is not a quotient of $I/J$, it follows that $Y$ is a non-split extension of $S_\un^G$ by $S_X^G$. Thus $Y\cong\jiota_{G/X}^G(S_\un^{G/X})$, by Lemma~\ref{ext1sxs1}. More precisely $T_X$ is isomorphic to the sequence $D_X^*$ of Lemma~\ref{D_Z}.\par 
Moreover, by Corollary~\ref{central p}, the isomorphism 
$$\Ext^{2j-2}(S_\un^G,S_\un^G)\cong \Ext^{2j-1}(S_X^G,S_\un^G)$$
is given by Yoneda composition with $D_X$. Hence the component index by $X\in\mathcal{S}$ of the morphism 
$$\gamma:\dirsum{X\in\mathcal{S}}\Ext^{2j-2}(S_\un^G,S_\un^G)\To \Ext^{2j}(S_\un^G,S_\un^G)$$ is given by composition with $\Gamma_X=D_X\circ D_X^*$, i.e. by Yoneda product with~$\gamma_X$, as claimed.\par
Now the morphism $r_H^G:\Ext^{2j}_{\comack_k(G)}(S_\un^G,S_\un^G)\to \Ext^{2j}_{\comack_k(H)}(S_\un^H,S_\un^H)$ is composed of two steps~: first, taking the image by the projection map $q:I\to S_\un^G$, i.e. taking the Yoneda product with $q\in \Ext^0_{\comack(G)}(\Ind_H^GS_\un^H,S_\un^G)$, and then
using the adjunction $(\Ind_H^G,\Res_H^G)$, which gives a map
$$\Ext^{2j}_{\comack_k(G)}(I,S_\un^G)\cong\Ext^{2j}_{\comack_k(H)}(S_\un^H,\Res_H^GS_\un^G)\cong\Ext^{2j}_{\comack_k(H)}(S_\un^H,S_\un^H)\mpoint$$
In other words, with the notation of Remark~\ref{Yoneda compatible}, for $u\in \Ext^{2j}_{\comack_k(G)}(S_\un^G,S_\un^G)$
$$r_H^G(u)=\alpha_{2j}(u\circ q)=\Res_H^G(u)\circ\alpha_0(q)=\Res_H^G(u)\mvirg$$
since the map $\alpha_0(q):S_\un^H\to \Res_H^GS_\un^G\cong S_\un^H$ obtained by adjunction is the identity map. Hence $r_H^G$ is induced by $\Res_H^G$, and this completes the proof.\findemo
The following corollary is Theorem~\ref{Poincare series} of Section~\ref{intro}~:
\begin{mth}{Corollary} \label{Ext generators}Let $G$ be an elementary abelian 2-group of rank $m$.
Then the algebra $\Ext_{\comack_k(G)}^*(S_\un^G,S_\un^G)$ is finitely generated by the elements $\gamma_X$, where $X$ is a subgroup of order 2 of~$G$. Its Poincar\'e series
$$P(t)=\sum_{j\in\N}\dim_k\Ext_{\comack_k(G)}^j(S_\un^G,S_\un^G)\;t^j$$
is equal to
$$P(t)=\frac{1}{(1-t^2)(1-3t^2)(1-7t^2)\ldots\big(1-(2^{m-1}-1)t^2\big)}\mpoint$$
\end{mth}
\pf By induction on $m$, the case $m=0$ being trivial. In Theorem~\ref{ses p=2}, one can assume that the algebra $\Ext^*(S_\un^H,S_\un^H)$ is generated by the elements~$\gamma_X^H$, where $X$ is a subgroup of order 2 of $H$. This means that for $n\in\N$, any element in $\Ext^{2n}(S_\un^H,S_\un^H)$ is a $k$-linear combination of products of the form
$$\gamma_{X_1}^H\gamma_{X_2}^H\cdots\gamma_{X_n}^H\mvirg$$
where $X_1$,\ldots, $X_n$ are subgroups of order 2 of $H$.  \par
Denote by $\Gamma$ the $k$-linear subspace of $\Ext^{2n}(S_\un^G,S_\un^G)$ generated by the similar products $\gamma_{X_1}^G\gamma_{X_2}^G\cdots\gamma_{X_n}^G$.\par
Now for $X\leq H$, by Proposition~\ref{composition}, there are isomorphisms of functors 
$$\Res_H^G\circ \iiota_{G/X}^G\cong \mathsf{L}_{G}\circ\mathsf{L}_{G/X}\cong \mathsf{L}_{G\times_G(G/X)}\mvirg$$
where $G$ is viewed as an $(H,G)$-biset and $G/X$ as a $(G,G/X)$-biset. Now $G\times_G(G/X)\cong (H/X)\times_{H/X}(G/X)$ as $(H,G/X)$-bisets, thus
$$\Res_H^G\circ \iiota_{G/X}^G\cong \iiota_{H/X}^H\circ\Res_{H/X}^{G/X}\mpoint$$
It follows that $\Res_H^G\bin{X}{\un}{G}\cong\bin{X}{\un}{H}$. Similarly $\Res_H^G\bin{\un}{X}{G}\cong\bin{\un}{X}{H}$. Thus $r_H^G\gamma_X^G=\gamma_X^H$, because $\Res_H^G$ respects the Yoneda composition.\par
It follows that 
$$\gamma_{X_1}^H\gamma_{X_2}^H\cdots\gamma_{X_n}^H=r_H^G(\gamma_{X_1}^G\gamma_{X_2}^G\cdots\gamma_{X_n}^G)\mvirg$$
thus $\Ext^{2n}(S_\un^G,S_\un^G)=\Gamma+\Ker\,r_H^G$. But
$$\Ker\,r_H^G=\gamma\big(\Ext^{2n-2}(S_\un^G,S_\un^G)\big)=\sum_{X\in\mathcal{S}}\Ext^{2n-2}(S_\un^G,S_\un^G)\circ\gamma_X^G\mvirg$$
and by induction on $n$, this shows that $\Gamma=\Ext^{2n}(S_\un^G,S_\un^G)$.\par
Denote by $P_m(t)$ the Poincar\'e series of the algebra $\Ext^*_{\comack_k(G)}(S_\un^G,S_\un^G)$, where~$G$ is an elementary abelian 2-group of rank~$m$. The short exact sequence of Theorem~\ref{ses p=2} gives the relation
$$P_m(t)=|\mathcal{S}|t^2P_m(t)+P_{m-1}(t)\mvirg$$
and moreover $|\mathcal{S}|=2^{m-1}-1$. Thus
$$P_m(t)=P_{m-1}(t)\cdot\frac{1}{1-(2^{m-1}-1)t^2}\mvirg$$
and the claimed formula follows by induction on $m$.\findemo
\pagebreak[3]
\begin{mth}{Theorem} Let $k$ be a field of characteristic 2, and $G$ be an elementary abelian 2-group of rank~$m$. 
\begin{enumerate}
\item {\rm [Samy Modeliar~\cite{samymodeliar}]} If $m\leq~2$, the group $G$ is a poco group over $k$.
\item If $m\geq 3$, the group $G$ is not a poco group over $k$. More precisely, the simple functor $S_\un^G$ of $\comack_k(G)$ has exponential growth.
\end{enumerate}
\end{mth}
For Assertion~2, observe that
$$\frac{1}{(1-t^2)(1-3t^2)\ldots\big(1-(2^{m-1}-1)t^2\big)}=\prod_{j=1}^{m-1}\left(\sum_{n_j=0}^{\infty}(2^j-1)^{n_j}t^{2n_j}\right)\mvirg$$
Thus
\begin{eqnarray*}
\dim_k\Ext^{2n}(S_\un^G,S_\un^G)&=&\sum_{n_1+\cdots+n_{m-1}=n}\prod_{j=1}^{m-1}(2^j-1)^{n_j}\\
&\geq&(2^{m-1}-1)^n\mpoint
\end{eqnarray*}
By Corollary~\ref{central p}, since moreover $\Ext^{2n+1}(S_Z^G,S_\un^G)\cong \Ext^{2n+1}(S_\un^G,S_Z^G)$ as simple functors are self dual, it follows that
$$\dim_k\Ext^{2n+1}(S_\un^G,S_Z^G)\geq (2^{m-1}-1)^n\mvirg$$
thus $S_\un^G$ has exponential growth, by Lemma~\ref{simple enough}, if $2^{m-1}-1>1$, i.e. if $m\geq 3$.\findemo
Assertion~1 is trivial if $m=0$. It is straightforward if $m=1$~: since every $kG$-module for a group $G$ of order 2 is a permutation module, it follows that all the fixed points functors are projective in this case, so $\comack_k(G)$ has finite projective dimension in this case. \par
For $m=2$, M. Samy Modeliar has described explicit eventually periodic resolutions of the fixed points functors $FP_V$,
 for all indecomposable $kG$-modules $V$. An alternative proof can be sketched as follows~: the Poincar\'e series of $\Ext^*(S_\un^G,S_\un^G)$ is equal to $\frac{1}{1-t^2}$ in this case. It follows that $\Ext^n(S_\un^G,S_\un^G)$ is equal to $\zero$ if $n$ is odd, and one dimensional if $n$ is even. If $Q$ is a non trivial subgroup of $G$, then there is an exact sequence
\begin{equation}\label{klein}
0\To \dirsum{X\in K_Q(Z)}S_X^G\To\iiota_{G/Z}^G(S_{Q/Z}^{G/Z})\To S_Q^G\To0\mvirg
\end{equation}
where $Z$ is a subgroup of order $2$ of $Q$. Moreover $\iiota_{G/Z}^G(S_{Q/Z}^{G/Z})$ has finite projective dimension, since $G/Z$ has order 2, and since $\iiota_{G/Z}^G$ is exact and preserves projectives. It follows easily by induction on $|Q|$ that there exists a constant $c_Q$ such that $\dim_k\Ext^j(S_Q^G,S_\un^G)\leq c_Q$ for any $j\in\N$. \par
Thus $\dim_k\Ext^j(S_\un^G,S_R^G)\leq c_R$ for any $j\in\N$, and any $R\leq G$. Using again the exact sequence~\ref{klein}, one can show by induction on $|Q|$ that for any $Q, R\leq G$, there is a constant $c_{Q,R}$ such that $\dim_k\Ext^j(S_Q^G,S_R^G)\leq c_{Q,R}$ for any $j\in\N$. By Lemma~\ref{simple enough} and Corollary~\ref{poco simple}, the group $G$ is a poco group.\findemo
\section{The case of 2-groups of sectional rank at most 2}\label{sectional 2-rank 2}
\begin{mth}{Proposition} \label{s1 pol}Let $k$ be a field of characteristic 2, let $G$ be a 2-group, and let $H$ be a subgroup of index 2 of $G$. If the functor $S_\un^H$ (over $k$) has polynomial growth, and if the functor $S_X^{C_G(X)}$ has polynomial growth, for any $X\in K_G(H)$, then the functor $S_\un^G$ has polynomial growth.
\end{mth}
\pf Consider the functor $I=\Ind_H^GS_\un^H$. By Proposition~\ref{index p}, there is a filtration
$$I\supset R\supseteq S\supset \zero\mvirg$$
where $R$ is the radical and $S$ is the socle of $I$, such that $I/R\cong S\cong S_\un^G$, and $R/S\cong \dirsum{X\in [G\dom K_G(H)]}S_X^G$. This gives two short exact sequences in $\comack_k(G)$
\begin{equation}\label{suite1}
0\To S_\un^G\To R\To \dirsum{X\in [G\dom K_G(H)]}S_X^G\To 0
\end{equation}\vspace{-.5cm}
\begin{equation}\label{suite2}
0\To R\To I\To S_\un^G\To 0\mpoint
\end{equation}
Let $M$ be a finitely generated cohomological Mackey functor for $G$ over~$k$. Applying the functor $\Hom_{\comack_k(G)}({-},M)$ to the sequence~(\ref{suite1}) gives the following long exact sequence of $\Ext$ groups in $\comack_k(G)$ 
$$\cdots\To\dirsum{X\in [G\dom K_G(H)]}\Ext^{n}(S_X^G,M)\To \Ext^n(R,M)\To \Ext^n(S_\un^G,M)\To\cdots$$
By Lemma~\ref{2 out of 3}, this gives
$$r_n\leq e_n+ \sum_{X\in [G\dom K_G(H)]}e_{X,n}\mvirg$$
where $r_n=\dim_k\Ext^n(R,M)$, $e_{X,n}=\dim_k\Ext^{n}(S_X^G,M)$, for $X\in K_G(H)$, and $e_n=\dim_k \Ext^n(S_\un^G,M)$.\par
Now applying the functor $\Hom_{\comack_k(G)}({-},M)$ to the sequence~(\ref{suite2}) gives the following long exact sequence of $\Ext$ groups in $\comack_k(G)$ 
$$\cdots\To \Ext^{n}(R,M)\To \Ext^{n+1}(S_\un^G,M)\To \Ext^{n+1}(I,M)\To\cdots$$
By Lemma~\ref{2 out of 3}, this gives
$$e_{n+1}\leq r_n+i_{n+1}\mvirg$$
where $i_{n}=\dim_k\Ext^n(I,M)$. Thus
\begin{equation}\label{bound}
e_{n+1}\leq e_n + (\sum_{X\in [G\dom K_G(H)]}e_{X,n}) + i_{n+1}\mpoint
\end{equation}
Now $S_X^G=\Ind_{C_G(X)}^GS_X^{C_G(X)}$, for any $X\in K_G(H)$, thus
$$e_{X,n}=\dim_k\Ext^n_{\comack_k(C_G(X))}(S_X^{C_G(X)},\Res_{C_G(X)}^GM)\mpoint$$
Hence if $S_X^{C_G(X)}$ has polynomial growth, there are constants $c_X$, $d_X$, and $e_X$, such that
$$ \forall n\in\N,\;\;e_{X,n}\leq c_Xn^{d_X}+e_X\mpoint$$
Thus
$$\forall n\in\N,\;\;\sum_{X\in [G\dom K_G(H)]}e_{X,n}\leq Cn^D+E\mvirg$$
where $C=\sum_{X\in [G\dom K_G(H)]}\limits c_X$, $D=\max_{X\in [G\dom K_G(H)]}\limits d_X$, and $E=\sum_{X\in [G\dom K_G(H)]}\limits e_X$.
Similarly, 
$$i_n=\dim_k\Ext^n_{\comack_k(H)}(S_\un^H,\Res_H^GM)\mvirg$$
so if $S_\un^H$ has polynomial growth, there are constants $c$, $d$, and $e$ such that
$$\forall n\in\N,\;\;i_n\leq c n^d+e\mpoint$$
Thus 
$$\forall n\in\N,\;\;i_{n+1}\leq c|(n+1)^d-1|+c+e\leq c(2n)^d+c+e=c2^dn^d+c+e\mpoint$$
Inequality \ref{bound} now gives
$$\forall n\in\N,\;\;e_{n+1}\leq e_n+ \gamma n^\delta+\epsilon$$
where $\gamma=c2^d+C$, $\delta=\max(D,d)$, and $\epsilon=E=c=e$. By induction, it follows that
$$\forall n\in\N,\;\;e_{n}\leq \gamma \sum_{j=1}^nj^\delta + n\epsilon+e_0\leq \gamma n^{\delta+1}+n\epsilon+e_0\leq (\gamma+\epsilon)n^{\delta+1}+e_0\mvirg$$ 
so $S_\un^G$ has polynomial growth.\findemo
\begin{mth}{Theorem} Let $k$ be a field of characteristic 2, and $G$ be a 2-group of sectional 2-rank at most 2. Then $G$ is a poco group over $k$.
\end{mth}
\pf By induction on the order of $G$, one may assume that for any section $(B,A)\neq (G,\un)$, the group $B/A$ is a poco group. By Corollary~\ref{poco simple}, proving that $G$ is a poco group is equivalent to proving that for any subgroup $Q$ of~$G$, the simple functor $S_Q^G$ has polynomial growth. First consider the case where $Q=\un$.\par
Let $H$ be a subgroup of index 2 of $G$. Then $H$ is a poco group by induction hypothesis. Similarly, if $X\in K_G(H)$ and $C_G(X)\neq G$, then $C_G(X)$ is a poco group. Thus, if $H$ has no central complement in $G$, by Proposition~\ref{s1 pol}, the functor $S_\un^G$ has polynomial growth. This holds in particular if $H$ has no complement at all in $G$, i.e. if $H$ contains all involutions of $G$. Thus, if $G$ is not generated by involutions, then $S_\un^G$ has polynomial growth. So one can assume that $G$ is generated by involutions.\par
Now if $H$ has some central complement $X$ in $G$, then $G\cong H\times X$. Since~$G$ has sectional 2-rank at most 2, the group $H$ has sectional 2-rank at most 1, hence it is cyclic. Thus $G$ is abelian, and generated by involutions, hence elementary abelian, and of sectional 2-rank at most 2. So the group $G$ is elementary abelian of rank at most 2, and the results of Section~\ref{elemab} show that $S_\un^G$ has polynomial growth.\par
In any case $S_\un^G$ has polynomial growth. By induction on the order of $Q$, it follows that $S_Q^G$ has polynomial growth, for any subgroup $Q$ of $G$~: 
indeed, setting $B=N_G(Q)$ and $A=\Phi(Q)$, Lemma~\ref{simple from elemab} shows that
$$S_Q^G\cong \mathsf{L}_{U}(S_{Q/A}^{B/A})\mvirg$$
where $U$ is the set $A\dom G$, for its natural structure of $(B/A,G)$-biset structure. If $(B,A)\neq (G,\un)$, then $S_{Q/A}^{B/A}$ has polynomial growth, hence $S_Q^G$ has polynomial growth, by Lemma~\ref{pol to pol}.\par
So it remains to consider the case where $(B,A)=(Q,\un)$, i.e. the case where $Q$ is a non trivial normal elementary abelian subgroup of $G$. Let $Z$ be a subgroup of order 2 of $Q\cap Z(G)$.\par
By Proposition~\ref{iota S_Q}, there is a short exact sequence
$$0\To \dirsum{X\in [G\dom K_Q(Z)]}S_{X}^G\To\iiota_{G/Z}^G(S_{Q/Z}^{G/Z})\To S_Q^G\To 0\mpoint
$$
By the induction hypothesis, all the functors $S_X^G$, for $X\in K_Q(Z)$, have polynomial growth. Moreover, the functor $\iiota_{G/Z}^G(S_{Q/Z}^{G/Z})$ has polynomial growth, by Lemma~\ref{pol to pol}. Hence $S_Q^G$ has polynomial growth, by Lemma~\ref{2 out of 3}, and this completes the proof.\findemo
\begin{rem}{Remark} The 2-groups of sectional rank at most 1 are the cyclic 2-groups. In a recent preprint, Andersen, Oliver and Ventura have shown that the 2-groups of sectional rank 2 are the metacyclic 2-groups (see Lemma~10.2 of~\cite{andersenoliverventura}). For 2-groups of order at least 32, this also follows from a theorem of Blackburn~(see \cite{blackburn}, or Satz 11.13 in~\cite{huppert}).
\end{rem}
\vspace{.8cm}
\pagebreak[3]
\centerline{\Large\bf \romain{2} - Cohomology}
\section{Extension of simple functors for elementary abelian $p$-groups}\label{extensions elemab}
In this section $k$ is a field of characteristic $p$, the group $G\cong (C_p)^m$ is an elementary abelian $p$-group of rank $m$, in additive notation, and all the cohomological Mackey functors have values in $k$-vector spaces.\par
Recall from Notation~\ref{gamma_X} that $\gamma_X\in\Ext^2_{\comack_k(G)}(S_\un^G,S_\un^G)$ is the element represented by the 2-fold extension
$$0\To S_\un^G\To \bin{X}{\un}{G}\To \bin{\un}{X}{G}\To S_\un^G\To 0\mpoint$$
\begin{mth}{Notation} \label{other gamma_x}If $G=(C_p)^m$ and $x\in G-\{0\}$, set $\gamma_x=\gamma_{{<}x{>}}^G$. Thus $\gamma_x\in \Ext^2(S_\un^G,S_\un^G)$.\par
\end{mth}
\begin{mth}{Lemma} \label{ext2 zero}Let $A$ be a ring. If $L\supseteq M\supseteq N\supseteq \zero$ is a filtration of an $A$-module $L$, then the exact sequence
$$0\To N\To M\To L/N\To L/M\To 0$$
obtained by splicing the short exact sequences
$$\xymatrix@R=1.5ex@C=1.5ex{
0\ar[rr]&& N\ar[rr]&& M\ar[rr]\ar[dr]&& L/N\ar[rr]&&L/M\ar[rr]&&0\\
&&&&&M/N\ar[dr]\ar[ur]&&&\\
&&&&0\ar[ur]&&0&&
}
$$
represents the zero class of $\Ext^2_A(L/M,N)$.
\end{mth}
\pf There is a commutative diagram
$$\xymatrix{
0\ar[r]& N\ar[r]\ar@{=}[d]& M\ar[r]\ar[d]^i& L/N\ar[r]^p\ar[d]^-{(\Id,p)}&L/M\ar[r]\ar@{=}[d]&0\\
0\ar[r]& N\ar[r]& L\ar[r]^-{(q,0)}& L/N\oplus L/M\ar[r]^-{\binom{0}{\Id}}&L/M\ar[r]&0\\
}
$$
where $i$ is the inclusion map, and $p$, $q$ are projection maps. Moreover, the bottom sequence represents zero in $\Ext^2_A(L/M,N)$, since the map $\binom{0}{\Id}$ is split surjective.\findemo
\begin{mth}{Lemma} \label{XY}Let $G$ be an elementary abelian $p$-group. If $X$ and $Y$ are distinct subgroups of order $p$ of $G$, set $Q=XY$. Then the sequences 
$$0\To S_X^G\To \bin{\un}{X}{G}\To \bin{Y}{\un}{G}\To S_Y^G\To 0$$
and
$$0\To S_X^G\To \bin{Q}{X}{G}\To \bin{Y}{Q}{G}\To S_Y^G\To 0$$
represent opposite elements of $\Ext^2_{\comack_k(G)}(S_Y^G,S_X^G)$.
\end{mth}
\pf Indeed, the sum of the corresponding elements in $\Ext^2_{\comack_k(G)}(S_Y^G,S_X^G)$ is represented by the sequence
\begin{equation}\label{zero}
0\To S_X^G\To \binom{S_\un^G\;S_Q^G}{S_X^G}\To \binom{S_Y^G}{S_\un^G\;S_Q^G}\To S_Y^G\To 0\mvirg
\end{equation}
where the functor $L_X=\binom{S_\un^G\;S_Q^G}{S_X^G}$ is the Mackey functor for $G$ over $k$ whose values at $\un$, $X$ and~$Q$ are equal to $k$, and other values are zero. The conjugation maps for this functor are all identity maps (possibly zero), and the possibly non zero transfer and restriction maps are given in the following diagram
$$\xymatrix{
L_X(Q)=k\ar@/_/[d]_-{r_X^Q=1}\\
L_X(X)=k\ar@/_/[u]_-{t_X^Q=0}\ar@/_/[d]_-{r_\un^X=0}\\
L_X(\un)=k\ar@/_/[u]_-{t_\un^X=1}
}
$$
It follows in particular that $S_X^G$ is isomorphic to a subfunctor of $L_X$, and this yields an exact sequence 
\begin{equation}\label{zero1}
0\To S_X^G\To \binom{S_\un^G\;S_Q^G}{S_X^G}\To S_\un^G\oplus S_Q^G\To 0
\end{equation}
in $\comack_k(G)$.\par
The functor $\binom{S_Y^G}{S_\un^G\;S_Q^G}$ is isomorphic to the dual $L_Y^*$ of $L_Y=\binom{S_\un^G\;S_Q^G}{S_Y^G}$~: its non zero values are at $\un$, $Y$ and $Q$, and they are equal to $k$. The conjugation maps for $L_Y^*$ are all identity maps (possibly zero), and the possibly non zero transfer and restriction maps are given in the following diagram
$$\xymatrix{
L_Y^*(Q)=k\ar@/_/[d]_-{r_Y^Q=0}\\
L_Y^*(Y)=k\ar@/_/[u]_-{t_Y^Q=1}\ar@/_/[d]_-{r_\un^Y=1}\\
L_Y^*(\un)=k\ar@/_/[u]_-{t_\un^Y=0}
}
$$
There is an exact sequence\begin{equation}\label{zero2}
0\To S_\un^G\oplus S_Q^G\To \binom{S_Y^G}{S_\un^G\;S_Q^G}\To S_Y^G\To 0
\end{equation}
in $\comack_k(G)$, and the exact sequence~\ref{zero} is obtained by splicing the exact sequences~\ref{zero1} and~\ref{zero2}.\par
Let $M_{X,Y}$ denote the functor defined for $H\leq G$ by $M_{X,Y}(H)=k$ if $H\in\{\un,X,Y,Q\}$, and by $M_{X,Y}(H)=\zero$ otherwise. The conjugation maps for $M_{X,Y}$ are all identity maps (possibly zero), and the possibly non zero transfer and restriction maps for $M_{X,Y}$ are given in the following diagram 
\begin{equation}\label{zero3}
\begin{array}{c}\xymatrix{
&M_{X,Y}(Q)=k\ar@/_/[dl]_-{r_X^Q=1}\ar@/_/[dr]_-{r_Y^Q=0}&\\
M_{X,Y}(X)=k\ar@/_/[ur]_-{t_X^Q=0}\ar@/_/[dr]_-{r_\un^X=0}&&M(Y)=k\ar@/_/[ul]_-{t_Y^Q=1}\ar@/_/[dl]_-{r_\un^Y=1}\\
&M_{X,Y}(\un)=k\ar@/_/[ul]_-{t_\un^X=1}\ar@/_/[ur]_-{t_\un^Y=0}&
}\end{array}
\end{equation}
One checks easily that $M_{X,Y}$ is a cohomological Mackey functor for $G$ over~$k$, that $L_X$ is a subfunctor of $M_{X,Y}$ (represented by the left half of diagram~\ref{zero3}), and that $L_Y^*$ is a quotient functor of $M_{X,Y}$ (represented by the right half of~\ref{zero3}). Hence there is a filtration 
$$M_{X,Y}\supset L_X\supset S_X^G\supset\zero\mvirg$$
such that $M_{X,Y}/L_X\cong S_Y^G$ and $M_{X,Y}/S_X^G\cong L_Y^*$. Now Lemma~\ref{ext2 zero} shows that the exact sequence~\ref{zero} represents zero in $\Ext^2_{\comack_k(G)}(S_Y^G,S_X^G)$, and this completes the proof.\findemo
\begin{mth}{Lemma} \label{not yet}Let $G$ be an elementary abelian $p$-group, and let $Q<R$ be subgroups of $G$ such that $|R:Q|=p^2$. Let $M$ denote the functor $\Sigma_{Q,R}^G$ of Corollary~\ref{SigmaQR}. Then there is a filtration $M\supset J\supset S\supset \zero$, where $J$ is the radical of $M$ and $S$ is its socle. Moreover
$$M/J\cong S_R^G\mvirg\;\;S\cong S_Q^G\mvirg\;\;J/M\cong\dirsum{Q<X<R}S_X^G\mpoint$$
\end{mth}
\pf Recall that $\Sigma_{Q,R}^G$ is the subquotient of $FP_k$ whose value at $H\leq G$ is equal to $k$ if $Q\leq H\leq R$, and to zero otherwise. By Corollary~\ref{SigmaQR}, the head of $\Sigma_{Q,R}^G$ is simple, isomorphic to $S_R^G$, and its socle is simple, isomorphic to $S_Q^G$. In particular $S\leq J$. Moreover, if $H\leq G$, then $(J/S)(H)$ is equal to zero, except if $Q<H<R$, and in this case $(J/S)(H)= k$. It follows that
$$J/M\cong\dirsum{Q<X<R}S_X^G\mvirg$$
as was to be shown.\findemo
\begin{mth}{Proposition} \label{commutators}Let $G$ be an elementary abelian $p$-group, and $Q$ be a subgroup of order $p^2$ of $G$. Then for any subgroup $Y$ of order $p$ of $Q$
$$[\sum_{\un<X<Q}\gamma_X,\gamma_Y]=0\mvirg$$
in $\Ext^4_{\comack_k(G)}(S_\un^G,S_\un^G)$, where $[a,b]=ab-ba$, for $a,b\in \Ext^2_{\comack_k(G)}(S_\un^G,S_\un^G)$.
\end{mth}
\pf Recall that if $X$ is a subgroup of order $p$ of $G$, the element $\gamma_X$ of $\Ext^2_{\comack_k(G)}(S_\un^G,S_\un^G)$ is represented by the sequence
$$0\To S_\un^G\To\bin{X}{\un}{G}\To\bin{\un}{X}{G}\To S_\un^G\To 0\mpoint$$
It follows that if $Y$ is another subgroup of order $p$ of $G$, the Yoneda product $\gamma_X\gamma_Y$ is represented by the sequence
$$0\To S_\un^G\To\bin{X}{\un}{G}\To\bin{\un}{X}{G}\To\bin{Y}{\un}{G}\To\bin{\un}{Y}{G}\To S_\un^G\To 0\mpoint$$
This sequence is obtained by splicing the following three exacts sequences~:
$$0\To S_\un^G\To\bin{X}{\un}{G}\To S_{X}^{G}\To 0$$
\begin{equation}\label{comm}0\To S_X^G\To\bin{\un}{X}{G}\To\bin{Y}{\un}{G}\To S_Y^G\To 0
\end{equation}
$$0\To S_Y^G\To\bin{Y}{\un}{G}\To S_{\un}^{G}\To 0\mpoint$$
By Lemma~\ref{XY}, if $XY=Q$, i.e. if $X$ and $Y$ are distinct subgroups of order~$p$ of $Q$, the sequence~\ref{comm} and the sequence
\begin{equation}\label{comm2}0\To S_X^G\To\bin{Q}{X}{G}\To\bin{Y}{Q}{G}\To S_Y^G\To 0
\end{equation}
represent opposite elements of $\Ext^2_{\comack_k(G)}(S_Y^G,S_X^G)$. It follows that for a given subgroup $Y$ of order $p$ of $Q$ 
$$(\sumb{\un<X<Q}{X\neq Y}\gamma_X)\gamma_Y=-(\sumb{\un<X<Q}{X\neq Y} f_{X,Y}^Q)\gamma_Y\mvirg$$
where $f_{X,Y}^Q$ is the element represented by the sequence
$$0\To S_\un^G\To\bin{X}{\un}{G}\To\bin{Q}{X}{G}\To\bin{Y}{Q}{G}\To\bin{\un}{Y}{G}\To S_\un^G\To 0\mpoint$$
This in turn is the splice of the following two exact sequences
\begin{eqnarray*}
U_X^Q\phantom{)^*}&:& 0\To S_\un^G\To\bin{X}{\un}{G}\To\bin{Q}{X}{G}\To S_Q^G\To 0\\
(U_Y^{Q})^*&:&0\To S_Q^G\To \bin{Y}{Q}{G}\To\bin{\un}{Y}{G}\To S_\un^G\To 0\mpoint
\end{eqnarray*}
The sum of the elements $u_X^Q$ of $\Ext^2_{\comack_k(G)}(S_Q^G,S_\un^G)$ represented by the sequences~$U_X^Q$, for $\un<X<Q$, is represented by the sequence
$$0\To S_\un^G\To \binom{\Sigma_Q^G}{S_\un^G}\To \binom{S_\un^G}{\Sigma_Q^G}\To S_Q^G\To 0\mvirg$$
where $\dsp{\binom{\Sigma_Q^G}{S_\un^G}}$ has simple socle isomorphic to $S_\un^G$, and head isomorphic to $\Sigma_Q^G=\dirsum{\un<X<Q}S_X^G$, and $\dsp{\binom{S_\un^G}{\Sigma_Q^G}}$ is isomorphic to the dual of $\dsp{\binom{\Sigma_Q^G}{S_\un^G}}$. By Lemma~\ref{ext2 zero} and Lemma~\ref{not yet}, it follows that
$$\sum_{\un<X<Q}u_X^Q=0\mpoint$$
Hence $\sumb{\un<X<Q}{X\neq Y}u_X^Q=-u_Y^Q$, so
$$(\sumb{\un<X<Q}{X\neq Y}\gamma_X)\gamma_Y=u_Y^Q\gamma_Y\mpoint$$
This is represented by the sequence
$$0\To S_\un^G\To\bin{Y}{\un}{G}\To\bin{Q}{Y}{G}\To \bin{Y}{Q}{G}\To \bin{\un}{Y}{G}\To S_\un^G\To 0\mvirg$$
which is obviously self dual. It follows that
$$(\sumb{\un<X<Q}{X\neq Y}\gamma_X)\gamma_Y=\Big((\sumb{\un<X<Q}{X\neq Y}\gamma_X)\gamma_Y\Big)^*=\gamma_Y(\sumb{\un<X<Q}{X\neq Y}\gamma_X)\mpoint$$
Adding $\gamma_Y^2$ to both side, this shows that $\gamma_Y$ and $\mathop{\sum}_{\un<X<Q}\limits\gamma_X$ commute, which completes the proof.\findemo
\section{The algebra $\Ext^*_{\comack_k(G)}(S_\un^G,S_\un^G)$ for $G\cong (C_2)^m$}\label{algebra}
In this section $k$ is a field of characteristic $2$, the group $G\cong (C_2)^m$ is an elementary abelian $2$-group of rank $m$, in additive notation, and all the cohomological Mackey functors have values in $k$-vector spaces.\par
By Corollary~\ref{Ext generators} and Notation~\ref{other gamma_x}, the algebra $\Ext^*_{\comack_k(G)}(S_\un^G,S_\un^G)$ is generated by the elements $\gamma_x$ of $\Ext^2(S_\un^G,S_\un^G)$, for $x\in G-\{0\}$.
\begin{mth}{Proposition} \label{linear relations}Let $H$ be a subgroup of index~2 of $G$. Then
$$\sum_{x\notin H}\gamma_x=0\;\;\hbox{in}\;\;\Ext^2_{\comack_k(G)}(S_\un^G,S_\un^G)\mpoint$$
\end{mth}
\pf By Proposition~\ref{index p elemab}, there is a filtration $I\supset R\supset S\supset \zero$, where $I=\Ind_H^GS_\un^H$, such that $S\cong I/R\cong S_\un^G$, and i
$$R/S\cong\dirsum{X\in\mathcal{K}}S_X^G\mvirg$$
where $\mathcal{K}=K_G(H)$. By Lemma~\ref{ext2 zero}, the sequence
$$0\To S_\un^G\To R\To I/S\To S_\un^G\To 0$$
represents 0 in $\Ext^2(S_\un^G,S_\un^G)$. \par
For each $X\in\mathcal{K}$, the functors $\rho_{G/X}^G(I)$, $\rho_{G/X}^G(R)$, and $\rho_{G/X}^G(I/S)$ are all isomorphic to $S_\un^{G/X}$, since their non zero evaluation is at the trivial subgroup, where it is isomorphic to $k$. By adjunction, this gives morphisms 
$$p_X:\bin{X}{\un}{G}\to R\;\;\hbox{and}\;\;q_X:I/S\to\bin{\un}{X}{G}\mpoint$$
This gives a commutative diagram with exact lines
$$\xymatrix{
0\ar[r]&*!U(0.3){\dirsum{X\in \mathcal{K}}S_\un^G}\ar[r]\ar[d]_-{\Sigma}&\dirsum{X\in \mathcal{K}}\bin{X}{\un}{G}\ar[r]\ar[d]_-{p}&\dirsum{X\in \mathcal{K}}\bin{\un}{X}{G}\ar[r]\ar@{=}[d]&*!U(0.3){\dirsum{X\in \mathcal{K}}S_\un^G}\ar[r]\ar@{=}[d]&0\\
0\ar[r]&S_\un^G\ar[r]&R\ar[r]&\dirsum{X\in \mathcal{K}}\bin{\un}{X}{G}\ar[r]&*!U(0.3){\dirsum{X\in \mathcal{K}}S_\un^G}\ar[r]&0\\
0\ar[r]&S_\un^G\ar[r]\ar@{=}[u]&R\ar[r]\ar@{=}[u]&I/S\ar[r]\ar[u]_-{q}&S_\un^G\ar[r]\ar[u]_-{\Delta}&0\mvirg\\
}
$$
where $p$ is the sum of the maps $p_X$ and $q$ is the sum of the maps $q_X$, for $X\in\mathcal{K}$, where $\Sigma$ is the summation map, and $\Delta$ the diagonal inclusion. The top line of this diagram is the direct sum of the sequences $\Gamma_X$, for $X\in\mathcal{K}$ (see Notation~\ref{gamma_X}).\par
One checks easily that the top left square in this diagram is cocartesian, and that the bottom right square is cartesian. It follows that the bottom line represents the sum in $\Ext^2(S_\un^G,S_\un^G)$ of the 2-fold extensions in the top line, i.e. the sum $\sum_{X\in\mathcal{K}}\limits\gamma_X$. \par
Hence this sum is equal to 0, i.e. equivalently $\sum_{x\notin H}\limits\gamma_x=0$.\findemo
\masubsect{Proof of Theorem~\ref{presentation}}
Let $G\cong (C_2)^m$ be an elementary abelian 2-group of rank $m$, in additive notation, and $k$ be a field of characteristic $2$. Denote by $\mathcal{E}$ the graded algebra $\Ext^*_{\comack_k(G)}(S_\un^G,S_\un^G)$. \spn
$\bullet$ By Corollary~\ref{Ext generators} and Notation~\ref{other gamma_x}, the algebra $\mathcal{E}$ is generated by the elements $\gamma_x$, for $x\in G-\zero$, where $\gamma_x$ has degree~2.\par
Let $x$ and $y$ be distinct elements in $G-\zero$. Then $x$ and $y$ generate a subgroup $Q$ of order 4 of $G$, and the non zero elements of $Q$ are $x$, $y$, and $x+y$. By Proposition~\ref{commutators}, the commutator
$$[\gamma_x+\gamma_y+\gamma_{x+y},\gamma_{x+y}]$$
in $\mathcal{E}$ is equal to 0. Equivalently
\begin{equation}\label{quadratic}[\gamma_x+\gamma_y,\gamma_{x+y}]=0\mpoint
\end{equation}
Now if $H$ is a subgroup of index 2 of $G$, Proposition~\ref{linear relations} shows that
\begin{equation}\label{linear}
\sum_{x\notin H}\gamma_x=0
\end{equation}
in $\mathcal{E}$.\spn
$\bullet$ Conversely, let $\tilde{\mathcal{E}}$ denote the graded associative $k$-algebra with generators~$\tilde{\gamma}_x$ in degree 2, for $x\in G-\zero$, subject to the relations
$$\left\{\begin{array}{l}\forall H\leq G,\;|G:H|=2,\;\;\sum_{x\notin H}\limits\tilde{\gamma}_x=0\mvirg\\
\forall x,y\in G-\zero, x\neq y,\;\;[\tilde{\gamma}_x+\tilde{\gamma}_y,\tilde{\gamma}_{x+y}]=0\mpoint\end{array}\right.$$
Then there is a unique surjective homomorphism $s:\tilde{\mathcal{E}}\to \mathcal{E}$ of graded $k$-algebras such that $s(\tilde{\gamma}_x)=\gamma_x$, for all $x\in G-\zero$. Thus showing that $s$ is an isomorphism is equivalent to showing that for any integer $n\in \N$, the restriction $s_n$ of $s$ to the subspace $\tilde{\mathcal{E}}_n$ of elements of degree $n$ in $\tilde{\mathcal{E}}$ is an isomorphism onto the corresponding subspace $\mathcal{E}_n$ of $\mathcal{E}$. Since $s_n$ is surjective, this amounts to showing that $\dim_k\tilde{\mathcal{E}}_n\leq \dim_k\mathcal{E}_n$.\par
\begin{mth}{Lemma} \label{basis E_2}Let $B$ be an $\F_2$-basis of $G$. Then the set 
$$\{\tilde{\gamma}_x\mid x\in G-(B\sqcup\zero)\}$$
is a $k$-basis of $\tilde{\mathcal{E}}_2$, and $s_2$ is an isomorphism.
\end{mth}
\pf By Corollary~\ref{Ext generators}
$$\dim_k\mathcal{E}_2=(2^1-1)+(2^2-1)+\cdots+(2^{m-1}-1)=2^m-m-1\mpoint$$
Thus $\dim_k\tilde{\mathcal{E}}_2\geq 2^m-m-1$, since $s_2$ is surjective. As $2^m-m-1$ is precisely equal to the cardinality of $G-(B\sqcup\zero)$, it is enough to show that the subspace $\tilde{\mathcal{E}}'_2$ of $\tilde{\mathcal{E}}_2$ generated by the elements $\tilde{\gamma}_x$, for $x\in G-(B\sqcup\zero)$, is equal to $\tilde{\mathcal{E}}_2$. \par
Let $b\in B$, and denote by $H$ the subgroup of $G$ generated by $B-\{b\}$. Then $|G:H|=2$, so
$$\sum_{x\notin H}\tilde{\gamma}_x=0\mpoint$$
Since $G-H=H+b$, the only element of $G-H$ which is also in $B$ is $b$ itself. This gives
$$\tilde{\gamma}_b=\sum_{y\in H-\{0\}}\tilde{\gamma}_{b+y}\mvirg$$
so $\tilde{\gamma}_b\in\tilde{\mathcal{E}}'_2$. Since this holds for any $b\in B$, and since $\tilde{\mathcal{E}}_2$ is the set of $k$-linear combinations of elements $\tilde{\gamma}_x$, for $x\in G-\zero$, it follows that $\tilde{\mathcal{E}}'_2=\tilde{\mathcal{E}}_2$, as was to be shown.\findemo
\begin{mth}{Notation} Fix a linear ordering $b_1<b_2<\cdots<b_m$ on $B$. If $x\in G-\zero$, let $d_B(x)$ denote the least integer $i\in\{1,\ldots,m\}$ such that $x\in\;{<}b_1,\ldots,b_i{>}$.
\end{mth}
Then $d_B(x)=i$ if and only if $x=y+b_i$, for some $y\in{<}b_1,\ldots,b_{i-1}{>}$.
\begin{mth}{Definition} Let $n\in\N$. If $M=\tilde{\gamma}_{x_1}\tilde{\gamma}_{x_2}\ldots\tilde{\gamma}_{x_n}\in\tilde{\mathcal{E}}_{2n}$, where $x_i\in G-\zero$ for $1\leq i\leq n$, set $w(M)=\sum_{i=1}^n\limits d_B(x_i)$. The monomial $M$ is called {\em special} if $x_i\in G-(B\sqcup\zero)$, for $1\leq i\leq n$, and {\em ordered} if $d_B(x_i)\leq d_B(x_{i+1})$, for $1\leq i<n$.
\end{mth}
\begin{mth}{Lemma} \label{exchange}Let $x_1$,\ldots,$x_n$ be elements of $G-(B\sqcup \zero)$. Then for any $i\in\{1,\dots,n-1\}$ such that $d_B(x_i)\neq d_B(x_{i+1})$, the sum
$$\tilde{\gamma}_{x_1}\ldots\tilde{\gamma}_{x_{i-1}}\tilde{\gamma}_{x_{i}}\tilde{\gamma}_{x_{i+1}}\tilde{\gamma}_{x_{i+2}}\ldots\tilde{\gamma}_{x_{n}}+\tilde{\gamma}_{x_1}\ldots\tilde{\gamma}_{x_{i-1}}\tilde{\gamma}_{x_{i+1}}\tilde{\gamma}_{x_{i}}\tilde{\gamma}_{x_{i+2}}\ldots\tilde{\gamma}_{x_{n}}$$
is a linear combination of special monomials $M'$ with $w(M')>w(M)$.
\end{mth}
\pf Since the assertion is symmetric in $x_i$ and $x_{i+1}$, one can assume that $r=d_B(x_i)>d_B(x_{i+1})=s$. In this case $x_i=u+b_r$, where $u\in{<}b_1,\ldots,b_{r-1}{>}-\zero$, and $x_{i+1}=v+b_s$, where $v\in{<}b_1,\ldots,b_{s-1}{>}-\zero$. Let $t=x_i+x_{i+1}=(v+b_{s}+u)+b_r$, so $d_B(t)=r$. Now the relation $[\tilde{\gamma}_{x_{i+1}}+\tilde{\gamma}_t,\tilde{\gamma}_{x_i}]=0$ gives
$$\tilde{\gamma}_{x_{i}}\tilde{\gamma}_{x_{i+1}}+\tilde{\gamma}_{x_{i+1}}\tilde{\gamma}_{x_{i}}=\tilde{\gamma}_t\tilde{\gamma}_{x_{i}}+\tilde{\gamma}_{x_{i}}\tilde{\gamma}_t\mvirg$$
so the sum in the lemma is equal to
\begin{equation}\label{sum}
S=\tilde{\gamma}_{x_1}\ldots\tilde{\gamma}_{x_{i-1}}\tilde{\gamma}_t\tilde{\gamma}_{x_{i}}\tilde{\gamma}_{x_{i+2}}\ldots\tilde{\gamma}_{x_{n}}+\tilde{\gamma}_{x_1}\ldots\tilde{\gamma}_{x_{i-1}}\tilde{\gamma}_{x_{i}}\tilde{\gamma}_t\tilde{\gamma}_{x_{i+2}}\ldots\tilde{\gamma}_{x_{n}}\mpoint
\end{equation}
If $t\notin B$, this is a sum of two special monomials $M'$ with $w(M')>w(M)$, for $d_B(t)=d_B(x_i)>d_B(x_{i+1})$. And if $t\in B$, then $t=b_r$, so
\begin{equation}\label{value}
\tilde{\gamma}_t=\sum_{u\in H-\zero}\tilde{\gamma}_{b_r+u}\mvirg
\end{equation}
where $H$ is the subgroup of index 2 of $G$ generated by $B-\{b_r\}$. But $d_B(b_r+u)\geq r$, for any $u\in H$. Replacing $\tilde{\gamma}_t$ in~\ref{sum} by the right hand side of~\ref{value} gives an expression of $S$ as a sum of special monomials $M'$ with $w(M')>w(M)$.\findemo

\begin{mth}{Proposition} \label{basis}The set of special ordered monomials is a $k$-basis of $\tilde{\mathcal{E}}$, and the map $s$ is an isomorphism.
\end{mth}
\pf This is equivalent to saying that for any $n\in\N$, the set of special ordered monomials of degree $2n$ is a $k$-basis of $\tilde{\mathcal{E}}_{2n}$, and that the map $s_{2n}$ is an isomorphism.\par
The first step consists in showing that the special ordered monomials of degree $2n$ generate $\tilde{\mathcal{E}}_{2n}$. In other words, using~(\ref{value}), any special monomial
$$M=\tilde{\gamma}_{x_1}\tilde{\gamma}_{x_2}\ldots \tilde{\gamma}_{x_n}$$
should be equal to a linear combination of special ordered monomials of degree $2n$. \par
For such an arbitrary monomial $w(M)\leq nm$, and this allows for a proof by induction on $j=nm-w(M)$~: if $j=0$, then $d_B(x_i)=m$, for $1\leq\nolinebreak i\leq\nolinebreak n$, and the monomial $M$ is a special ordered monomial, so there is nothing to prove. \par
Otherwise, if $M$ is not ordered, there is a least integer $i\in\{1,\ldots,n-1\}$ such that $d_B(x_i)>d_B(x_{i+1})$. By Lemma~\ref{exchange}, the monomial $M$ is equal to the monomial obtained by exchanging $\tilde{\gamma}_{x_i}$ and $\tilde{\gamma}_{x_{i+1}}$, up to a linear combination of monomials $M'$ with $w(M')>w(M)$, which are equal to a linear combination of special ordered monomials, by induction hypothesis. \par
By repeated application of this procedure, the term $\tilde{\gamma}_{x_{i+1}}$ can be moved to the left, until it sits between $\tilde{\gamma}_{x_{j}}$ and $\tilde{\gamma}_{x_{j+1}}$ such that
$$d_B(x_j)\leq d_B(x_{i+1})< d_B(x_{j+1})\mvirg$$
(possibly $j=0$, in which case $\tilde{\gamma}_{x_{i+1}}$ is moved to the first place on the left), and the monomial $M$ is equal to the monomial
$$\tilde{\gamma}_{x_1}\tilde{\gamma}_{x_2}\ldots \tilde{\gamma}_{x_j}\tilde{\gamma}_{x_{i+1}}\tilde{\gamma}_{x_{j+1}}\ldots\tilde{\gamma}_{x_{i-1}}\tilde{\gamma}_{x_i}\tilde{\gamma}_{x_{i+2}}\ldots\tilde{\gamma}_{x_n}$$
up to a linear combination of special ordered monomials. In this monomial, the $i+1$ first values 
$$d_B(x_1),\; d_B(x_2),\;\ldots ,\; d_B(x_j),\; d_B(x_{i+1}),\; d_B(x_{j+1})\ldots,\; d_B(x_{i-1}),\; d_B(x_i)$$
are linearly ordered. By induction on $n-i$, this monomial is is equal to a linear combination of special ordered monomials. This shows that the special ordered monomials of degree $2n$ generate $\tilde{\mathcal{E}}_{2n}$.\par
The second step consists in counting the special ordered monomials~: such a monomial is a product 
$$\underbrace{\tilde{\gamma}_{x_1}\ldots\tilde{\gamma}_{x_{j_1}}}_{d_B(x_i)=1}\underbrace{\tilde{\gamma}_{x_{j_1+1}}\ldots\tilde{\gamma}_{x_{j_1+j_2}}}_{d_B(x_i)=2}\ldots\underbrace{\tilde{\gamma}_{x_{j_1+\cdots+j_{m-1}+1}}\ldots\tilde{\gamma}_{x_{j_1+\cdots+j_m}}}_{d_B(x_i)=m}$$
of $j_1$ elements $\gamma_x$ with $x\in G-(B\sqcup\zero)$ and $d_B(x)=1$, followed by $j_2$ elements $\gamma_x$ with $x\in G-(B\sqcup\zero)$ and $d_B(x)=2$, and so on, up to a product of $j_m$ elements $\gamma_x$ with $x\in G-(B\sqcup\zero)$ and $d_B(x)=m$, where $j_1+j_2+\cdots+j_m=n$. An element $x$ of $G-(B\sqcup\zero)$ such that $d_B(x)=i$ is an element of the form $y+b_i$, where $y$ is a non zero element of ${<}b_1,\ldots,b_{i-1}{>}-\zero$. It follows that the number of special ordered monomials of degree $2n$ is equal to
$$\sum_{j_1+j_2+\cdots+j_m=n}\prod_{i=1}^m(2^{i-1}-1)^{j_i}\mpoint$$
The only element $x$ in $G-\zero$ with $d_B(x)=1$ is $b_1\in B$, and this forces $j_1=0$. So the number of ordered monomials of degree $2n$ is equal to
$$\sum_{l_1+l_2+\cdots+l_{m-1}=n}\prod_{i=1}^{m-1}(2^{i}-1)^{l_i}\mvirg$$
where $l_i=j_{i+1}$, for $1\leq i\leq m-1$. But this is precisely equal to the coefficient of the term of degree $2n$ in the Poincar\'e series for $\mathcal{E}$, by Corollary~\ref{Ext generators}. It follows that $\dim_k\tilde{\mathcal{E}}_{2n}\leq \dim_k\mathcal{E}_{2n}$, hence $\dim_k\tilde{\mathcal{E}}_{2n}=\dim_k\mathcal{E}_{2n}$, and the map $s_{2n}$ is an isomorphism. This completes the proof of Corollary~\ref{basis}, and also the proof of Theorem~\ref{presentation}.\findemo
\section{Partial results for $G\cong (C_p)^m$, $p>2$}\label{algebra p odd}
Let $G\cong (C_p)^m$ be an odd order elementary abelian $p$-group. The main difference with the case $p=2$ is that Theorem~\ref{odd Ext} no longer holds, as can be seen from Theorem~\ref{ext1s1s1}~: if $\varphi:G\to k^+$ is a group homomorphism, let $E_\varphi^G$ denote the $kG$-module $k\oplus k$, where the $G$-action is defined by
$$\forall g\in G,\;\;\forall (x,y)\in k^2,\;\;g(x,y)=\big(x+y\varphi(g),y\big)\mvirg$$
and let $T_\varphi^G$ denote the unique Mackey functor for $G$ over $k$ such that $T_\varphi(H)$ is equal to zero if $\un\neq H\leq G$, and such that $T_\varphi(\un)\cong E_\varphi$. Then $T_\varphi^G$ is a cohomological Mackey functor, and there is a non split exact sequence
$$0\To S_\un^G\To T_\varphi^G\To S_\un^G\To 0\mvirg$$
whose class is an element $\tau_\varphi^G\in \Ext^1_{\comack_k(G)}(S_\un^G,S_\un^G)$.\par
If $X$ is a subgroup of order $p$ of $G$, recall that there is a 2-fold extension $\gamma_X^G\in\Ext^2_{\comack_k(G)}(S_\un^G,S_\un^G)$ represented by the short exact sequence
$$0\To S_\un^G\To\bin{X}{\un}{G}\To\bin{\un}{X}{G}\To S_\un^G\To 0\mpoint$$
In this case, I propose the following conjecture~:
\begin{mth}{Conjecture} \label{conjecture}Let $k$ be a field of odd characteristic $p$, and $G\cong (C_p)^m$. Then~:
\begin{enumerate} 
\item The algebra $\mathcal{E}=\Ext^*_{\comack_k(G)}(S_\un^G,S_\un^G)$ is generated by the elements $\tau_\varphi^G$ in degree 1, for $\varphi\in \Hom_\Z(G,k^+)$, and by the elements $\gamma_X^G$ in degree~2, for $X\leq G$ with $|X|=p$.
\item The Poincar\'e series for $\mathcal{E}$ is equal to
$$\frac{1}{(1-t)\big(1-t-(p\!-\!1)t^2\big)\big(1-t-(p^2\!-\!1)t^2\big)\ldots\big(1-t-(p^{m-1}\!-\!1)t^2\big)}\mpoint$$
\end{enumerate}
\end{mth}
\begin{mth}{Theorem} \label{p=3}Conjecture~\ref{conjecture} is true for $p=3$.
\end{mth}
\pf If $G=\un$, i.e. if $m=0$, there is nothing to prove. By induction on~$m$, one can assume that the result holds for any elementary abelian $p$-group of rank smaller than $m$. Let $H$ be a subgroup of index $p$ of $G$, let $T$ be a complement of $H$ in $G$, and set $\mathcal{X}=K_G(H)-\{T\}$. By Corollary~\ref{long exact sequence}, there is a long exact sequence in $\comack_k(G)$ of the form
\begin{equation}\label{long}
\xymatrix{
\cdots\ar[r]&*!U(0.6){L(n-1)\oplus\dirsum{X\in\mathcal{X}}E_G(n-2)}\ar[r]&E_G(n)\ar[r]^-{r_n}&E_H(n)\ar`r[d]`[l]`[dlll]`[dll][dll]&\\
&*!U(0.6){L(n)\oplus\dirsum{X\in \mathcal{X}}E_G(n-1)}\ar[r]&E_G(n+1)\ar[r]^-{r_{n+1}}&E_H(n+1)\ar[r]& \cdots\;,\\
}
\end{equation}
where $E_G(n)=\Ext^n_{\comack_k(G)}(S_\un^G,S_\un^G)$, $E_H(n)=\Ext^n_{\comack_k(H)}(S_\un^H,S_\un^H)$, $L(n)=\Ext^n_{\comack_k(G)}(L,S_\un^G)$, and $\mathcal{X}=K_G(H)-\{T\}$. Recall from Proposition~\ref{index p} that $L$ is a functor all of whose composition factors are isomorphic to $S_\un^G$, with multiplicity $p-2$. Thus if $p=3$, the functor $L$ is isomorphic to $S_\un^G$, and it follows that $L(n)\cong E_G(n)$, for any $n\in\N$.\par
It is easy to see that the map $L(n-1)\to E_G(n)$ in Sequence~\ref{long} consists in taking Yoneda product with the sequence
\begin{equation}
\label{short}
0\To L\To M\To S_\un^G\To 0
\end{equation}
obtained by taking the image of the sequence
$$0\To R/J\To I/J\To S_\un^G\To 0$$
of Proposition~\ref{index p elemab} under the split surjection $R/J\to L$.\par
On the other hand, the map $E_G(n-2)\to E_G(n)$ from the component indexed by $X\in\mathcal{X}$ is equal to the Yoneda product by $\gamma_X^G$, by the argument given in the proof Theorem~\ref{ses p=2}.\par
The map $r_n$ in Sequence~\ref{long} is induced by restriction from $G$ to $H$, also by the argument used in the proof of Theorem~\ref{ses p=2}. In particular, it is compatible with the Yoneda product. By induction hypothesis, any element in $E_H(n)$ is a linear combination of Yoneda products of elements $\tau_\psi^H$ and elements $\gamma_Y^H$, where $\psi\in\Hom_\Z(H,k^+)$ and $Y\leq H$ with $|Y|=p$.\par
To show that $r_n$ is surjective, it suffices to show that for any such $\psi$, there exists $\varphi\in\Hom_\Z(G,k^+)$ such that $\tau_\psi^H=r_1(\tau_\varphi^G)$, and that for any such $Y$, there exists $X\leq G$ with $|X|=p$, such that $\gamma_Y^H=r_2(\gamma_X^G)$. For the latter, the argument of the proof of Corollary~\ref{Ext generators} applies, and one can take $X=Y$. Now if $\psi\in\Hom_\Z(H,k^+)$, there exists $\varphi\in\Hom_\Z(G,k^+)$ whose restriction to $H$ is equal to $\psi$, and it is straightforward to check that the restriction to $H$ of the extension $T_\varphi^G$ defining $\tau_\varphi^G$ is isomorphic to $T_\psi^H$. It follows that $r_1(\tau_\varphi^G)=\tau_\psi^H$, so $r_n$ is surjective, for any $n\in\N$. Actually this proves more~: the submodule $E'_G(n)$ of $E_G(n)$ generated by products of elements $\tau_\varphi^G$ and~$\gamma_X^G$ maps surjectively by $r_n$ on $E_H(n)$. Thus $E'_G(n)+\Ker\,r_n=E_G(n)$.\par
Finally, in the case $p=3$, the long exact sequence~\ref{long} splits as a series of short ones
\begin{equation}
\label{shorts}
0\To E_G(n-1)\oplus\dirsum{X\in\mathcal{X}}E_G(n-2)\To E_G(n)\stackrel{r_n}{\To}E_H(n)\To 0\mpoint
\end{equation}
The image of each component $E_G(n-2)$ in $E_G(n)$ is obtained by taking Yoneda product with some $\gamma_X^G$, and the image of $E_G(n-1)$ is obtained by taking the Yoneda product with the sequence~\ref{short}, which for $p=3$, is of the form
$$0\To S_\un^G\To M\To S_\un^G\To 0\mpoint$$  
It follows easily that it is isomorphic to a sequence $T_\varphi$, where $\varphi\in\Hom_\Z(G,k^+)$ has kernel $H$. So the image of $E_G(n-1)$ in $E_G(n)$ is obtained by taking Yoneda product with $\tau_\varphi^G$.\par
An easy induction argument now shows that the kernel of $r_n$ is contained in $E'_G(n)$, so $E'_G(n)=E_G(n)$, completing the inductive step for Assertion~1 of Conjecture~\ref{conjecture}.\par
Assertion~2 now follows from Sequence~\ref{shorts}~: this sequence shows that if $P_m(t)$ denotes the Poincar\'e series for the algebra $\Ext^*_{\comack_k(G)}(S_\un^G,S_\un^G)$, where $G\cong(C_p)^m$, then
$$P_m(t)=tP_m(t)+|\mathcal{X}|t^2P_m(t)+P_{m-1}(t)\mpoint$$
Thus
$$P_m(t)=\frac{P_{m-1}(t)}{1-t-\big(p^{m-1}-1)t^2}\mvirg$$
and Assertion~2 follows easily by induction.\findemo
\begin{rem}{Remark} The main reason for proposing Conjecture~\ref{conjecture}, apart from Theorem~\ref{p=3}, is a computer calculation for $p=5$ or $p=7$, using GAP software ({\tt http://www.gap-system.org}), showing that the terms of lower degree of the Poincar\'e series for the algebra $\Ext^*_{\comack_k(G)}(S_\un^G,S_\un^G)$ (up to degree~5 for $p=5$ and $m=2$) are as predicted by Assertion~2 of Conjecture~\ref{conjecture}.
\end{rem}
\section{More on extension of simple functors}\label{more extensions}
\begin{mth}{Proposition} \label{centralize each other}Let $k$ be a field of characteristic $p$, and $G$ be a finite $p$-group. Let $Q$ and~$R$ be normal subgroups of $G$, and set $N=[Q,R]$. Then for each $n\in\N$, the functor $\rho_{G/N}^G$ induces an isomorphism 
$$\Ext^n_{\comack_k(G)}(S_Q^G,S_R^G)\stackrel{\cong}{\To} \Ext^{n}_{\comack_k(G/N)}(S_{Q/N}^{G/N},S_{R/N}^{G/N})\mpoint$$
\end{mth}
\pf By induction on the order of $N$~: if $N$ is trivial, there is nothing to show. Otherwise, the subgroup $N$ is a non-trivial normal subgroup of $G$, contained in $Q\cap R$, so $N$ contains a central subgroup $Z$ of order $p$ of $G$. If $X$ is a complement of $Z$ in $Q$, then for any $n\in \N$
\begin{eqnarray*}
\Ext^n_{\comack_k(G)}(S_X^G,S_R^G)&\cong&\Ext^n_{\comack_k(G)}(\Ind_{N_G(X)}^GS_X^{N_G(X)},S_R^G)\\
&\cong&\Ext^n_{\comack_k(N_G(X))}(S_X^{N_G(X)},\Res_{N_G(X)}^GS_R^G)\mvirg
\end{eqnarray*}
and $\Res_{N_G(X)}^GS_R^G=\zero$ if $R\not\leq N_G(X)$, by Lemma~\ref{rho simple}. But if $R\leq N_G(X)$, then $[R,X]\leq X$, hence $[R,Q]\leq X$, for $Q=X\cdot Z$ and $Z\leq Z(G)$. It follows that $Z\leq X$, and this contradiction shows that $\Ext^n_{\comack_k(G)}(S_X^G,S_R^G)=\zero$ for any $X\in K_Q(Z)$. Then by Theorem~\ref{quotient by central}, for any $n\in\N$, the map
$$\pi_n:\Ext^{n}_{\comack_k(G)}(S_Q^G,S_R^G){\To} \Ext^{n}_{\comack_k(\sur{G})}(S_{\sur{Q}}^{\sur{G}},S_{\sur{R}}^{\sur{G}})$$
induced by $\rho_{G/Z}^G$ is an isomorphism. Now $\sur{Q}$ and $\sur{R}$ are normal subgroups of $\sur{G}$, and $[\sur{Q},\sur{R}]=[Q,R]/Z$ has order smaller that $|N|$. By induction, the functor $\rho_{G/N}^{G/Z}$ induces an isomorphism
$$\Ext^{n}_{\comack_k(\sur{G})}(S_{\sur{Q}}^{\sur{G}},S_{\sur{R}}^{\sur{G}})\To\Ext^{n}_{\comack_k(G/N)}(S_{Q/N}^{G/N},S_{R/N}^{G/N})\mvirg$$
and the result follows by composition with $\pi_n$, since $\rho_{G/N}^{G/Z}\circ\rho_{G/Z}^G\cong\rho_{G/N}^G$.\findemo
The following result shows that for a $p$-group $G$ and a field $k$ of characteristic $p$, the computation of extension groups $\Ext^j_{\comack_k(G)}(S_Q^G,S_R^G)$ comes down to the case where $Q$ and $R$ are elementary abelian normal subgroups of $G$, which centralize each other. First a notation~:
\begin{mth}{Notation} Let $G$ be a finite $p$-group. If $Q$ and $R$ are subgroups of~$G$, set
$$\mathcal{S}_{R,Q}^G=\{g\in G\mid \Phi(R)\cdot{^g\Phi(Q)}\cdot[R,{^gQ}]\leq R\cap{^gQ}\}\mpoint$$
Denote by $\sur{\mathcal{S}}_{R,Q}^G$ the set of $\big(N_G(R),N_G(Q)\big)$-double cosets in $\mathcal{S}_{R,Q}^G$, and by $[\sur{\mathcal{S}}_{R,Q}^G]$ a set of representatives of $\sur{\mathcal{S}}_{R,Q}^G$.\par
Let $\sur{\N}=\N\sqcup\{\infty\}$ be the linearly ordered set obtained by adding to $\N$ a largest element~$\infty$. Denote by $\nu_G(R,Q)\in \sur{\N}$ the element defined by
$$p^{\nu_G(R,Q)}=\min_{g\in\mathcal{S}_{R,Q}^G}\limits|R{\cdot}{^gQ}:R\cap{^gQ}|=\min_{g\in[\sur{\mathcal{S}}_{R,Q}^G]}\limits|R{\cdot}{^gQ}:R\cap{^gQ}|\mvirg$$
if $\mathcal{S}_{R,Q}^G\neq\emptyset$, and by $\nu_G(R,Q)=\infty$ otherwise.
\end{mth}
\begin{mth}{Theorem} \label{general Ext}Let $k$ be a field of characteristic $p$, let $G$ be a finite $p$-group, and let $Q$ and $R$ be subgroups of~$G$. Then for $g\in\mathcal{S}_{R,Q}^G$~:
\begin{enumerate}
\item the groups $R$ and $^gQ$ normalize each other.
\item the group $\Phi(R)\cdot{^g\Phi(Q)}\cdot[R,{^gQ}]$ is equal to $\Phi(R{\cdot}{^gQ})$.
\item the groups $\widehat{R}=R/\Phi(R{\cdot}{^gQ})$ and $\widehat{^gQ}={^gQ}/\Phi(R{\cdot}{^gQ})$ are elementary abelian normal subgroups of $\widehat{N}_G(R,{^gQ})=N_G(R,{^gQ})/\Phi(R{\cdot}{^gQ})$, which centralize each other. 
\end{enumerate}
Moreover for any $j\in\N$,
$$\Ext^j_{\comack_k(G)}(S_Q,S_R)\cong\dirsum{g\in[\sur{\mathcal{S}}_{R,Q}^G]}\Ext^j_{\comack_k(\widehat{N}_G(R,{^gQ}))}\big(S_{\widehat{^gQ}}^{\widehat{N}_G(R,{^gQ})},S_{\widehat{R}}^{\widehat{N}_G(R,{^gQ})}\big)\mvirg$$
\end{mth}
\pf If $g\in\mathcal{S}_{R,Q}^G$, then $[R,{^gQ}]\leq R\cap{^gQ}$, so $R$ and $^gQ$ normalize each other. Let $T=\Phi(R)\cdot{^g\Phi(Q)}\cdot[R,{^gQ}]$. Then the groups $R/T$ and $^gQ/T$ are elementary abelian subgroups of $N_G(R,{^gQ})/T$, which centralize each other. Thus $(R{\cdot}{^gQ})/T$ is elementary abelian, and $T\geq \Phi(R{\cdot}{^gQ})$. \par
Conversely $\Phi(R)\leq\Phi(R{\cdot}{^gQ})$, since $R{\cdot}{^gQ}$ is a $p$-group, and similarly $\Phi({^gQ})={^g\Phi(Q)}\leq \Phi(R{\cdot}{^gQ})$. Moreover $[R,{^gQ}]\leq[R{\cdot}{^gQ},R{\cdot}{^gQ}]\leq \Phi(R{\cdot}{^gQ})$, so finally $T\leq \Phi(R{\cdot}{^gQ})$, hence $T=\Phi(R{\cdot}{^gQ})$. Now $\widehat{R}=R/T$ and $\widehat{^gQ}={^gQ}/T$ are elementary abelian normal subgroups of $\widehat{N}_G(R,{^gQ})$, which centralize each other.\par
Set $A=\Phi(Q)$ and $B=N_G(Q)$. Set similarly $C=\Phi(R)$ and $D=N_G(R)$. Then $S_Q^G\cong \Ind_B^G\iiota_{B/A}^B(S_{Q/A}^{B/A})$ and $S_R^G=\Ind_{D}^G\jiota_{D/C}^D(S_{R/C}^{D/C})$, by Lemma~\ref{simple from elemab} and Corollary~\ref{simple from elemab 2}. Thus, by Proposition~\ref{ext adjunction}, for any $j\in\N$,
$$\Ext^j_{\comack_k(G)}(S_Q,S_R)\cong\Ext^j_{\comack_k(D/C)}\big(\rho_{D/C}^D\Res_D^G\Ind_B^G\iiota_{B/A}^B(S_{Q/A}^{B/A}),S_{R/C}^{D/C}\big)\mpoint
$$
Now by Proposition~\ref{general Mackey}, 
$$\rho_{D/C}^D\Res_D^G\Ind_B^G\iiota_{B/A}^B(S_{Q/A}^{B/A})\!\cong\!\dirsum{g\in[D\dom G/B]}\Ind_{\sur{D}_g}^{\sur{D}}\iiota_{\sur{D}_g/\sur{C}_g}^{\sur{D}_g}\Iso(f_g)\rho_{\sur{B}_g/\sur{A}_g}^{\sur{B}_g}\Res_{\sur{B}_g}^{\sur{B}}(S_{Q/A}^{B/A})\mvirg$$
where 
$$\sur{D}=D/C, \;\;\sur{D}_g=(D\cap{^gB})C/C, \;\;\sur{C}_g=(D\cap{^gA})C/C\mvirg$$
$$\sur{B}=B/A, \;\;\sur{B}_g=(D^g\cap B)A/A, \;\;\sur{A}_g=(C^g\cap B)A/A\mvirg$$
and where $f_g:\sur{B}_g/\sur{A}_g\To\sur{D}_g/\sur{C}_g$ is the group isomorphism sending $x\sur{A}_g$ to ${^gx}\sur{C}_g$, for $x\in D^g\cap B$.\par
By Lemma~\ref{rho simple}, the functor $\rho_{\sur{B}_g/\sur{A}_g}^{\sur{B}_g}\Res_{\sur{B}_g}^{\sur{B}}(S_{Q/A}^{B/A})$ is isomorphic to the direct sum of the simple functors $S_{X/\sur{A}_g}^{\sur{B}_g/\sur{A}_g}$ corresponding to subgroup $X$ such that $X/A$ is conjugate to $\sur{Q}=Q/A$ in $B/A$, up to conjugation by $\sur{B}_g/\sur{A}_g$. Since $Q\normal B$, the only possible such subgroup is $Q$ itself, if $\sur{A}_g\leq \sur{Q}\leq \sur{B}_g$. This gives
$$\Ext^j_{\comack_k(G)}(S_Q,S_R)\cong\oplusb{g\in[D\dom G/B]}{\sur{A}_g\leq \sur{Q}\leq \sur{B}_g}\Ext^j_{\comack_k(\sur{D})}\big(\Ind_{\sur{D}_g}^{\sur{D}}\iiota_{\sur{D}_g/\sur{C}_g}^{\sur{D}_g}\Iso(f_g)(S_{\sur{Q}/\sur{A}_g}^{\sur{B}_g/\sur{A}_g}),S_{R/C}^{D/C}\big)\mpoint$$
By adjunction, this gives
$$\Ext^j_{\comack_k(G)}(S_Q,S_R)\!\cong\!\oplusb{g\in[D\dom G/B]}{\sur{A}_g\leq \sur{Q}\leq \sur{B}_g}\Ext^j_{\comack_k(\sur{D})}\big(\Iso(f_g)(S_{\sur{Q}/\sur{A}_g}^{\sur{B}_g/\sur{A}_g}),\rho_{\sur{D}_g/\sur{C}_g}^{\sur{D}_g}\Res_{\sur{D}_g}^{\sur{D}}(S_{\sur{R}}^{\sur{D}})\big)\mvirg$$
where $\sur{R}=R/C$. By the same argument, the functor $\rho_{\sur{D}_g/\sur{C}_g}^{\sur{D}_g}\Res_{\sur{D}_g}^{\sur{D}}(S_{\sur{R}}^{\sur{D}})$ is equal to zero, unless $\sur{C}_g\leq \sur{R}\leq \sur{D}_g$, in which case it is isomorphic to $S_{\sur{R}/\sur{C}_g}^{\sur{D}_g/\sur{C}_g}$.  Thus
$$\Ext^j_{\comack_k(G)}(S_Q,S_R)\!\cong\!\oplusc{g\in[D\dom G/B]}{\sur{A}_g\leq \sur{Q}\leq \sur{B}_g}{\sur{C}_g\leq \sur{R}\leq \sur{D}_g}\Ext^j_{\comack_k(\sur{D}_g/\sur{C}_g)}\big(\Iso(f_g)(S_{\sur{Q}/\sur{A}_g}^{\sur{B}_g/\sur{A}_g}),S_{\sur{R}/\sur{C}_g}^{\sur{D}_g/\sur{C}_g}\big)\mpoint$$
Now the condition $\sur{Q}\leq \sur{B}_g$ is equivalent to $Q\leq (D^g\cap B)A$, i.e. to $Q=A(D^g\cap Q)$, i.e. to $D^g\cap Q=Q$, since $A=\Phi(Q)$. In other words $Q\leq D^g$. Similarly, the condition $\sur{R}\leq \sur{D}_g$ is equivalent to $R\leq (D\cap {^gB})C$, i.e. to $R=(R\cap {^gB})C$, i.e. to $R=R\cap{^gB}$, since $C=\Phi(R)$. In other words $R\leq {^gB}$.\par
Then the condition $\sur{A}_g\leq \sur{Q}$ is equivalent to $C^g\cap B\leq Q$. Since moreover $C^g\leq R^g\leq B$, it follows that $C^g\cap B=C^g\leq Q$. \par
Similarly, the condition $\sur{C}_g\leq \sur{R}$ is equivalent to $D\cap{^gA}\leq R$. But ${^gA}\leq{^gQ}\leq D$, thus $D\cap{^gA}={^gA}\leq R$.\par
In this situation 
$$\sur{D}_g=(D\cap{^gB})/C,\;\;\sur{C}_g={^gA}{\cdot}C/C,\;\;\sur{B}_g=(D^g\cap B)/A,\;\;\sur{A}_g=C^g{\cdot}A/A\mvirg$$
and the isomorphism 
$$f_g:\sur{B}_g/\sur{A}_g=(D^g\cap B)/(C^g{\cdot}A)\To (D\cap {^gB})/(C{\cdot}{^gA})=\sur{D}_g/\sur{C}_g$$
is induced by conjugation by $g$. In particular $f_g(\sur{Q}/\sur{A}_g)={^gQ}/(C{\cdot} {^gA})$, thus finally
$$\Ext^j_{\comack_k(G)}(S_Q,S_R)\cong\dirsum{g\in \mathcal{S} }\Ext^j_{\comack_k((D\cap{^gB})/(C{\cdot}{^gA}))}\big(S_{^gQ/(C{\cdot} {^gA})}^{(D\cap{^gB})/(C{\cdot}{^gA})},S_{R/(C{\cdot} {^gA})}^{(D\cap{^gB})/(C{\cdot}{^gA})}\big)\mvirg$$
where 
\begin{eqnarray*}
 \mathcal{S} &=&\{g\in[D\dom G/B]\mid C{\cdot}{^gA}\leq R\cap{^gQ}\leq R{\cdot}^gQ\leq D\cap{^gB}\}\\
&=&\{g\in[D\dom G/B]\mid C{\cdot}{^gA}\cdot[R,{^gQ}]\leq R\cap{^gQ}\}=[\sur{\mathcal{S}}_{R,Q}^G]\mpoint
\end{eqnarray*}
Moreover, by Proposition~\ref{centralize each other}, for $g\in\mathcal{S}$, the group 
$$\Ext^j_{\comack_k((D\cap{^gB})/(C{\cdot}{^gA}))}\big(S_{^gQ/(C{\cdot} {^gA})}^{(D\cap{^gB})/(C{\cdot}{^gA})},S_{R/(C{\cdot} {^gA})}^{(D\cap{^gB})/(C{\cdot}{^gA})}\big)$$
is isomorphic to
$$\Ext^j_{\comack_k(\widehat{N}_G(R,{^gQ}))}\big(S_{\widehat{^gQ}}^{\widehat{N}_G(R,{^gQ})},S_{\widehat{R}}^{\widehat{N}_G(R,{^gQ})}\big)\mvirg$$
and this completes the proof.\findemo
\begin{mth}{Proposition} \label{Ext zero}Let $k$ be a field of characteristic $p$, let $G$ be a finite $p$-group, and let $Q$ and $R$ be subgroups of $G$. Then
$$n<\nu_G(R,Q)\;\Rightarrow\;\Ext^n_{\comack_k(G)}(S_Q^G,S_R^G)=\zero\mpoint$$
\end{mth}
\pf If $\nu_G(R,Q)=\infty$, i.e. if $\mathcal{S}^G_{R,Q}=\emptyset$, then $\Ext^n_{\comack_k(G)}(S_Q^G,S_R^G)=\zero$ for any $n\in\N$, by Theorem~\ref{general Ext}. So one can assume $\nu_G(R,Q)\in\N$.\par
The result is trivial if $G=\un$, so by induction on $|G|$, one can assume that it holds for any group of order less than $|G|$. The group $G$ being given with this property, one can proceed by induction on $n$, to show that $\Ext^n_{\comack_k(G)}(S_Q^G,S_R^G)=\zero$ if $Q$ and $R$ are subgroups of $G$ with $\nu_G(R,Q)>n$. This starts with the case $n=-1$, where there is nothing to prove, since $\nu_G(R,Q)\geq 0$ and $\Ext^n_{\comack_k(G)}(S_Q^G,S_R^G)=\zero$ for $n<0$.\par
For the inductive step, let $n\geq 0$, let $Q$ and $R$ be subgroups of $G$ such that $\nu_G(R,Q)>n$, and let $g\in \mathcal{S}_{R,Q}^G$. Then the index $|R{\cdot}{^gQ}:R\cap{^gQ}|$ is equal to the index $|R'{\cdot}Q':R'\cap Q'|$, where $R'=R/\big(\Phi(R){\cdot}{^g\Phi(Q)}\big)$ and $Q'={^gQ}/\big(\Phi(R){\cdot}{^g\Phi(Q)}\big)$. By Theorem~\ref{general Ext}, and by induction hypothesis on~$G$, it suffices to consider the case where $Q$ and $R$ are elementary abelian normal subgroups of $G$, which centralize each other. \par
In this case, set $j=\nu_G(R,Q)$, so $p^j=|R{\cdot}Q:R\cap Q|$. If $Q\cap R\neq \un$, let $Z$ be a subgroup of order $p$ of $Q\cap R\cap Z(G)$. By Theorem~\ref{quotient by central}, for any integer~$n$, there is a short exact sequence
$$0\to \dirsum{X\in\mathcal{K}}\Ext^{n-1}_{\comack_k(G)}(S_X^G,S_R^G)\to \Ext^{n}_{\comack_k(G)}(S_Q^G,S_R^G)\to \Ext^{n}_{\comack_k(\sur{G})}(S_{\sur{Q}}^{\sur{G}},S_{\sur{R}}^{\sur{G}})\to 0\mvirg$$
where $\mathcal{K}$ is a set of representatives of $G$-conjugacy classes of complements of~$Z$ in $Q$, and $\sur{Q}$ and $\sur{R}$ denote $Q/Z$ and $R/Z$, respectively.\par
If $X\in\mathcal{K}$, since $R\normal G$, the set $N_G(X)\dom G/N_G(R)$ has
cardinality~1. Moreover $X\cap R\neq Q\cap R$, since $Z\leq Q\cap R$ and
$Z\not\leq X$, so $|Q\cap R:X\cap R|=p$, and
$$|X{\cdot}R:X\cap R|=\frac{|Q||R|}{p|X\cap R|^2}=p^{j+1}\mpoint$$
In other words $\nu_G(R,X)=j+1>n-1$, so $\Ext^{n-1}_{\comack_k(G)}(S_X^G,S_R^G)=\zero$ by induction hypothesis on $n$. \par
On the other hand $\Ext^{n}_{\comack_k(\sur{G})}(S_{\sur{Q}}^{\sur{G}},S_{\sur{R}}^{\sur{G}})=0$ since $n<j$, by induction hypothesis on~$|G|$. It follows that $\Ext^{n}_{\comack_k(G)}(S_Q^G,S_R^G)=\zero$.\par
It remains to consider the case where $Q\cap R=\un$. Since $\nu_G(R,Q)>n\geq0$, it follows that at least one of the groups $Q$ or $R$ is non trivial.
By symmetry on $Q$ and $R$, one can assume that $Q\neq \un$, and choose a subgroup $Z$ of order~$p$ in $Q\cap Z(G)$. \par
Now for any $n\in\N$
$$\Ext^n_{\comack_k(G)}\big(\iiota_{G/Z}^G(S_{Q/Z}^{G/Z}),S_R^Q\big)\cong\Ext^n_{\comack_k(G/Z)}\big(S_{Q/Z}^{G/Z},\rho_{G/Z}^G(S_R^G)\big)=\zero\mvirg$$
by Proposition~\ref{rho simple}, since $Z\not\subseteq R$. Applying $\Hom_{\comack_k(G)}({-},S_R^G)$ to the first exact sequence of Proposition~\ref{iota S_Q} gives the isomorphisms
\begin{equation}
\label{iso}
\Ext^n_{\comack_k(G)}(S_Q^G,S_R^G)\cong\dirsum{X\in\mathcal{K}}\Ext^{n-1}_{\comack_k(G)}(S_X^G,S_R^G)\mvirg
\end{equation}
for any $n\in \N$, where $\mathcal{K}$ is a set of representatives of $G$-conjugacy classes of complements of $Z$ in $Q$.\par
Let $X\in\mathcal{K}$, and set $l=\nu_{G}(R,X)$. Then
$$p^l=|X{\cdot}R:X\cap R|=|X{\cdot}R|=|X||R|\mvirg$$
so $l=j-1$. \par
If $n<j$, then $n-1<l$, hence $\Ext^{n-1}_{\comack_k(G)}(S_X^{G},S_R^{G})=\zero$ by induction hypothesis on $n$.
It follows that $\Ext^n_{\comack_k(G)}(S_Q^G,S_R^G)=\zero$, as was to be shown. \findemo
\begin{mth}{Notation} Let $Q$ be an elementary abelian $p$-group. Denote by~${\rm St}(Q)$ the only non-zero reduced integral homology group of the poset $]\un,Q[$ of proper non trivial subgroups of $Q$, ordered by inclusion.
\end{mth}
It is well known (see e.g. \cite{homol}) that, if $Q$ has rank $r$, then ${\rm St}(Q)=\tilde{H}_{r-2}(]\un,Q[,\Z)$ is isomorphic to $\Z^{p^{\binom{r}{2}}}$. The automorphism group $A\cong {\rm GL}(r,\F_p)$ of $Q$ acts on ${\rm St}(Q)$, and the module ${\rm St}(Q)$ is isomorphic to {\em the Steinberg module} of $A$. The restriction of ${\rm St}(Q)$ to a Sylow $p$-subgroup $P$ of $A$ is a free $\Z P$-module of rank~1.\par
When $G$ is an elementary abelian $p$-group, the following theorem was proved by Tambara (\cite{tambara} Theorem~4.1)~:
\begin{mth}{Theorem} \label{Steinberg}Let $k$ be a field of characteristic $p$, and $G$ be a finite $p$-group. Let $Q$ be an elementary abelian normal subgroup of rank $q$ of $G$, and $H$ be a subgroup of rank $h$ of $Q$. Then $\nu_G(H,Q)=q-h$ and~:
\begin{enumerate}
\item There is a group isomorphism 
$$\Ext^{q-h}_{\comack_k(G)}(S_Q^G,S_H^G)\cong\Ext^{q-h}_{\comack_k(N_G(H)/H)}(S_{Q/H}^{N_G(H)/H},S_{H/H}^{N_G(H)/H})\mpoint$$
\item These groups are isomorphic to the space of coinvariants $k{\rm St}(Q/H)_{N_G(H)}$ of $N_G(H)$ on $k\otimes_\Z{\rm St}(Q/H)$. 
\item Let $H=H_0<H_1<\ldots<H_{q-h}=Q$ be a maximal $N_G(H)$-invariant flag. Then $\Ext^{q-h}_{\comack_k(G)}(S_Q^G,S_H^G)$ has a $k$-basis indexed by the $N_G(H)$-conjugacy classes of flags $X_1>X_2>\ldots>X_{q-h-1}$ of $Q$ such that $X_i+ H_i=Q$ and $X_i\cap H_i=H$, for $1\leq i\leq q-h-1$. In particular
$$\dim_k\Ext^{q-h}_{\comack_k(G)}(S_Q^G,S_H^G)=\frac{p^{\binom{q-h}{2}}}{|N_G(H):Z_G(H,Q)|}\mvirg$$
where $Z_G(H,Q)=\{g\in G\mid [g,Q]\leq H\}$ (in other words $Z_G(H,Q)$ is the preimage in $N_G(H)$ of the centralizer of $Q/H$ in $N_G(H)/H$).
\end{enumerate}
\end{mth}
\pf Since $H\leq Q$ and $\Phi(Q)=\un$, it follows that $\Phi(H)\Phi(Q)[H,Q]=\un\leq H\cap Q$, so $\nu_G(H,Q)<\infty$. Moreover $p^{\nu_G(H,Q)}=|Q{\cdot}H:Q\cap H|=|Q:H|$, since $Q\normal G$, so $\nu_G(H,Q)=q-h$.\par
For Assertion~1, observe first that 
$$\Ext^n_{\comack_k(G)}(S_Q^G,S_H^G)\cong \Ext^n_{\comack_k(N_G(H))}(S_Q^{N_G(H)},S_H^{N_G(H)})\mvirg$$
for any $n\in\N$, and $\nu_{N_G(H)}(H,Q)=q-h$ since $Q\normal N_G(H)$. By induction on~$|G|$, if $N_G(H)<G$, there is an isomorphism
$$\Ext^n_{\comack_k(N_G(H))}(S_Q^{N_G(H)},S_H^{N_G(H)})\cong \Ext^n_{\comack_k(N_G(H)/H)}(S_{Q/H}^{N_G(H)/H},S_{H/H}^{N_G(H)/H})\mvirg$$
so one can suppose $H\normal G$.\par
If $H=\un$, there is nothing to prove. Otherwise, let $Z$ be a subgroup of order $p$ of $H\cap Z(G)$. 
By Theorem~\ref{quotient by central}, there is an exact sequence 
$$0\to \dirsum{X\in\mathcal{K}}\Ext^{q-h-1}_{\comack_k(G)}(S_X^G,S_H^G)\to \Ext^{q-h}_{\comack_k(G)}(S_Q^G,S_H^G)\stackrel{\pi_{q-h}}{\to} \Ext^{q-h}_{\comack_k(\sur{G})}(S_{\sur{Q}}^{\sur{G}},S_{\sur{H}}^{\sur{G}})\to 0\mvirg$$
in $\comack_k(G)$, where $\mathcal{K}$ is a set of representatives of $G$-conjugacy classes of complements of $Z$ in $Q$, and overlines denote quotients by $Z$. If $X\in\mathcal{K}$, then
$$p^{\nu_G(H,X)}=|X{\cdot} H:X\cap H|=\frac{|X||H|}{|X\cap H|^2}= \frac{|Q||H|}{p|X\cap H|^2}=p^{q-h-1}|H:X\cap H|^2\mpoint$$
Since $Z\leq H$ and $Z\cap X=\un$, it follows that $H\not\leq X$, thus $|H:X\cap H|=p$, and $\nu_G(H,X)=q-h+1$. By Proposition~\ref{Ext zero}, the group $\Ext^{q-h-1}_{\comack_k(G)}(S_X^G,S_R^G)$ is equal to zero. Hence the map $\pi_{q-h}$ is an isomorphism. By induction on the the order of $G$, there is an isomorphism
$$\Ext^{q-h}_{\comack_k(\sur{G})}(S_{\sur{Q}}^{\sur{G}},S_{\sur{H}}^{\sur{G}})\cong \Ext^{q-h}_{\comack_k(G/H)}(S_{Q/H}^{G/H},S_{H/H}^{G/H})\mvirg$$
which completes the proof of Assertion~1.\par
It follows that for Assertions~2 and~3, it suffices to consider the case $H=\un$, i.e. $h=0$. In this case $H_i\normal G$, for $0\leq i\leq q$. If $q=0$, i.e. if $Q=\un$, there is nothing to prove, since $\Ext^0_{\comack_k(G)}(S_\un^G,S_\un^G)\cong k$. This allows to start a proof by induction on $q$.\par
If $q>0$, by Proposition~\ref{iota S_Q}, there is an exact sequence
$$0\To \dirsum{X_\un\in \mathcal{K}}S_{X_1}^G\To \iota_{G/H_1}^G(S_{Q/H_1}^{G/H_1})\To S_Q^G\To 0\mvirg$$
where $\mathcal{K}$ is a set of representatives of $G$-conjugacy classes of complements of~$H_1$ in $Q$. Since $\rho_{G/H_1}^G(S_\un^G)=\zero$, this gives an isomorphism
$$\Ext^{q}_{\comack_k(G)}(S_Q^G,S_\un^G) \cong \dirsum{X_1\in \mathcal{K}}\Ext^{q-1}_{\comack_k(G)}(S_{X_1}^G,S_\un^G)\mpoint$$
Moreover, for $X_\un\in\mathcal{K}$,
$$\Ext^{q-1}_{\comack_k(G)}(S_{X_1}^G,S_\un^G)\cong \Ext^{q-1}_{\comack_k(N_G(X_1))}(S_{X_1}^{N_G(X_1)},S_\un^{N_G(X_1)})\mpoint$$
The flag $\un<H_2\cap X_1<H_3\cap X_1<\ldots<H_{q-1}\cap X_1<X_1$ is a maximal $N_G(X_1)$-invariant flag in $X_1$. By induction hypothesis, the space 
$$\Ext^{q-1}_{\comack_k(N_G(X_1))}(S_{X_1}^{N_G(X_1)},S_\un^{N_G(X_1)})$$
has a $k$-basis indexed by $N_G(X_1)$-conjugacy classes of flags 
$$X_2>X_3>\ldots>X_{q-1}$$
such that $X_i\oplus (H_{i}\cap X_1)=X_1$, for $2\leq i\leq q-1$. Equivalently $X_2<X_1$, and $X_i\oplus H_{i}=Q$, for $2\leq i\leq q-1$. \par
Now $X_1$ runs through a set of representatives of $G$-conjugacy classes of complements of $H_1$ in $Q$, and $X_2>\ldots>X_{q-1}$ runs through a set of representatives of $N_G(X_1)$-conjugacy classes of flags such that $X_i\oplus H_{i}=Q$, for $2\leq i\leq q-1$. Equivalently $X_1>X_2>\ldots>X_{q-1}$ runs through a set of conjugacy classes of flags of $Q$ such that $X_i\oplus H_i=Q$, for $1\leq i\leq q-1$. This completes the inductive step, and the first part of Assertion~3 follows.\par
Now the stabilizer $C$ in $G$ of a flag $X_1>X_2>\ldots>X_{q-1}$ such that $X_i\oplus H_i=Q$, for $1\leq i\leq q-1$ stabilizes two opposite maximal flags of $Q$. Since $G$ is a $p$-group, it follows that the image of $C$ in the automorphism group of $Q$ is trivial (since it consists of matrices which are both upper and lower triangular with 1's on the diagonal). In other words $C=C_G(Q)$.\par
There are $p^{\binom{q}{2}}$ maximal flags opposite to a given maximal flag in $Q$, and the cardinality of each $G$-conjugacy class of such flags is $|G:C_G(Q)|$. This completes the proof of Assertion~3.\par
Assertion~2 follows from the fact that $k{\rm St}(Q)$ has a $k$-basis indexed by maximal flags opposite to a given $G$-invariant maximal flag of $Q$, and this basis is permuted by $G$ (see \cite{homol} for details). In other words $k{\rm St}(Q)$ is a free $k\big(G/C_G(Q)\big)$-module.\findemo
Together with Theorem~\ref{general Ext}, the following theorem allows to compute any extension group between simple functors for a $p$-group $G$, from the self extension groups for the simple functor $S_\un^C$, where $C$ runs through some subquotients of~$G$~:
\begin{mth}{Theorem} \label{factor Ext}Let $k$ be a field of characteristic $p$, let $G$ be a finite $p$-group, and let $Q$ and $R$ be elementary abelian normal subgroups of $G$ which centralize each other. Let $q$ denote the rank of $Q$ and $r$ the rank of $R$. 
\begin{enumerate}
\item For any $n\in\N$, the group $\Ext^n_{\comack_k(G)}(S_Q^G,S_R^G)$ is isomorphic to
$$\left(\dirsumc{h\in\N}{\rule{0ex}{1.5ex}H\leq Q\cap R}{\rule{0ex}{1.5ex}{\rm rank}(H)=h}k{\rm St}(R/H)\otimes_k\Ext^{n+2h-q-r}_{\comack_k(Z_H)}(S_\un^{Z_H},S_\un^{Z_H})\otimes_kk{\rm St}(Q/H)\right)_G\mvirg$$
where $Z_H=Z_G(H,Q{\cdot}R)/H=\{g\in G\mid [g,Q{\cdot}R]\leq H\}/H$, and the $G$-subscript denotes the space of coinvariants. \par
\item In particular, the dimension of the space $\Ext^n_{\comack_k(G)}(S_Q^G,S_R^G)$ is equal to
$$\sumc{h\in\N}{\rule{0ex}{1.5ex}H\leq Q\cap R,\;{\rm mod.}G}{\rule{0ex}{1.5ex}{\rm rank}(H)=h}\frac{p^{\binom{q-h}{2}+\binom{r-h}{2}}}{|N_G(H):Z_G(H,Q{\cdot}R)|}\dim_k\Ext^{n+2h-q-r}_{\comack_k(Z_H)}(S_\un^{Z_H},S_\un^{Z_H})\mpoint$$
\end{enumerate}
\end{mth}
\pf For Assertion~1, when $Q$ and $R$ are both trivial, there is nothing to prove. This allows for a proof by induction on $|Q||R|$. By symmetry on $Q$ and $R$, one can assume that $Q\neq\un$, and choose a subgroup $Z$ of order $p$ of $Q\cap Z(G)$. There are two cases~:\spn
$\bullet$ if $Z\not\leq R$, then $\rho_{G/Z}^G(S_R^G)=\zero$ by Lemma~\ref{rho simple}, and the exact sequence of Proposition~\ref{iota S_Q} gives isomorphisms
$$\Ext^n_{\comack_k(G)}(S_Q^G,S_R^G)\cong\dirsum{X\in\mathcal{K}}\Ext^{n-1}_{\comack_k(G)}(S_X^G,S_R^G)\mvirg$$
for any $n\in\N$, where $\mathcal{K}=[G\dom K_Q(Z)]$ is a set of representatives of $G$ conjugacy classes of complements of $Z$ in $Q$. 
Now for each $X\in\mathcal{K}$
$$\Ext^{n-1}_{\comack_k(G)}(S_X^G,S_R^G)\cong \Ext^{n-1}_{\comack_k(N_G(X))}(S_X^{N_G(X)},S_R^{N_G(X)})\mpoint$$
Since $|X||R|=|Q||R|/p$, by induction hypothesis, this group 
is isomorphic to
$$\left(\dirsumc{h\in \N}{H\leq X\cap R}{{\rm rank}(H)=h}k{\rm St}(R/H)\otimes_k\Ext^{n-1+2h-x-r}_{\comack_k(Z_{H}^X)}(S_\un^{Z_{H}^X},S_\un^{Z_{H}^X})\otimes_kk{\rm St}(X/H)\right)_{N_G(X)}\mvirg$$
where $x=q-1$ is the rank of $X$, and 
$$Z_{H}^X=\{g\in N_G(X)\mid [g,X{\cdot}R]\leq H\}/H\mpoint$$
But if $g\in G$, then $[g,X{\cdot}R]\leq H$ if and only if $[g,Q{\cdot}R]\leq H$, since $Q=X{\cdot}Z$ and $Z\leq Z(G)$. And if $[g,Q{\cdot}R]\leq H$, then 
$$[g,X]\leq [g,Q]\leq[g,Q{\cdot}R]\leq H\leq X\mvirg$$
so $g\in N_G(X)$. It follows that 
$$Z_H^X=\{g\in G\mid [g,Q{\cdot}R]\leq H\}/H=Z_H\mpoint$$
Finally, the group $\Ext^n_{\comack_k(G)}(S_Q^G,S_R^G)$ is isomorphic to
$$\left(\dirsum{X\in K_Q(Z)}\!\!\!\dirsumc{h\in \N}{H\leq X\cap R}{{\rm rank}(H)=h}k{\rm St}(R/H)\otimes_k\Ext^{n+2h-q-r}_{\comack_k(Z_{H})}(S_\un^{Z_{H}},S_\un^{Z_{H}})\otimes_kk{\rm St}(X/H)\right)_G\mpoint$$
Exchanging the order of summation, this gives
$$\left(\dirsumc{h\in \N}{H\leq Q\cap R}{{\rm rank}(H)=h}\!\!\!k{\rm St}(R/H)\otimes_k\Ext^{n+2h-q-r}_{\comack_k(Z_{H})}(S_\un^{Z_{H}},S_\un^{Z_{H}})\otimes_k\big(\!\!\!\!\!\dirsumb{\rule{0ex}{2ex}X\in K_Q(Z)}{X\geq H}k{\rm St}(X/H)\big)\right)_G\mpoint$$
Now $Z\cap H=\un$ since $Z\cap R=\un$, hence the conditions $X\in K_Q(Z)$ and $X\geq H$ are equivalent to $X{\cdot}(ZH)=Q$ and $X\cap (ZH)=H$, i.e. to $H\leq X\leq Q$ and $X/H\in K_{Q/H}(ZH/H)$. By classical combinatorial results,
$$\dirsum{X/H\in K_{Q/H}(ZH/H)}k{\rm St}(X/H)\cong {\rm St}(Q/H)\mvirg$$
thus $\Ext^n_{\comack_k(G)}(S_Q^G,S_R^G)$ is isomorphic to
$$\left(\dirsumc{h\in \N}{H\leq Q\cap R}{{\rm rank}(H)=h}\!\!\!k{\rm St}(R/H)\otimes_k\Ext^{n+2h-q-r}_{\comack_k(Z_{H})}(S_\un^{Z_{H}},S_\un^{Z_{H}})\otimes_kk{\rm St}(Q/H)\right)_G\mvirg$$
as was to be shown.\spn
$\bullet$ If $Z\leq R$, then by Theorem~\ref{quotient by central}, for any $n\in \N$, there is an exact sequence 
$$0\to \underbrace{\dirsum{X\in\mathcal{K}}\Ext^{n-1}_{\comack_k(G)}(S_X^G,S_R^G)}_\Sigma\to \Ext^{n}_{\comack_k(G)}(S_Q^G,S_R^G){\to} \Ext^{n}_{\comack_k(\sur{G})}(S_{\sur{Q}}^{\sur{G}},S_{\sur{R}}^{\sur{G}})\to 0\mvirg$$
in $\comack_k(G)$, where $\mathcal{K}=[G\dom K_Q(Z)]$ as above, and overlines denote quotients by $Z$. It follows that $\Ext^{n}_{\comack_k(G)}(S_Q^G,S_R^G)\cong \Sigma\oplus \Ext^{n}_{\comack_k(\sur{G})}(S_{\sur{Q}}^{\sur{G}},S_{\sur{R}}^{\sur{G}})$. \par
The same argument as above shows that the space $\Sigma$ is isomorphic to
\begin{equation}\label{part1}
\left(\dirsumc{h\in \N}{H\leq Q\cap R,\, H\cap Z=\un}{{\rm rank}(H)=h}\!\!\!k{\rm St}(R/H)\otimes_k\Ext^{n+2h-q-r}_{\comack_k(Z_{H})}(S_\un^{Z_{H}},S_\un^{Z_{H}})\otimes_kk{\rm St}(Q/H)\right)_G\mpoint
\end{equation}
On the other hand, by induction hypothesis, the space $\Ext^{n}_{\comack_k(\sur{G})}(S_{\sur{Q}}^{\sur{G}},S_{\sur{R}}^{\sur{G}})$ is isomorphic to
$$\left(\dirsumc{\sur{h}\in\N}{\sur{H}\leq \sur{Q}\cap \sur{R}}{{\rm rank}(\sur{H})=\sur{h}}k{\rm St}(\sur{R}/\sur{H})\otimes_k\Ext^{n+2\sur{h}-\sur{q}-\sur{r}}_{\comack_k(Z_{\sur{H}})}(S_\un^{Z_{\sur{H}}},S_\un^{Z_{\sur{H}}})\otimes_kk{\rm St}(\sur{Q}/\sur{H})\right)_{\sur{G}}\mvirg$$
where $\sur{q}=q-1$ and $\sur{r}=r-1$. Summing over $\sur{H}$ is equivalent to summing over $H\leq Q\cap R$ with $H\geq Z$. If $H$ has rank $h$, then $\sur{H}=H/Z$ has rank $\sur{h}=h-1$, and $n+2\sur{h}-\sur{q}-\sur{r}=n+2h-q-r$. Moreover $\sur{R}/\sur{H}\cong R/H$, and $\sur{Q}/\sur{H}\cong Q/H$. Also
$$Z_{\sur{H}}=\{\sur{g}\in\sur{G}\mid [\sur{g},\sur{Q}{\cdot}\sur{R}]\leq\sur{H}\}/\sur{H}\cong \{g\in G\mid [g,Q{\cdot}R]\leq H\}/H=Z_H\mpoint$$
Thus $\Ext^{n}_{\comack_k(\sur{G})}(S_{\sur{Q}}^{\sur{G}},S_{\sur{R}}^{\sur{G}})$ is isomorphic to
\begin{equation}\label{part2}
\left(\dirsumc{h\in \N}{Z\leq H\leq Q\cap R}{{\rm rank}(H)=h}\!\!\!k{\rm St}(R/H)\otimes_k\Ext^{n+2h-q-r}_{\comack_k(Z_{H})}(S_\un^{Z_{H}},S_\un^{Z_{H}})\otimes_kk{\rm St}(Q/H)\right)_G\mpoint
\end{equation}
Now the direct sum of~\ref{part1} and~\ref{part2} is isomorphic to the expression in Theorem~\ref{factor Ext}, as was to be shown for Assertion~1.\par
It follows that the dimension of the space $\Ext^n_{\comack_k(G)}(S_Q^G,S_R^G)$ is equal to
$$\sumc{h\in\N}{\rule{0ex}{1.5ex}H\leq Q\cap R,\;{\rm mod.}G}{\rule{0ex}{1.5ex}{\rm rank}(H)=h}\hspace{-3ex}\dim_k \left(k{\rm St}(R/H)\otimes_k\Ext^{n+2h-q-r}_{\comack_k(Z_H)}(S_\un^{Z_H},S_\un^{Z_H})\otimes_kk{\rm St}(G/H)\right)_{N_G(H)}\;.$$
The image of the group $N_G(H)$ in the automorphism group of $R/H$ is isomorphic to $N_G(H)/Z_G(H,R)$, and $k{\rm St}(R/H)$ is a free $kN_G(H)/Z_G(H,R)$-module. Similarly the module $k{\rm St}(Q/H)$ is a free $kN_G(H)/Z_G(H,Q)$-module. It follows that the module $V=k{\rm St}(R/H)\otimes_k k{\rm St}(Q/H)$ is a free $kN_G(H)/Z$-module, where $Z=Z_G(H,R)\cap Z_G(H,Q)=Z_G(H,Q{\cdot}R)$.\par
Moreover since $Z_H=Z/H$, the group $Z$ acts trivially on 
$$W=\Ext^{n+2h-q-r}_{\comack_k(Z_H)}(S_\un^{Z_H},S_\un^{Z_H})\mvirg$$
so $W$ is a $kN_G(H)/Z$-module. Hence $V\otimes W$ is a free $kN_G(H)/Z$-module, and it follows that $\dim_k(V\otimes W)_{N_G(H)}=\frac{\dsp\dim_k(V\otimes W)}{\dsp|N_G(H)/Z|}$, which completes the proof of Assertion~2.\findemo
In the case where $G$ is elementary abelian, the following has been proved by Tambara~(\cite{tambara} Theorem~B)~:

\begin{mth}{Corollary} Let $k$ be a field of characteristic $p$, let $G$ be a finite $p$-group, and let $Q$ and $R$ be subgroups of $G$, such that $\nu_G(R,Q)<\infty$. Then $\Ext_{\comack_k(G)}^{\nu_G(R,Q)}(S_Q^G,S_R^G)$ is isomorphic to
$$\dirsum{g\in\mathcal{M}_{Q,R}^G}\Big(k{\rm St}\big(R/(R\cap {^gQ})\big)\otimes_kk{\rm St}\big({^gQ}/(R\cap {^gQ})\big)\Big)_{N_G(R,{^gQ})}\mvirg$$
where $\mathcal{M}_{Q,R}^G=\{g\in[\sur{\mathcal{S}}_{R,Q}^G]\mid |R{\cdot}{^gQ}:R\cap {^gQ}|=p^{\nu_G(R,Q)}\}$.\medskip\par
In particular $\nu_G(R,Q)=\min\{n\in\N\mid\Ext^{n}_{\comack_k(G)}(S_Q^G,S_R^G)\neq\zero\}$.
\end{mth}
\pf If $n<\nu_G(R,Q)$, then $\Ext^{n}_{\comack_k(G)}(S_Q^G,S_R^G)=\zero$, by Proposition~\ref{zero}. Suppose that $n=\nu_G(R,Q)$, and let $g\in [\sur{\mathcal{S}}_{R,Q}^G]$ such that, with the notation of Theorem~\ref{general Ext},
$$\Ext^n_{\comack_k(\widehat{N}_G(R,{^gQ}))}\big(S_{\widehat{^gQ}}^{\widehat{N}_G(R,{^gQ})},S_{\widehat{R}}^{\widehat{N}_G(R,{^gQ})}\big)\neq 0\mpoint$$
It follows that $p^n\geq p^{\nu_{\widehat{N}_G(R,{^gQ})}(\widehat{R},\widehat{^gQ})}=|R{\cdot}{^gQ}:R\cap {^gQ}|\geq p^{\nu_G(R,Q)}$. Hence $n=\nu_G(R,Q)=|R{\cdot}{^gQ}:R\cap {^gQ}|$. \par
Moreover, in the expression of $\Ext^n_{\comack_k(\widehat{N}_G(R,{^gQ}))}\big(S_{\widehat{^gQ}}^{\widehat{N}_G(R,{^gQ})},S_{\widehat{R}}^{\widehat{N}_G(R,{^gQ})}\big)$ given by Theorem~\ref{factor Ext}, the non zero terms correspond to subgroups $H$ of $\widehat{R}\cap \widehat{^gQ}$ such that $n+2h-q-r\geq 0$, where $h$, $q$, and $r$ are the ranks of $H$, $\widehat{^gQ}$ and~$\widehat{R}$ respectively. Thus 
$$|H|^2\geq |\widehat{^gQ}||\widehat{R}|/p^n=|\widehat{^gQ}||\widehat{R}||/|\widehat{R}{\cdot}\widehat{^gQ}|=|\widehat{R}\cap\widehat{^gQ}|^2\mpoint$$
It follows that $H=\widehat{R}\cap\widehat{^gQ}$, that $\widehat{R}/H\cong R/R\cap{^gQ}$, that $\widehat{^gQ}/H\cong {^gQ}/R\cap{^gQ}$, and that $n+2h-q-r=0$, so $\Ext^0_{\comack_k(Z_H)}(S_\un^{Z_H},S_\un^{Z_H})\cong k$. This completes the proof. \findemo
\begin{rem}{Remark} Let $G$ be an abelian $p$-group, and $Q$, $R$ be elementary abelian subgroups of $G$, of rank $q$ and $r$, respectively. In this case, by Assertion~1 of Theorem~\ref{general Ext}, the group $\Ext^n_{\comack_k(G)}(S_Q^G,S_R^G)$ is isomorphic to
$$\dirsumc{h\in\N}{\rule{0ex}{1.5ex}H\leq Q\cap R}{\rule{0ex}{1.5ex}{\rm rank}(H)=h}k{\rm St}(R/H)\otimes_k\Ext^{n+2h-q-r}_{\comack_k(G/H)}(S_\un^{G/H},S_\un^{G/H})\otimes_kk{\rm St}(Q/H)\mpoint$$
Fix $H\leq Q\cap R$, of rank $h$, and denote by $i$ the inclusion map 
$$k{\rm St}(R/H)\otimes_k\Ext^{n+2h-q-r}_{\comack_k(G/H)}(S_\un^{G/H},S_\un^{G/H})\otimes_kk{\rm St}(Q/H)\hookrightarrow \Ext^n_{\comack_k(G)}(S_Q^G,S_R^G).$$
It is easy to check that this map can be explicitly described as follows~: if $y\in k{\rm St}(R/H)\cong\Ext^{r-h}_{\comack_k(G)}(S_H^G,S_R^G)$, if $e\in \Ext^{n+2h-q-r}_{\comack_k(G/H)}(S_\un^{G/H},S_\un^{G/H})$, and if $x\in k{\rm St}(Q/H)\cong\Ext^{q-h}_{\comack_k(G)}(S_Q^G,S_H^G)$, then $i(y\otimes e\otimes x)$ is equal to the Yoneda composition $y\circ \Inf_{G/H}^Ge\circ x$.\par
In terms of morphisms in the derived category of $\comack_k(G)$, it means that any morphism from $S_Q^G$ to some translate of $S_R^G$ is a linear combination of morphisms which factor through a morphism of the form $\Inf_{G/H}^Ge$, where $H\leq Q\cap R$, and $e$ is a morphism from $S_\un^{G/H}$ to some of its translate.\par
In the case $p=2$ and $G$ is elementary abelian, this leads to an explicit description of all the Yoneda products of extensions of simple cohomological Mackey functors.
\end{rem}
%

\vspace{3ex}
\noindent Serge Bouc - CNRS - LAMFA\\
Universit\'e de Picardie Jules Verne\\
33, rue St Leu - 80039 - Amiens cedex 1\\
FRANCE\\
email~: {\tt serge.bouc@u-picardie.fr}
\end{document}